%
%

\font \aufont=  cmr10 at 12pt 
\font\titfont=  cmr10 at 18pt 
\font\footfont=cmr10 at 8  pt
\font\ftfont=cmr10   

 \magnification=1200 

%
%

\def\fmtversion{2.0}
\catcode`\@=11
\ifx\amstexloaded@\relax\catcode`\@=\active
 \endinput\else\let\amstexloaded@\relax\fi
\def\W@{\immediate\write\sixt@@n}
\def\CR@{\W@{}\W@{AmS-TeX - Version \fmtversion}\W@{}
\W@{COPYRIGHT 1985, 1990 - AMERICAN MATHEMATICAL SOCIETY}
\W@{Use of this macro package is not restricted provided}
\W@{each use is acknowledged upon publication.}\W@{}}
\CR@
\everyjob{\CR@}
\toksdef\toks@@=2
\long\def\rightappend@#1\to#2{\toks@{\\{#1}}\toks@@
 =\expandafter{#2}\xdef#2{\the\toks@@\the\toks@}\toks@{}\toks@@{}}
\def\alloclist@{}
\newif\ifalloc@
\def\showallocations{{\def\\{\immediate\write\m@ne}\alloclist@}\alloc@true}
\def\alloc@#1#2#3#4#5{\global\advance\count1#1by\@ne
 \ch@ck#1#4#2\allocationnumber=\count1#1
 \global#3#5=\allocationnumber
 \edef\next@{\string#5=\string#2\the\allocationnumber}%
 \expandafter\rightappend@\next@\to\alloclist@}
\newcount\count@@
\newcount\count@@@
\def\FN@{\futurelet\next}
\def\DN@{\def\next@}
\def\DNii@{\def\nextii@}
\def\RIfM@{\relax\ifmmode}
\def\RIfMIfI@{\relax\ifmmode\ifinner}
\def\setboxz@h{\setbox\z@\hbox}
\def\wdz@{\wd\z@}
\def\boxz@{\box\z@}
\def\setbox@ne{\setbox\@ne}
\def\wd@ne{\wd\@ne}
\def\iterate{\body\expandafter\iterate\else\fi}
\newlinechar=`\^^J
\def\err@#1{\errmessage{AmS-TeX error: #1}}
\newhelp\defaulthelp@{Sorry, I already gave what help I could...^^J
Maybe you should try asking a human?^^J
An error might have occurred before I noticed any problems.^^J
``If all else fails, read the instructions.''}
\def\Err@#1{\errhelp\defaulthelp@\errmessage{AmS-TeX error: #1}}
\def\eat@#1{}
\def\in@#1#2{\def\in@@##1#1##2##3\in@@{\ifx\in@##2\in@false\else\in@true\fi}%
 \in@@#2#1\in@\in@@}
\newif\ifin@
\def\space@.{\futurelet\space@\relax}
\space@. %
\newhelp\athelp@
{Only certain combinations beginning with @ make sense to me.^^J
Perhaps you wanted \string\@\space for a printed @?^^J
I've ignored the character or group after @.}
\def\futureletnextat@{\futurelet\next\at@}
{\catcode`\@=\active
\lccode`\Z=`\@ \lccode`\I=`\I \lowercase{%
\gdef@{\csname futureletnextatZ\endcsname}
\expandafter\gdef\csname atZ\endcsname                                      
 {\ifcat\noexpand\next a\def\next{\csname atZZ\endcsname}\else
 \ifcat\noexpand\next0\def\next{\csname atZZ\endcsname}\else
 \ifcat\noexpand\next\relax\def\next{\csname atZZZ\endcsname}\else
 \def\next{\csname atZZZZ\endcsname}\fi\fi\fi\next}
\expandafter\gdef\csname atZZ\endcsname#1{\expandafter                      
 \ifx\csname #1Zat\endcsname\relax\def\next
 {\errhelp\expandafter=\csname athelpZ\endcsname
 \csname errZ\endcsname{Invalid use of \string@}}\else
 \def\next{\csname #1Zat\endcsname}\fi\next}
\expandafter\gdef\csname atZZZ\endcsname#1{\expandafter                     
 \ifx\csname \string#1ZZat\endcsname\relax\def\next
 {\errhelp\expandafter=\csname athelpZ\endcsname
 \csname errZ\endcsname{Invalid use of \string@}}\else
 \def\next{\csname \string#1ZZat\endcsname}\fi\next}
\expandafter\gdef\csname atZZZZ\endcsname#1{\errhelp                        
 \expandafter=\csname athelpZ\endcsname
 \csname errZ\endcsname{Invalid use of \string@}}}}                         
\def\atdef@#1{\expandafter\def\csname #1@at\endcsname}
\def\atdef@@#1{\expandafter\def\csname \string#1@@at\endcsname}
\newhelp\defahelp@{If you typed \string\define\space cs instead of
\string\define\string\cs\space^^J
I've substituted an inaccessible control sequence so that your^^J
definition will be completed without mixing me up too badly.^^J
If you typed \string\define{\string\cs} the inaccessible control sequence^^J
was defined to be \string\cs, and the rest of your^^J
definition appears as input.}
\newhelp\defbhelp@{I've ignored your definition, because it might^^J
conflict with other uses that are important to me.}
\def\define{\FN@\define@}
\def\define@{\ifcat\noexpand\next\relax
 \expandafter\define@@\else\errhelp\defahelp@                               
 \err@{\string\define\space must be followed by a control
 sequence}\expandafter\def\expandafter\nextii@\fi}                          
\def\undefined@@@@@@@@@@{}
\def\preloaded@@@@@@@@@@{}
\def\next@@@@@@@@@@{}
\def\define@@#1{\ifx#1\relax\errhelp\defbhelp@                              
 \err@{\string#1\space is already defined}\DN@{\DNii@}\else
 \expandafter\ifx\csname\expandafter\eat@\string                            
 #1@@@@@@@@@@\endcsname\undefined@@@@@@@@@@\errhelp\defbhelp@
 \err@{\string#1\space can't be defined}\DN@{\DNii@}\else
 \expandafter\ifx\csname\expandafter\eat@\string#1\endcsname\relax          
 \global\let#1\undefined\DN@{\def#1}\else\errhelp\defbhelp@
 \err@{\string#1\space is already defined}\DN@{\DNii@}\fi
 \fi\fi\next@}
\let\redefine\def
\def\predefine#1#2{\let#1#2}
\def\undefine#1{\let#1\undefined}
\newdimen\captionwidth@
\captionwidth@\hsize
\advance\captionwidth@-1.5in
\def\pagewidth#1{\hsize#1\relax
 \captionwidth@\hsize\advance\captionwidth@-1.5in}
\def\pageheight#1{\vsize#1\relax}
\def\hcorrection#1{\advance\hoffset#1\relax}
\def\vcorrection#1{\advance\voffset#1\relax}

\let\graveaccent\`
\let\acuteaccent\'
\let\tildeaccent\~
\let\hataccent\^
\let\underscore\_
\let\B\=
\let\D\.
\let\ic@\/
\def\/{\unskip\ic@}
\def\textfonti{\the\textfont\@ne}
\def\t#1#2{{\edef\next@{\the\font}\textfonti\accent"7F \next@#1#2}}
\def~{\unskip\nobreak\ \ignorespaces}
\def\.{.\spacefactor\@m}
\atdef@;{\leavevmode\null;}
\atdef@:{\leavevmode\null:}
\atdef@?{\leavevmode\null?}
\def\@{\char64 }
\def\relaxnext@{\let\next\relax}
\atdef@-{\relaxnext@\leavevmode
 \DN@{\ifx\next-\DN@-{\FN@\nextii@}\else
  \DN@{\leavevmode\hbox{-}}\fi\next@}%
 \DNii@{\ifx\next-\DN@-{\leavevmode\hbox{---}}\else
  \DN@{\leavevmode\hbox{--}}\fi\next@}%
 \FN@\next@}
\def\srdr@{\kern.16667em}
\def\drsr@{\kern.02778em}
\def\sldl@{\kern.02778em}
\def\dlsl@{\kern.16667em}
\atdef@"{\unskip\relaxnext@
 \DN@{\ifx\next\space@\DN@. {\FN@\nextii@}\else
  \DN@.{\FN@\nextii@}\fi\next@.}%
 \DNii@{\ifx\next`\DN@`{\FN@\nextiii@}\else
  \ifx\next\lq\DN@\lq{\FN@\nextiii@}\else
  \DN@####1{\FN@\nextiv@}\fi\fi\next@}%
 \def\nextiii@{\ifx\next`\DN@`{\sldl@``}\else\ifx\next\lq
  \DN@\lq{\sldl@``}\else\DN@{\dlsl@`}\fi\fi\next@}%
 \def\nextiv@{\ifx\next'\DN@'{\srdr@''}\else
  \ifx\next\rq\DN@\rq{\srdr@''}\else\DN@{\drsr@'}\fi\fi\next@}%
 \FN@\next@}

\def\textfontii{\the\textfont\tw@}
\def\lbrace@{\delimiter"4266308 }
\def\rbrace@{\delimiter"5267309 }
\def\{{\RIfM@\lbrace@\else{\textfontii f}\spacefactor\@m\fi}
\def\}{\RIfM@\rbrace@\else
 \let\@sf\empty\ifhmode\edef\@sf{\spacefactor\the\spacefactor}\fi
 {\textfontii g}\@sf\relax\fi}
\let\lbrace\{
\let\rbrace\}
\def\AmSTeX{{\textfontii A}\kern-.1667em\lower.5ex\hbox
 {\textfontii M}\kern-.125em{\textfontii S}-\TeX}
\def\vmodeerr@#1{\Err@{\string#1\space not allowed between paragraphs}}
\def\mathmodeerr@#1{\Err@{\string#1\space not allowed in math mode}}
\def\linebreak{\RIfM@\mathmodeerr@\linebreak\else
 \ifhmode\unskip\unkern\break\else\vmodeerr@\linebreak\fi\fi}

\newskip\saveskip@
\def\allowlinebreak{\RIfM@\mathmodeerr@\allowlinebreak\else
 \ifhmode\saveskip@\lastskip\unskip
 \allowbreak\ifdim\saveskip@>\z@\hskip\saveskip@\fi
 \else\vmodeerr@\allowlinebreak\fi\fi}
\def\nolinebreak{\RIfM@\mathmodeerr@\nolinebreak\else
 \ifhmode\saveskip@\lastskip\unskip
 \nobreak\ifdim\saveskip@>\z@\hskip\saveskip@\fi
 \else\vmodeerr@\nolinebreak\fi\fi}
\def\newline{\relaxnext@
 \DN@{\RIfM@\expandafter\mathmodeerr@\expandafter\newline\else
  \ifhmode\ifx\next\par\else
  \expandafter\unskip\expandafter\null\expandafter\hfill\expandafter\break\fi
  \else
  \expandafter\vmodeerr@\expandafter\newline\fi\fi}%
 \FN@\next@}
\def\dmatherr@#1{\Err@{\string#1\space not allowed in display math mode}}
\def\nondmatherr@#1{\Err@{\string#1\space not allowed in non-display math
 mode}}
\def\onlydmatherr@#1{\Err@{\string#1\space allowed only in display math mode}}
\def\nonmatherr@#1{\Err@{\string#1\space allowed only in math mode}}
\def\mathbreak{\RIfMIfI@\break\else
 \dmatherr@\mathbreak\fi\else\nonmatherr@\mathbreak\fi}
\def\nomathbreak{\RIfMIfI@\nobreak\else
 \dmatherr@\nomathbreak\fi\else\nonmatherr@\nomathbreak\fi}
\def\allowmathbreak{\RIfMIfI@\allowbreak\else
 \dmatherr@\allowmathbreak\fi\else\nonmatherr@\allowmathbreak\fi}
\def\pagebreak{\RIfM@
 \ifinner\nondmatherr@\pagebreak\else\postdisplaypenalty-\@M\fi
 \else\ifvmode\removelastskip\break\else\vadjust{\break}\fi\fi}
\def\nopagebreak{\RIfM@
 \ifinner\nondmatherr@\nopagebreak\else\postdisplaypenalty\@M\fi
 \else\ifvmode\nobreak\else\vadjust{\nobreak}\fi\fi}
\def\nonvmodeerr@#1{\Err@{\string#1\space not allowed within a paragraph
 or in math}}
\def\vnonvmode@#1#2{\relaxnext@\DNii@{\ifx\next\par\DN@{#1}\else
 \DN@{#2}\fi\next@}%
 \ifvmode\DN@{#1}\else
 \DN@{\FN@\nextii@}\fi\next@}
\def\newpage{\vnonvmode@{\vfill\break}{\nonvmodeerr@\newpage}}
\def\smallpagebreak{\vnonvmode@\smallbreak{\nonvmodeerr@\smallpagebreak}}
\def\medpagebreak{\vnonvmode@\medbreak{\nonvmodeerr@\medpagebreak}}
\def\bigpagebreak{\vnonvmode@\bigbreak{\nonvmodeerr@\bigpagebreak}}
\def\NoBlackBoxes{\global\overfullrule\z@}
\def\BlackBoxes{\global\overfullrule5\p@}
\def\Invalid@#1{\def#1{\Err@{\Invalid@@\string#1}}}
\def\Invalid@@{Invalid use of }
\Invalid@\caption
\Invalid@\captionwidth
\newdimen\smallcaptionwidth@
\def\topspace{\mid@false\ins@}
\def\midspace{\mid@true\ins@}
\newif\ifmid@
\def\captionfont@{}
\def\ins@#1{\relaxnext@\allowbreak
 \smallcaptionwidth@\captionwidth@\gdef\thespace@{#1}%
 \DN@{\ifx\next\space@\DN@. {\FN@\nextii@}\else
  \DN@.{\FN@\nextii@}\fi\next@.}%
 \DNii@{\ifx\next\caption\DN@\caption{\FN@\nextiii@}%
  \else\let\next@\nextiv@\fi\next@}%
 \def\nextiv@{\vnonvmode@
  {\ifmid@\expandafter\midinsert\else\expandafter\topinsert\fi
   \vbox to\thespace@{}\endinsert}
  {\ifmid@\nonvmodeerr@\midspace\else\nonvmodeerr@\topspace\fi}}%
 \def\nextiii@{\ifx\next\captionwidth\expandafter\nextv@
  \else\expandafter\nextvi@\fi}%
 \def\nextv@\captionwidth##1##2{\smallcaptionwidth@##1\relax\nextvi@{##2}}%
 \def\nextvi@##1{\def\thecaption@{\captionfont@##1}%
  \DN@{\ifx\next\space@\DN@. {\FN@\nextvii@}\else
   \DN@.{\FN@\nextvii@}\fi\next@.}%
  \FN@\next@}%
 \def\nextvii@{\vnonvmode@
  {\ifmid@\expandafter\midinsert\else
  \expandafter\topinsert\fi\vbox to\thespace@{}\nobreak\smallskip
  \setboxz@h{\noindent\ignorespaces\thecaption@\unskip}%
  \ifdim\wdz@>\smallcaptionwidth@\centerline{\vbox{\hsize\smallcaptionwidth@
   \noindent\ignorespaces\thecaption@\unskip}}%
  \else\centerline{\boxz@}\fi\endinsert}
  {\ifmid@\nonvmodeerr@\midspace
  \else\nonvmodeerr@\topspace\fi}}%
 \FN@\next@}
\def\newcodes@{\catcode`\\=12 \catcode`\{=12 \catcode`\}=12 \catcode`\#=12
 \catcode`\%=12\relax}
\def\oldcodes@{\catcode`\\=0 \catcode`\{=1 \catcode`\}=2 \catcode`\#=6
 \catcode`\%=14\relax}
\def\comment{\newcodes@\endlinechar=10 \comment@}
{\lccode`\0=`\\
\lowercase{\gdef\comment@#1^^J{\comment@@#10endcomment\comment@@@}%
\gdef\comment@@#10endcomment{\FN@\comment@@@}%
\gdef\comment@@@#1\comment@@@{\ifx\next\comment@@@\let\next\comment@
 \else\def\next{\oldcodes@\endlinechar=`\^^M\relax}%
 \fi\next}}}
\def\pr@m@s{\ifx'\next\DN@##1{\prim@s}\else\let\next@\egroup\fi\next@}
\def\prime{{\null\prime@\null}}
\mathchardef\prime@="0230
\let\dsize\displaystyle
\let\tsize\textstyle
\let\ssize\scriptstyle

\def\,{\RIfM@\mskip\thinmuskip\relax\else\kern.16667em\fi}
\def\!{\RIfM@\mskip-\thinmuskip\relax\else\kern-.16667em\fi}
\let\thinspace\,
\let\negthinspace\!
\def\medspace{\RIfM@\mskip\medmuskip\relax\else\kern.222222em\fi}
\def\negmedspace{\RIfM@\mskip-\medmuskip\relax\else\kern-.222222em\fi}
\def\thickspace{\RIfM@\mskip\thickmuskip\relax\else\kern.27777em\fi}
\let\;\thickspace
\def\negthickspace{\RIfM@\mskip-\thickmuskip\relax\else
 \kern-.27777em\fi}
\atdef@,{\RIfM@\mskip.1\thinmuskip\else\leavevmode\null,\fi}
\atdef@!{\RIfM@\mskip-.1\thinmuskip\else\leavevmode\null!\fi}
\atdef@.{\RIfM@&&\else\leavevmode.\spacefactor3000 \fi}
\def\and{\DOTSB\;\mathchar"3026 \;}

\def\frac#1#2{{#1\over#2}}

\newdimen\ex@
\ex@.2326ex
\Invalid@\thickness
\def\thickfrac{\relaxnext@
 \DN@{\ifx\next\thickness\let\next@\nextii@\else
 \DN@{\nextii@\thickness1}\fi\next@}%
 \DNii@\thickness##1##2##3{{##2\above##1\ex@##3}}%
 \FN@\next@}

\def\thickfracwithdelims#1#2{\relaxnext@\def\ldelim@{#1}\def\rdelim@{#2}%
 \DN@{\ifx\next\thickness\let\next@\nextii@\else
 \DN@{\nextii@\thickness1}\fi\next@}%
 \DNii@\thickness##1##2##3{{##2\abovewithdelims
 \ldelim@\rdelim@##1\ex@##3}}%
 \FN@\next@}

\def\:{\nobreak\hskip.1111em\mathpunct{}\nonscript\mkern-\thinmuskip{:}\hskip
 .3333emplus.0555em\relax}
\def\snug{\unskip\kern-\mathsurround}
\def\topsmash{\top@true\bot@false\smash@}
\def\botsmash{\top@false\bot@true\smash@}
\newif\iftop@
\newif\ifbot@
\def\smash{\top@true\bot@true\smash@}
\def\smash@{\RIfM@\expandafter\mathpalette\expandafter\mathsm@sh\else
 \expandafter\makesm@sh\fi}
\def\finsm@sh{\iftop@\ht\z@\z@\fi\ifbot@\dp\z@\z@\fi\leavevmode\boxz@}
\def\LimitsOnSums{\global\let\slimits@\displaylimits}
\def\NoLimitsOnSums{\global\let\slimits@\nolimits}
\LimitsOnSums
\mathchardef\coprod@="1360       \def\coprod{\DOTSB\coprod@\slimits@}
\mathchardef\bigvee@="1357       \def\bigvee{\DOTSB\bigvee@\slimits@}
\mathchardef\bigwedge@="1356     \def\bigwedge{\DOTSB\bigwedge@\slimits@}
\mathchardef\biguplus@="1355     \def\biguplus{\DOTSB\biguplus@\slimits@}
\mathchardef\bigcap@="1354       \def\bigcap{\DOTSB\bigcap@\slimits@}
\mathchardef\bigcup@="1353       \def\bigcup{\DOTSB\bigcup@\slimits@}
\mathchardef\prod@="1351         \def\prod{\DOTSB\prod@\slimits@}
\mathchardef\sum@="1350          \def\sum{\DOTSB\sum@\slimits@}
\mathchardef\bigotimes@="134E    \def\bigotimes{\DOTSB\bigotimes@\slimits@}
\mathchardef\bigoplus@="134C     \def\bigoplus{\DOTSB\bigoplus@\slimits@}
\mathchardef\bigodot@="134A      \def\bigodot{\DOTSB\bigodot@\slimits@}
\mathchardef\bigsqcup@="1346     \def\bigsqcup{\DOTSB\bigsqcup@\slimits@}
\def\LimitsOnInts{\global\let\ilimits@\displaylimits}
\def\NoLimitsOnInts{\global\let\ilimits@\nolimits}
\NoLimitsOnInts
\def\int{\DOTSI\intop\ilimits@}
\def\oint{\DOTSI\ointop\ilimits@}
\def\intic@{\mathchoice{\hskip.5em}{\hskip.4em}{\hskip.4em}{\hskip.4em}}
\def\negintic@{\mathchoice
 {\hskip-.5em}{\hskip-.4em}{\hskip-.4em}{\hskip-.4em}}
\def\intkern@{\mathchoice{\!\!\!}{\!\!}{\!\!}{\!\!}}
\def\intdots@{\mathchoice{\plaincdots@}
 {{\cdotp}\mkern1.5mu{\cdotp}\mkern1.5mu{\cdotp}}
 {{\cdotp}\mkern1mu{\cdotp}\mkern1mu{\cdotp}}
 {{\cdotp}\mkern1mu{\cdotp}\mkern1mu{\cdotp}}}
\newcount\intno@
\def\iint{\DOTSI\intno@\tw@\FN@\ints@}
\def\iiint{\DOTSI\intno@\thr@@\FN@\ints@}
\def\iiiint{\DOTSI\intno@4 \FN@\ints@}
\def\idotsint{\DOTSI\intno@\z@\FN@\ints@}
\def\ints@{\findlimits@\ints@@}
\newif\iflimtoken@
\newif\iflimits@
\def\findlimits@{\limtoken@true\ifx\next\limits\limits@true
 \else\ifx\next\nolimits\limits@false\else
 \limtoken@false\ifx\ilimits@\nolimits\limits@false\else
 \ifinner\limits@false\else\limits@true\fi\fi\fi\fi}
\def\multint@{\int\ifnum\intno@=\z@\intdots@                                
 \else\intkern@\fi                                                          
 \ifnum\intno@>\tw@\int\intkern@\fi                                         
 \ifnum\intno@>\thr@@\int\intkern@\fi                                       
 \int}                                                                      
\def\multintlimits@{\intop\ifnum\intno@=\z@\intdots@\else\intkern@\fi
 \ifnum\intno@>\tw@\intop\intkern@\fi
 \ifnum\intno@>\thr@@\intop\intkern@\fi\intop}
\def\ints@@{\iflimtoken@                                                    
 \def\ints@@@{\iflimits@\negintic@\mathop{\intic@\multintlimits@}\limits    
  \else\multint@\nolimits\fi                                                
  \eat@}                                                                    
 \else                                                                      
 \def\ints@@@{\iflimits@\negintic@
  \mathop{\intic@\multintlimits@}\limits\else
  \multint@\nolimits\fi}\fi\ints@@@}
\def\LimitsOnNames{\global\let\nlimits@\displaylimits}
\def\NoLimitsOnNames{\global\let\nlimits@\nolimits@}
\LimitsOnNames
\def\nolimits@{\relaxnext@
 \DN@{\ifx\next\limits\DN@\limits{\nolimits}\else
  \let\next@\nolimits\fi\next@}%
 \FN@\next@}
\def\newmcodes@{\mathcode`\'="0027 \mathcode`\*="002A \mathcode`\.="613A
 \mathcode`\-="002D \mathcode`\/="002F \mathcode`\:="603A }
\def\operatorname#1{\mathop{\newmcodes@\kern\z@\fam\z@#1}\nolimits@}
\def\operatornamewithlimits#1{\mathop{\newmcodes@\kern\z@\fam\z@#1}\nlimits@}
\def\qopname@#1{\mathop{\fam\z@#1}\nolimits@}
\def\qopnamewl@#1{\mathop{\fam\z@#1}\nlimits@}
\def\arccos{\qopname@{arccos}}
\def\arcsin{\qopname@{arcsin}}
\def\arctan{\qopname@{arctan}}
\def\arg{\qopname@{arg}}
\def\cos{\qopname@{cos}}
\def\cosh{\qopname@{cosh}}
\def\cot{\qopname@{cot}}
\def\coth{\qopname@{coth}}
\def\csc{\qopname@{csc}}
\def\deg{\qopname@{deg}}
\def\det{\qopnamewl@{det}}
\def\dim{\qopname@{dim}}
\def\exp{\qopname@{exp}}
\def\gcd{\qopnamewl@{gcd}}
\def\hom{\qopname@{hom}}
\def\inf{\qopnamewl@{inf}}
\def\injlim{\qopnamewl@{inj\,lim}}
\def\ker{\qopname@{ker}}
\def\lg{\qopname@{lg}}
\def\lim{\qopnamewl@{lim}}
\def\liminf{\qopnamewl@{lim\,inf}}
\def\limsup{\qopnamewl@{lim\,sup}}
\def\ln{\qopname@{ln}}
\def\log{\qopname@{log}}
\def\max{\qopnamewl@{max}}
\def\min{\qopnamewl@{min}}
\def\Pr{\qopnamewl@{Pr}}
\def\projlim{\qopnamewl@{proj\,lim}}
\def\sec{\qopname@{sec}}
\def\sin{\qopname@{sin}}
\def\sinh{\qopname@{sinh}}
\def\sup{\qopnamewl@{sup}}
\def\tan{\qopname@{tan}}
\def\tanh{\qopname@{tanh}}
\def\varinjlim{\mathop{\vtop{\ialign{##\crcr
 \hfil\rm lim\hfil\crcr\noalign{\nointerlineskip}\rightarrowfill\crcr
 \noalign{\nointerlineskip\kern-\ex@}\crcr}}}}
\def\varprojlim{\mathop{\vtop{\ialign{##\crcr
 \hfil\rm lim\hfil\crcr\noalign{\nointerlineskip}\leftarrowfill\crcr
 \noalign{\nointerlineskip\kern-\ex@}\crcr}}}}
\def\varliminf{\mathop{\underline{\vrule height\z@ depth.2exwidth\z@
 \hbox{\rm lim}}}}

\newdimen\buffer@
\buffer@\fontdimen13 \tenex
\newdimen\buffer
\buffer\buffer@

\def\ResetBuffer{\fontdimen13 \tenex\buffer@\global\buffer\buffer@}
\def\shave#1{\mathop{\hbox{$\m@th\fontdimen13 \tenex\z@                     
 \displaystyle{#1}$}}\fontdimen13 \tenex\buffer}

\Invalid@\\
\def\Let@{\relax\iffalse{\fi\let\\=\cr\iffalse}\fi}
\Invalid@\vspace
\def\vspace@{\def\vspace##1{\crcr\noalign{\vskip##1\relax}}}
\def\multilimits@{\bgroup\vspace@\Let@
 \baselineskip\fontdimen10 \scriptfont\tw@
 \advance\baselineskip\fontdimen12 \scriptfont\tw@
 \lineskip\thr@@\fontdimen8 \scriptfont\thr@@
 \lineskiplimit\lineskip
 \vbox\bgroup\ialign\bgroup\hfil$\m@th\scriptstyle{##}$\hfil\crcr}
\def\Sb{_\multilimits@}
\def\endSb{\crcr\egroup\egroup\egroup}
\def\Sp{^\multilimits@}

\def\spreadlines#1{\RIfMIfI@\onlydmatherr@\spreadlines\else
 \openup#1\relax\fi\else\onlydmatherr@\spreadlines\fi}
\def\Mathstrut@{\copy\Mathstrutbox@}
\newbox\Mathstrutbox@
\setbox\Mathstrutbox@\null
\setboxz@h{$\m@th($}
\ht\Mathstrutbox@\ht\z@
\dp\Mathstrutbox@\dp\z@
\newdimen\spreadmlines@
\def\spreadmatrixlines#1{\RIfMIfI@
 \onlydmatherr@\spreadmatrixlines\else
 \spreadmlines@#1\relax\fi\else\onlydmatherr@\spreadmatrixlines\fi}
\def\matrix{\null\,\vcenter\bgroup\Let@\vspace@
 \normalbaselines\openup\spreadmlines@\ialign
 \bgroup\hfil$\m@th##$\hfil&&\quad\hfil$\m@th##$\hfil\crcr
 \Mathstrut@\crcr\noalign{\kern-\baselineskip}}
\def\endmatrix{\crcr\Mathstrut@\crcr\noalign{\kern-\baselineskip}\egroup
 \egroup\,}
\def\format{\crcr\egroup\iffalse{\fi\ifnum`}=0 \fi\format@}
\newtoks\hashtoks@
\hashtoks@{#}
\def\format@#1\\{\def\preamble@{#1}%
 \def\l{$\m@th\the\hashtoks@$\hfil}%
 \def\c{\hfil$\m@th\the\hashtoks@$\hfil}%
 \def\r{\hfil$\m@th\the\hashtoks@$}%
 \edef\Preamble@{\preamble@}\ifnum`{=0 \fi\iffalse}\fi
 \ialign\bgroup\span\Preamble@\crcr}
\def\smallmatrix{\null\,\vcenter\bgroup\vspace@\Let@
 \baselineskip9\ex@\lineskip\ex@
 \ialign\bgroup\hfil$\m@th\scriptstyle{##}$\hfil&&\thickspace\hfil
 $\m@th\scriptstyle{##}$\hfil\crcr}
\def\endsmallmatrix{\crcr\egroup\egroup\,}

\newmuskip\dotsspace@
\dotsspace@1.5mu
\def\strip@#1 {#1}
\def\spacehdots#1\for#2{\multispan{#2}\xleaders
 \hbox{$\m@th\mkern\strip@#1 \dotsspace@.\mkern\strip@#1 \dotsspace@$}\hfill}
\def\hdotsfor#1{\spacehdots\@ne\for{#1}}
\def\multispan@#1{\omit\mscount#1\unskip\loop\ifnum\mscount>\@ne\sp@n\repeat}
\def\spaceinnerhdots#1\for#2\after#3{\multispan@{\strip@#2 }#3\xleaders
 \hbox{$\m@th\mkern\strip@#1 \dotsspace@.\mkern\strip@#1 \dotsspace@$}\hfill}
\def\innerhdotsfor#1\after#2{\spaceinnerhdots\@ne\for#1\after{#2}}
\def\cases{\bgroup\spreadmlines@\jot\left\{\,\matrix\format\l&\quad\l\\}
\def\endcases{\endmatrix\right.\egroup}
\newif\ifinany@
\newif\ifinalign@
\newif\ifingather@
\def\strut@{\copy\strutbox@}
\newbox\strutbox@
\setbox\strutbox@\hbox{\vrule height8\p@ depth3\p@ width\z@}
\def\topaligned{\null\,\vtop\aligned@}
\def\botaligned{\null\,\vbox\aligned@}
\def\aligned{\null\,\vcenter\aligned@}
\def\aligned@{\bgroup\vspace@\Let@
 \ifinany@\else\openup\jot\fi\ialign
 \bgroup\hfil\strut@$\m@th\displaystyle{##}$&
 $\m@th\displaystyle{{}##}$\hfil\crcr}
\def\endaligned{\crcr\egroup\egroup}

\def\alignedat#1{\null\,\vcenter\bgroup\doat@{#1}\vspace@\Let@
 \ifinany@\else\openup\jot\fi\ialign\bgroup\span\preamble@@\crcr}
\newcount\atcount@
\def\doat@#1{\toks@{\hfil\strut@$\m@th
 \displaystyle{\the\hashtoks@}$&$\m@th\displaystyle
 {{}\the\hashtoks@}$\hfil}
 \atcount@#1\relax\advance\atcount@\m@ne                                    
 \loop\ifnum\atcount@>\z@\toks@=\expandafter{\the\toks@&\hfil$\m@th
 \displaystyle{\the\hashtoks@}$&$\m@th
 \displaystyle{{}\the\hashtoks@}$\hfil}\advance
  \atcount@\m@ne\repeat                                                     
 \xdef\preamble@{\the\toks@}\xdef\preamble@@{\preamble@}}

\def\gathered{\null\,\vcenter\bgroup\vspace@\Let@
 \ifinany@\else\openup\jot\fi\ialign
 \bgroup\hfil\strut@$\m@th\displaystyle{##}$\hfil\crcr}
\def\endgathered{\crcr\egroup\egroup}
\newif\iftagsleft@
\def\TagsOnLeft{\global\tagsleft@true}
\def\TagsOnRight{\global\tagsleft@false}
\TagsOnLeft
\newif\ifmathtags@
\def\TagsAsMath{\global\mathtags@true}
\def\TagsAsText{\global\mathtags@false}
\TagsAsText
\def\tagform@#1{\hbox{\rm(\ignorespaces#1\unskip)}}
\def\thetag{\leavevmode\tagform@}
\def\tag#1$${\iftagsleft@\leqno\else\eqno\fi                                
 \maketag@#1\maketag@                                                       
 $$}                                                                        
\def\maketag@{\FN@\maketag@@}
\def\maketag@@{\ifx\next"\expandafter\maketag@@@\else\expandafter\maketag@@@@
 \fi}
\def\maketag@@@"#1"#2\maketag@{\hbox{\rm#1}}                                
\def\maketag@@@@#1\maketag@{\ifmathtags@\tagform@{$\m@th#1$}\else
 \tagform@{#1}\fi}
\interdisplaylinepenalty\@M
\def\allowdisplaybreaks{\RIfMIfI@
 \onlydmatherr@\allowdisplaybreaks\else
 \interdisplaylinepenalty\z@\fi\else\onlydmatherr@\allowdisplaybreaks\fi}
\Invalid@\allowdisplaybreak
\Invalid@\displaybreak
\Invalid@\intertext
\def\allowdisplaybreak@{\def\allowdisplaybreak{\crcr\noalign{\allowbreak}}}
\def\displaybreak@{\def\displaybreak{\crcr\noalign{\break}}}
\def\intertext@{\def\intertext##1{\crcr\noalign{\vskip\belowdisplayskip
 \vbox{\normalbaselines\noindent##1}\vskip\abovedisplayskip}}}
\newskip\centering@
\centering@\z@ plus\@m\p@
\def\align{\relax\ifingather@\DN@{\csname align (in
  \string\gather)\endcsname}\else
 \ifmmode\ifinner\DN@{\onlydmatherr@\align}\else
  \let\next@\align@\fi
 \else\DN@{\onlydmatherr@\align}\fi\fi\next@}
\newhelp\andhelp@
{An extra & here is so disastrous that you should probably exit^^J
and fix things up.}
\newif\iftag@
\newcount\and@
\def\align@{\inalign@true\inany@true
 \vspace@\allowdisplaybreak@\displaybreak@\intertext@
 \def\tag{\global\tag@true\ifnum\and@=\z@\DN@{&&}\else
          \DN@{&}\fi\next@}%
 \iftagsleft@\DN@{\csname align \endcsname}\else
  \DN@{\csname align \space\endcsname}\fi\next@}
\def\Tag@{\iftag@\else\errhelp\andhelp@\err@{Extra & on this line}\fi}
\newdimen\lwidth@
\newdimen\rwidth@
\newdimen\maxlwidth@
\newdimen\maxrwidth@
\newdimen\totwidth@
\def\measure@#1\endalign{\lwidth@\z@\rwidth@\z@\maxlwidth@\z@\maxrwidth@\z@
 \global\and@\z@                                                            
 \setbox@ne\vbox                                                            
  {\everycr{\noalign{\global\tag@false\global\and@\z@}}\Let@                
  \halign{\setboxz@h{$\m@th\displaystyle{\@lign##}$}
   \global\lwidth@\wdz@                                                     
   \ifdim\lwidth@>\maxlwidth@\global\maxlwidth@\lwidth@\fi                  
   \global\advance\and@\@ne                                                 
   &\setboxz@h{$\m@th\displaystyle{{}\@lign##}$}\global\rwidth@\wdz@        
   \ifdim\rwidth@>\maxrwidth@\global\maxrwidth@\rwidth@\fi                  
   \global\advance\and@\@ne                                                
   &\Tag@
   \eat@{##}\crcr#1\crcr}}
 \totwidth@\maxlwidth@\advance\totwidth@\maxrwidth@}                       
\def\displ@y@{\global\dt@ptrue\openup\jot
 \everycr{\noalign{\global\tag@false\global\and@\z@\ifdt@p\global\dt@pfalse
 \vskip-\lineskiplimit\vskip\normallineskiplimit\else
 \penalty\interdisplaylinepenalty\fi}}}
\def\black@#1{\noalign{\ifdim#1>\displaywidth
 \dimen@\prevdepth\nointerlineskip                                          
 \vskip-\ht\strutbox@\vskip-\dp\strutbox@                                   
 \vbox{\noindent\hbox to#1{\strut@\hfill}}
 \prevdepth\dimen@                                                          
 \fi}}
\expandafter\def\csname align \space\endcsname#1\endalign
 {\measure@#1\endalign\global\and@\z@                                       
 \ifingather@\everycr{\noalign{\global\and@\z@}}\else\displ@y@\fi           
 \Let@\tabskip\centering@                                                   
 \halign to\displaywidth
  {\hfil\strut@\setboxz@h{$\m@th\displaystyle{\@lign##}$}
  \global\lwidth@\wdz@\boxz@\global\advance\and@\@ne                        
  \tabskip\z@skip                                                           
  &\setboxz@h{$\m@th\displaystyle{{}\@lign##}$}
  \global\rwidth@\wdz@\boxz@\hfill\global\advance\and@\@ne                  
  \tabskip\centering@                                                       
  &\setboxz@h{\@lign\strut@\maketag@##\maketag@}
  \dimen@\displaywidth\advance\dimen@-\totwidth@
  \divide\dimen@\tw@\advance\dimen@\maxrwidth@\advance\dimen@-\rwidth@     
  \ifdim\dimen@<\tw@\wdz@\llap{\vtop{\normalbaselines\null\boxz@}}
  \else\llap{\boxz@}\fi                                                    
  \tabskip\z@skip                                                          
  \crcr#1\crcr                                                             
  \black@\totwidth@}}                                                      
\newdimen\lineht@
\expandafter\def\csname align \endcsname#1\endalign{\measure@#1\endalign
 \global\and@\z@
 \ifdim\totwidth@>\displaywidth\let\displaywidth@\totwidth@\else
  \let\displaywidth@\displaywidth\fi                                        
 \ifingather@\everycr{\noalign{\global\and@\z@}}\else\displ@y@\fi
 \Let@\tabskip\centering@\halign to\displaywidth
  {\hfil\strut@\setboxz@h{$\m@th\displaystyle{\@lign##}$}%
  \global\lwidth@\wdz@\global\lineht@\ht\z@                                 
  \boxz@\global\advance\and@\@ne
  \tabskip\z@skip&\setboxz@h{$\m@th\displaystyle{{}\@lign##}$}%
  \global\rwidth@\wdz@\ifdim\ht\z@>\lineht@\global\lineht@\ht\z@\fi         
  \boxz@\hfil\global\advance\and@\@ne
  \tabskip\centering@&\kern-\displaywidth@                                  
  \setboxz@h{\@lign\strut@\maketag@##\maketag@}%
  \dimen@\displaywidth\advance\dimen@-\totwidth@
  \divide\dimen@\tw@\advance\dimen@\maxlwidth@\advance\dimen@-\lwidth@
  \ifdim\dimen@<\tw@\wdz@
   \rlap{\vbox{\normalbaselines\boxz@\vbox to\lineht@{}}}\else
   \rlap{\boxz@}\fi
  \tabskip\displaywidth@\crcr#1\crcr\black@\totwidth@}}
\expandafter\def\csname align (in \string\gather)\endcsname
  #1\endalign{\vcenter{\align@#1\endalign}}
\Invalid@\endalign
\newif\ifxat@
\def\alignat{\RIfMIfI@\DN@{\onlydmatherr@\alignat}\else
 \DN@{\csname alignat \endcsname}\fi\else
 \DN@{\onlydmatherr@\alignat}\fi\next@}
\newif\ifmeasuring@
\newbox\savealignat@
\expandafter\def\csname alignat \endcsname#1#2\endalignat                   
 {\inany@true\xat@false
 \def\tag{\global\tag@true\count@#1\relax\multiply\count@\tw@
  \xdef\tag@{}\loop\ifnum\count@>\and@\xdef\tag@{&\tag@}\advance\count@\m@ne
  \repeat\tag@}%
 \vspace@\allowdisplaybreak@\displaybreak@\intertext@
 \displ@y@\measuring@true                                                   
 \setbox\savealignat@\hbox{$\m@th\displaystyle\Let@
  \attag@{#1}
  \vbox{\halign{\span\preamble@@\crcr#2\crcr}}$}%
 \measuring@false                                                           
 \Let@\attag@{#1}
 \tabskip\centering@\halign to\displaywidth
  {\span\preamble@@\crcr#2\crcr                                             
  \black@{\wd\savealignat@}}}                                               
\Invalid@\endalignat
\def\xalignat{\RIfMIfI@
 \DN@{\onlydmatherr@\xalignat}\else
 \DN@{\csname xalignat \endcsname}\fi\else
 \DN@{\onlydmatherr@\xalignat}\fi\next@}
\expandafter\def\csname xalignat \endcsname#1#2\endxalignat
 {\inany@true\xat@true
 \def\tag{\global\tag@true\def\tag@{}\count@#1\relax\multiply\count@\tw@
  \loop\ifnum\count@>\and@\xdef\tag@{&\tag@}\advance\count@\m@ne\repeat\tag@}%
 \vspace@\allowdisplaybreak@\displaybreak@\intertext@
 \displ@y@\measuring@true\setbox\savealignat@\hbox{$\m@th\displaystyle\Let@
 \attag@{#1}\vbox{\halign{\span\preamble@@\crcr#2\crcr}}$}%
 \measuring@false\Let@
 \attag@{#1}\tabskip\centering@\halign to\displaywidth
 {\span\preamble@@\crcr#2\crcr\black@{\wd\savealignat@}}}
\def\attag@#1{\let\Maketag@\maketag@\let\TAG@\Tag@                          
 \let\Tag@=0\let\maketag@=0
 \ifmeasuring@\def\llap@##1{\setboxz@h{##1}\hbox to\tw@\wdz@{}}%
  \def\rlap@##1{\setboxz@h{##1}\hbox to\tw@\wdz@{}}\else
  \let\llap@\llap\let\rlap@\rlap\fi                                         
 \toks@{\hfil\strut@$\m@th\displaystyle{\@lign\the\hashtoks@}$\tabskip\z@skip
  \global\advance\and@\@ne&$\m@th\displaystyle{{}\@lign\the\hashtoks@}$\hfil
  \ifxat@\tabskip\centering@\fi\global\advance\and@\@ne}
 \iftagsleft@
  \toks@@{\tabskip\centering@&\Tag@\kern-\displaywidth
   \rlap@{\@lign\maketag@\the\hashtoks@\maketag@}%
   \global\advance\and@\@ne\tabskip\displaywidth}\else
  \toks@@{\tabskip\centering@&\Tag@\llap@{\@lign\maketag@
   \the\hashtoks@\maketag@}\global\advance\and@\@ne\tabskip\z@skip}\fi      
 \atcount@#1\relax\advance\atcount@\m@ne
 \loop\ifnum\atcount@>\z@
 \toks@=\expandafter{\the\toks@&\hfil$\m@th\displaystyle{\@lign
  \the\hashtoks@}$\global\advance\and@\@ne
  \tabskip\z@skip&$\m@th\displaystyle{{}\@lign\the\hashtoks@}$\hfil\ifxat@
  \tabskip\centering@\fi\global\advance\and@\@ne}\advance\atcount@\m@ne
 \repeat                                                                    
 \xdef\preamble@{\the\toks@\the\toks@@}
 \xdef\preamble@@{\preamble@}
 \let\maketag@\Maketag@\let\Tag@\TAG@}                                      
\Invalid@\endxalignat
\def\xxalignat{\RIfMIfI@
 \DN@{\onlydmatherr@\xxalignat}\else\DN@{\csname xxalignat
  \endcsname}\fi\else
 \DN@{\onlydmatherr@\xxalignat}\fi\next@}
\expandafter\def\csname xxalignat \endcsname#1#2\endxxalignat{\inany@true
 \vspace@\allowdisplaybreak@\displaybreak@\intertext@
 \displ@y\setbox\savealignat@\hbox{$\m@th\displaystyle\Let@
 \xxattag@{#1}\vbox{\halign{\span\preamble@@\crcr#2\crcr}}$}%
 \Let@\xxattag@{#1}\tabskip\z@skip\halign to\displaywidth
 {\span\preamble@@\crcr#2\crcr\black@{\wd\savealignat@}}}
\def\xxattag@#1{\toks@{\tabskip\z@skip\hfil\strut@
 $\m@th\displaystyle{\the\hashtoks@}$&%
 $\m@th\displaystyle{{}\the\hashtoks@}$\hfil\tabskip\centering@&}%
 \atcount@#1\relax\advance\atcount@\m@ne\loop\ifnum\atcount@>\z@
 \toks@=\expandafter{\the\toks@&\hfil$\m@th\displaystyle{\the\hashtoks@}$%
  \tabskip\z@skip&$\m@th\displaystyle{{}\the\hashtoks@}$\hfil
  \tabskip\centering@}\advance\atcount@\m@ne\repeat
 \xdef\preamble@{\the\toks@\tabskip\z@skip}\xdef\preamble@@{\preamble@}}
\Invalid@\endxxalignat
\newdimen\gwidth@
\newdimen\gmaxwidth@
\def\gmeasure@#1\endgather{\gwidth@\z@\gmaxwidth@\z@\setbox@ne\vbox{\Let@
 \halign{\setboxz@h{$\m@th\displaystyle{##}$}\global\gwidth@\wdz@
 \ifdim\gwidth@>\gmaxwidth@\global\gmaxwidth@\gwidth@\fi
 &\eat@{##}\crcr#1\crcr}}}
\def\gather{\RIfMIfI@\DN@{\onlydmatherr@\gather}\else
 \ingather@true\inany@true\def\tag{&}%
 \vspace@\allowdisplaybreak@\displaybreak@\intertext@
 \displ@y\Let@
 \iftagsleft@\DN@{\csname gather \endcsname}\else
  \DN@{\csname gather \space\endcsname}\fi\fi
 \else\DN@{\onlydmatherr@\gather}\fi\next@}
\expandafter\def\csname gather \space\endcsname#1\endgather
 {\gmeasure@#1\endgather\tabskip\centering@
 \halign to\displaywidth{\hfil\strut@\setboxz@h{$\m@th\displaystyle{##}$}%
 \global\gwidth@\wdz@\boxz@\hfil&
 \setboxz@h{\strut@{\maketag@##\maketag@}}%
 \dimen@\displaywidth\advance\dimen@-\gwidth@
 \ifdim\dimen@>\tw@\wdz@\llap{\boxz@}\else
 \llap{\vtop{\normalbaselines\null\boxz@}}\fi
 \tabskip\z@skip\crcr#1\crcr\black@\gmaxwidth@}}
\newdimen\glineht@
\expandafter\def\csname gather \endcsname#1\endgather{\gmeasure@#1\endgather
 \ifdim\gmaxwidth@>\displaywidth\let\gdisplaywidth@\gmaxwidth@\else
 \let\gdisplaywidth@\displaywidth\fi\tabskip\centering@\halign to\displaywidth
 {\hfil\strut@\setboxz@h{$\m@th\displaystyle{##}$}%
 \global\gwidth@\wdz@\global\glineht@\ht\z@\boxz@\hfil&\kern-\gdisplaywidth@
 \setboxz@h{\strut@{\maketag@##\maketag@}}%
 \dimen@\displaywidth\advance\dimen@-\gwidth@
 \ifdim\dimen@>\tw@\wdz@\rlap{\boxz@}\else
 \rlap{\vbox{\normalbaselines\boxz@\vbox to\glineht@{}}}\fi
 \tabskip\gdisplaywidth@\crcr#1\crcr\black@\gmaxwidth@}}
\newif\ifctagsplit@
\def\CenteredTagsOnSplits{\global\ctagsplit@true}
\def\TopOrBottomTagsOnSplits{\global\ctagsplit@false}
\TopOrBottomTagsOnSplits
\def\split{\relax\ifinany@\let\next@\insplit@\else
 \ifmmode\ifinner\def\next@{\onlydmatherr@\split}\else
 \let\next@\outsplit@\fi\else
 \def\next@{\onlydmatherr@\split}\fi\fi\next@}
\def\insplit@{\global\setbox\z@\vbox\bgroup\vspace@\Let@\ialign\bgroup
 \hfil\strut@$\m@th\displaystyle{##}$&$\m@th\displaystyle{{}##}$\hfill\crcr}
\def\endsplit{\crcr\egroup\egroup\iftagsleft@\expandafter\lendsplit@\else
 \expandafter\rendsplit@\fi}
\def\rendsplit@{\global\setbox9 \vbox
 {\unvcopy\z@\global\setbox8 \lastbox\unskip}
 \setbox@ne\hbox{\unhcopy8 \unskip\global\setbox\tw@\lastbox
 \unskip\global\setbox\thr@@\lastbox}
 \global\setbox7 \hbox{\unhbox\tw@\unskip}
 \ifinalign@\ifctagsplit@                                                   
  \gdef\split@{\hbox to\wd\thr@@{}&
   \vcenter{\vbox{\moveleft\wd\thr@@\boxz@}}}
 \else\gdef\split@{&\vbox{\moveleft\wd\thr@@\box9}\crcr
  \box\thr@@&\box7}\fi                                                      
 \else                                                                      
  \ifctagsplit@\gdef\split@{\vcenter{\boxz@}}\else
  \gdef\split@{\box9\crcr\hbox{\box\thr@@\box7}}\fi
 \fi
 \split@}                                                                   
\def\lendsplit@{\global\setbox9\vtop{\unvcopy\z@}
 \setbox@ne\vbox{\unvcopy\z@\global\setbox8\lastbox}
 \setbox@ne\hbox{\unhcopy8\unskip\setbox\tw@\lastbox
  \unskip\global\setbox\thr@@\lastbox}
 \ifinalign@\ifctagsplit@                                                   
  \gdef\split@{\hbox to\wd\thr@@{}&
  \vcenter{\vbox{\moveleft\wd\thr@@\box9}}}
  \else                                                                     
  \gdef\split@{\hbox to\wd\thr@@{}&\vbox{\moveleft\wd\thr@@\box9}}\fi
 \else
  \ifctagsplit@\gdef\split@{\vcenter{\box9}}\else
  \gdef\split@{\box9}\fi
 \fi\split@}
\def\outsplit@#1$${\align\insplit@#1\endalign$$}
\newdimen\multlinegap@
\multlinegap@1em
\newdimen\multlinetaggap@
\multlinetaggap@1em
\def\MultlineGap#1{\global\multlinegap@#1\relax}
\def\multlinegap#1{\RIfMIfI@\onlydmatherr@\multlinegap\else
 \multlinegap@#1\relax\fi\else\onlydmatherr@\multlinegap\fi}
\def\nomultlinegap{\multlinegap{\z@}}
\def\multline{\RIfMIfI@
 \DN@{\onlydmatherr@\multline}\else
 \DN@{\multline@}\fi\else
 \DN@{\onlydmatherr@\multline}\fi\next@}
\newif\iftagin@
\def\tagin@#1{\tagin@false\in@\tag{#1}\ifin@\tagin@true\fi}
\def\multline@#1$${\inany@true\vspace@\allowdisplaybreak@\displaybreak@
 \tagin@{#1}\iftagsleft@\DN@{\multline@l#1$$}\else
 \DN@{\multline@r#1$$}\fi\next@}
\newdimen\mwidth@
\def\rmmeasure@#1\endmultline{%
 \def\shoveleft##1{##1}\def\shoveright##1{##1}
 \setbox@ne\vbox{\Let@\halign{\setboxz@h
  {$\m@th\@lign\displaystyle{}##$}\global\mwidth@\wdz@
  \crcr#1\crcr}}}
\newdimen\mlineht@
\newif\ifzerocr@
\newif\ifonecr@
\def\lmmeasure@#1\endmultline{\global\zerocr@true\global\onecr@false
 \everycr{\noalign{\ifonecr@\global\onecr@false\fi
  \ifzerocr@\global\zerocr@false\global\onecr@true\fi}}
  \def\shoveleft##1{##1}\def\shoveright##1{##1}%
 \setbox@ne\vbox{\Let@\halign{\setboxz@h
  {$\m@th\@lign\displaystyle{}##$}\ifonecr@\global\mwidth@\wdz@
  \global\mlineht@\ht\z@\fi\crcr#1\crcr}}}
\newbox\mtagbox@
\newdimen\ltwidth@
\newdimen\rtwidth@
\def\multline@l#1$${\iftagin@\DN@{\lmultline@@#1$$}\else
 \DN@{\setbox\mtagbox@\null\ltwidth@\z@\rtwidth@\z@
  \lmultline@@@#1$$}\fi\next@}
\def\lmultline@@#1\endmultline\tag#2$${%
 \setbox\mtagbox@\hbox{\maketag@#2\maketag@}
 \lmmeasure@#1\endmultline\dimen@\mwidth@\advance\dimen@\wd\mtagbox@
 \advance\dimen@\multlinetaggap@                                            
 \ifdim\dimen@>\displaywidth\ltwidth@\z@\else\ltwidth@\wd\mtagbox@\fi       
 \lmultline@@@#1\endmultline$$}
\def\lmultline@@@{\displ@y
 \def\shoveright##1{##1\hfilneg\hskip\multlinegap@}%
 \def\shoveleft##1{\setboxz@h{$\m@th\displaystyle{}##1$}%
  \setbox@ne\hbox{$\m@th\displaystyle##1$}%
  \hfilneg
  \iftagin@
   \ifdim\ltwidth@>\z@\hskip\ltwidth@\hskip\multlinetaggap@\fi
  \else\hskip\multlinegap@\fi\hskip.5\wd@ne\hskip-.5\wdz@##1}
  \halign\bgroup\Let@\hbox to\displaywidth
   {\strut@$\m@th\displaystyle\hfil{}##\hfil$}\crcr
   \hfilneg                                                                 
   \iftagin@                                                                
    \ifdim\ltwidth@>\z@                                                     
     \box\mtagbox@\hskip\multlinetaggap@                                    
    \else
     \rlap{\vbox{\normalbaselines\hbox{\strut@\box\mtagbox@}%
     \vbox to\mlineht@{}}}\fi                                               
   \else\hskip\multlinegap@\fi}                                             
\def\multline@r#1$${\iftagin@\DN@{\rmultline@@#1$$}\else
 \DN@{\setbox\mtagbox@\null\ltwidth@\z@\rtwidth@\z@
  \rmultline@@@#1$$}\fi\next@}
\def\rmultline@@#1\endmultline\tag#2$${\ltwidth@\z@
 \setbox\mtagbox@\hbox{\maketag@#2\maketag@}%
 \rmmeasure@#1\endmultline\dimen@\mwidth@\advance\dimen@\wd\mtagbox@
 \advance\dimen@\multlinetaggap@
 \ifdim\dimen@>\displaywidth\rtwidth@\z@\else\rtwidth@\wd\mtagbox@\fi
 \rmultline@@@#1\endmultline$$}
\def\rmultline@@@{\displ@y
 \def\shoveright##1{##1\hfilneg\iftagin@\ifdim\rtwidth@>\z@
  \hskip\rtwidth@\hskip\multlinetaggap@\fi\else\hskip\multlinegap@\fi}%
 \def\shoveleft##1{\setboxz@h{$\m@th\displaystyle{}##1$}%
  \setbox@ne\hbox{$\m@th\displaystyle##1$}%
  \hfilneg\hskip\multlinegap@\hskip.5\wd@ne\hskip-.5\wdz@##1}%
 \halign\bgroup\Let@\hbox to\displaywidth
  {\strut@$\m@th\displaystyle\hfil{}##\hfil$}\crcr
 \hfilneg\hskip\multlinegap@}
\def\endmultline{\iftagsleft@\expandafter\lendmultline@\else
 \expandafter\rendmultline@\fi}
\def\lendmultline@{\hfilneg\hskip\multlinegap@\crcr\egroup}
\def\rendmultline@{\iftagin@                                                
 \ifdim\rtwidth@>\z@                                                        
  \hskip\multlinetaggap@\box\mtagbox@                                       
 \else\llap{\vtop{\normalbaselines\null\hbox{\strut@\box\mtagbox@}}}\fi     
 \else\hskip\multlinegap@\fi                                                
 \hfilneg\crcr\egroup}
\def\bmod{\mskip-\medmuskip\mkern5mu\mathbin{\fam\z@ mod}\penalty900
 \mkern5mu\mskip-\medmuskip}
\def\pmod#1{\allowbreak\ifinner\mkern8mu\else\mkern18mu\fi
 ({\fam\z@ mod}\,\,#1)}
\def\pod#1{\allowbreak\ifinner\mkern8mu\else\mkern18mu\fi(#1)}
\def\mod#1{\allowbreak\ifinner\mkern12mu\else\mkern18mu\fi{\fam\z@ mod}\,\,#1}
\newcount\cfraccount@
\def\cfrac{\bgroup\bgroup\advance\cfraccount@\@ne\strut
 \iffalse{\fi\def\\{\over\displaystyle}\iffalse}\fi}
\def\lcfrac{\bgroup\bgroup\advance\cfraccount@\@ne\strut
 \iffalse{\fi\def\\{\hfill\over\displaystyle}\iffalse}\fi}
\def\rcfrac{\bgroup\bgroup\advance\cfraccount@\@ne\strut\hfill
 \iffalse{\fi\def\\{\over\displaystyle}\iffalse}\fi}
\def\gloop@#1\repeat{\gdef\body{#1}\iterate}
\def\endcfrac{\gloop@\ifnum\cfraccount@>\z@\global\advance\cfraccount@\m@ne
 \egroup\hskip-\nulldelimiterspace\egroup\repeat}
\def\binrel@#1{\setboxz@h{\thinmuskip0mu
  \medmuskip\m@ne mu\thickmuskip\@ne mu$#1\m@th$}%
 \setbox@ne\hbox{\thinmuskip0mu\medmuskip\m@ne mu\thickmuskip
  \@ne mu${}#1{}\m@th$}%
 \setbox\tw@\hbox{\hskip\wd@ne\hskip-\wdz@}}
\def\overset#1\to#2{\binrel@{#2}\ifdim\wd\tw@<\z@
 \mathbin{\mathop{\kern\z@#2}\limits^{#1}}\else\ifdim\wd\tw@>\z@
 \mathrel{\mathop{\kern\z@#2}\limits^{#1}}\else
 {\mathop{\kern\z@#2}\limits^{#1}}{}\fi\fi}
\def\underset#1\to#2{\binrel@{#2}\ifdim\wd\tw@<\z@
 \mathbin{\mathop{\kern\z@#2}\limits_{#1}}\else\ifdim\wd\tw@>\z@
 \mathrel{\mathop{\kern\z@#2}\limits_{#1}}\else
 {\mathop{\kern\z@#2}\limits_{#1}}{}\fi\fi}
\def\oversetbrace#1\to#2{\overbrace{#2}^{#1}}
\def\undersetbrace#1\to#2{\underbrace{#2}_{#1}}
\def\sideset#1\and#2\to#3{%
 \setbox@ne\hbox{$\dsize{\vphantom{#3}}#1{#3}\m@th$}%
 \setbox\tw@\hbox{$\dsize{#3}#2\m@th$}%
 \hskip\wd@ne\hskip-\wd\tw@\mathop{\hskip\wd\tw@\hskip-\wd@ne
  {\vphantom{#3}}#1{#3}#2}}
\def\rightarrowfill@#1{$#1\m@th\mathord-\mkern-6mu\cleaders
 \hbox{$#1\mkern-2mu\mathord-\mkern-2mu$}\hfill
 \mkern-6mu\mathord\rightarrow$}
\def\leftarrowfill@#1{$#1\m@th\mathord\leftarrow\mkern-6mu\cleaders
 \hbox{$#1\mkern-2mu\mathord-\mkern-2mu$}\hfill\mkern-6mu\mathord-$}
\def\leftrightarrowfill@#1{$#1\m@th\mathord\leftarrow\mkern-6mu\cleaders
 \hbox{$#1\mkern-2mu\mathord-\mkern-2mu$}\hfill
 \mkern-6mu\mathord\rightarrow$}
\def\overrightarrow{\mathpalette\overrightarrow@}
\def\overrightarrow@#1#2{\vbox{\ialign{##\crcr\rightarrowfill@#1\crcr
 \noalign{\kern-\ex@\nointerlineskip}$\m@th\hfil#1#2\hfil$\crcr}}}

\def\overleftarrow{\mathpalette\overleftarrow@}
\def\overleftarrow@#1#2{\vbox{\ialign{##\crcr\leftarrowfill@#1\crcr
 \noalign{\kern-\ex@\nointerlineskip}$\m@th\hfil#1#2\hfil$\crcr}}}
\def\overleftrightarrow{\mathpalette\overleftrightarrow@}
\def\overleftrightarrow@#1#2{\vbox{\ialign{##\crcr\leftrightarrowfill@#1\crcr
 \noalign{\kern-\ex@\nointerlineskip}$\m@th\hfil#1#2\hfil$\crcr}}}
\def\underrightarrow{\mathpalette\underrightarrow@}
\def\underrightarrow@#1#2{\vtop{\ialign{##\crcr$\m@th\hfil#1#2\hfil$\crcr
 \noalign{\nointerlineskip}\rightarrowfill@#1\crcr}}}

\def\underleftarrow{\mathpalette\underleftarrow@}
\def\underleftarrow@#1#2{\vtop{\ialign{##\crcr$\m@th\hfil#1#2\hfil$\crcr
 \noalign{\nointerlineskip}\leftarrowfill@#1\crcr}}}
\def\underleftrightarrow{\mathpalette\underleftrightarrow@}
\def\underleftrightarrow@#1#2{\vtop{\ialign{##\crcr$\m@th\hfil#1#2\hfil$\crcr
 \noalign{\nointerlineskip}\leftrightarrowfill@#1\crcr}}}
\let\DOTSI\relax
\let\DOTSB\relax

\newif\ifmath@
{\uccode`7=`\\ \uccode`8=`m \uccode`9=`a \uccode`0=`t \uccode`!=`h
 \uppercase{\gdef\math@#1#2#3#4#5#6\math@{\global\math@false\ifx 7#1\ifx 8#2%
 \ifx 9#3\ifx 0#4\ifx !#5\xdef\meaning@{#6}\global\math@true\fi\fi\fi\fi\fi}}}
\newif\ifmathch@
{\uccode`7=`c \uccode`8=`h \uccode`9=`\"
 \uppercase{\gdef\mathch@#1#2#3#4#5#6\mathch@{\global\mathch@false
  \ifx 7#1\ifx 8#2\ifx 9#5\global\mathch@true\xdef\meaning@{9#6}\fi\fi\fi}}}
\newcount\classnum@
\def\getmathch@#1.#2\getmathch@{\classnum@#1 \divide\classnum@4096
 \ifcase\number\classnum@\or\or\gdef\thedots@{\dotsb@}\or
 \gdef\thedots@{\dotsb@}\fi}
\newif\ifmathbin@
{\uccode`4=`b \uccode`5=`i \uccode`6=`n
 \uppercase{\gdef\mathbin@#1#2#3{\relaxnext@
  \DNii@##1\mathbin@{\ifx\space@\next\global\mathbin@true\fi}%
 \global\mathbin@false\DN@##1\mathbin@{}%
 \ifx 4#1\ifx 5#2\ifx 6#3\DN@{\FN@\nextii@}\fi\fi\fi\next@}}}
\newif\ifmathrel@
{\uccode`4=`r \uccode`5=`e \uccode`6=`l
 \uppercase{\gdef\mathrel@#1#2#3{\relaxnext@
  \DNii@##1\mathrel@{\ifx\space@\next\global\mathrel@true\fi}%
 \global\mathrel@false\DN@##1\mathrel@{}%
 \ifx 4#1\ifx 5#2\ifx 6#3\DN@{\FN@\nextii@}\fi\fi\fi\next@}}}
\newif\ifmacro@
{\uccode`5=`m \uccode`6=`a \uccode`7=`c
 \uppercase{\gdef\macro@#1#2#3#4\macro@{\global\macro@false
  \ifx 5#1\ifx 6#2\ifx 7#3\global\macro@true
  \xdef\meaning@{\macro@@#4\macro@@}\fi\fi\fi}}}
\def\macro@@#1->#2\macro@@{#2}
\newif\ifDOTS@
\newcount\DOTSCASE@
{\uccode`6=`\\ \uccode`7=`D \uccode`8=`O \uccode`9=`T \uccode`0=`S
 \uppercase{\gdef\DOTS@#1#2#3#4#5{\global\DOTS@false\DN@##1\DOTS@{}%
  \ifx 6#1\ifx 7#2\ifx 8#3\ifx 9#4\ifx 0#5\let\next@\DOTS@@\fi\fi\fi\fi\fi
  \next@}}}
{\uccode`3=`B \uccode`4=`I \uccode`5=`X
 \uppercase{\gdef\DOTS@@#1{\relaxnext@
  \DNii@##1\DOTS@{\ifx\space@\next\global\DOTS@true\fi}%
  \DN@{\FN@\nextii@}%
  \ifx 3#1\global\DOTSCASE@\z@\else
  \ifx 4#1\global\DOTSCASE@\@ne\else
  \ifx 5#1\global\DOTSCASE@\tw@\else\DN@##1\DOTS@{}%
  \fi\fi\fi\next@}}}
\newif\ifnot@
{\uccode`5=`\\ \uccode`6=`n \uccode`7=`o \uccode`8=`t
 \uppercase{\gdef\not@#1#2#3#4{\relaxnext@
  \DNii@##1\not@{\ifx\space@\next\global\not@true\fi}%
 \global\not@false\DN@##1\not@{}%
 \ifx 5#1\ifx 6#2\ifx 7#3\ifx 8#4\DN@{\FN@\nextii@}\fi\fi\fi
 \fi\next@}}}
\newif\ifkeybin@
\def\keybin@{\keybin@true
 \ifx\next+\else\ifx\next=\else\ifx\next<\else\ifx\next>\else\ifx\next-\else
 \ifx\next*\else\ifx\next:\else\keybin@false\fi\fi\fi\fi\fi\fi\fi}
\def\dots{\RIfM@\expandafter\mdots@\else\expandafter\tdots@\fi}
\def\tdots@{\unskip\relaxnext@
 \DN@{$\m@th\mathinner{\ldotp\ldotp\ldotp}\,
   \ifx\next,\,$\else\ifx\next.\,$\else\ifx\next;\,$\else\ifx\next:\,$\else
   \ifx\next?\,$\else\ifx\next!\,$\else$ \fi\fi\fi\fi\fi\fi}%
 \ \FN@\next@}
\def\mdots@{\FN@\mdots@@}
\def\mdots@@{\gdef\thedots@{\dotso@}
 \ifx\next\boldkey\gdef\thedots@\boldkey{\boldkeydots@}\else                
 \ifx\next\boldsymbol\gdef\thedots@\boldsymbol{\boldsymboldots@}\else       
 \ifx,\next\gdef\thedots@{\dotsc}
 \else\ifx\not\next\gdef\thedots@{\dotsb@}
 \else\keybin@
 \ifkeybin@\gdef\thedots@{\dotsb@}
 \else\xdef\meaning@{\meaning\next..........}\xdef\meaning@@{\meaning@}
  \expandafter\math@\meaning@\math@
  \ifmath@
   \expandafter\mathch@\meaning@\mathch@
   \ifmathch@\expandafter\getmathch@\meaning@\getmathch@\fi                 
  \else\expandafter\macro@\meaning@@\macro@                                 
  \ifmacro@                                                                
   \expandafter\not@\meaning@\not@\ifnot@\gdef\thedots@{\dotsb@}
  \else\expandafter\DOTS@\meaning@\DOTS@
  \ifDOTS@
   \ifcase\number\DOTSCASE@\gdef\thedots@{\dotsb@}%
    \or\gdef\thedots@{\dotsi}\else\fi                                      
  \else\expandafter\math@\meaning@\math@                                   
  \ifmath@\expandafter\mathbin@\meaning@\mathbin@
  \ifmathbin@\gdef\thedots@{\dotsb@}
  \else\expandafter\mathrel@\meaning@\mathrel@
  \ifmathrel@\gdef\thedots@{\dotsb@}
  \fi\fi\fi\fi\fi\fi\fi\fi\fi\fi\fi\fi
 \thedots@}
\def\plainldots@{\mathinner{\ldotp\ldotp\ldotp}}
\def\plaincdots@{\mathinner{\cdotp\cdotp\cdotp}}
\def\dotsi{\!\plaincdots@}
\let\dotsb@\plaincdots@
\newif\ifextra@
\newif\ifrightdelim@
\def\rightdelim@{\global\rightdelim@true                                    
 \ifx\next)\else                                                            
 \ifx\next]\else
 \ifx\next\rbrack\else
 \ifx\next\}\else
 \ifx\next\rbrace\else
 \ifx\next\rangle\else
 \ifx\next\rceil\else
 \ifx\next\rfloor\else
 \ifx\next\rgroup\else
 \ifx\next\rmoustache\else
 \ifx\next\right\else
 \ifx\next\bigr\else
 \ifx\next\biggr\else
 \ifx\next\Bigr\else                                                        
 \ifx\next\Biggr\else\global\rightdelim@false
 \fi\fi\fi\fi\fi\fi\fi\fi\fi\fi\fi\fi\fi\fi\fi}
\def\extra@{%
 \global\extra@false\rightdelim@\ifrightdelim@\global\extra@true            
 \else\ifx\next$\global\extra@true                                          
 \else\xdef\meaning@{\meaning\next..........}
 \expandafter\macro@\meaning@\macro@\ifmacro@                               
 \expandafter\DOTS@\meaning@\DOTS@
 \ifDOTS@
 \ifnum\DOTSCASE@=\tw@\global\extra@true                                    
 \fi\fi\fi\fi\fi}
\newif\ifbold@
\def\dotso@{\relaxnext@
 \ifbold@
  \let\next\delayed@
  \DNii@{\extra@\plainldots@\ifextra@\,\fi}%
 \else
  \DNii@{\DN@{\extra@\plainldots@\ifextra@\,\fi}\FN@\next@}%
 \fi
 \nextii@}
\def\extrap@#1{%
 \ifx\next,\DN@{#1\,}\else
 \ifx\next;\DN@{#1\,}\else
 \ifx\next.\DN@{#1\,}\else\extra@
 \ifextra@\DN@{#1\,}\else
 \let\next@#1\fi\fi\fi\fi\next@}
\def\ldots{\DN@{\extrap@\plainldots@}%
 \FN@\next@}
\def\cdots{\DN@{\extrap@\plaincdots@}%
 \FN@\next@}

\def\dotsc{\relaxnext@
 \DN@{\ifx\next;\plainldots@\,\else
  \ifx\next.\plainldots@\,\else\extra@\plainldots@
  \ifextra@\,\fi\fi\fi}%
 \FN@\next@}
\def\cdot{\mathchar"2201 }
\def\longrightarrow{\DOTSB\relbar\joinrel\rightarrow}

\def\longleftarrow{\DOTSB\leftarrow\joinrel\relbar}

\def\mapsto{\DOTSB\mapstochar\rightarrow}

\def\hookrightarrow{\DOTSB\lhook\joinrel\rightarrow}
\def\hookleftarrow{\DOTSB\leftarrow\joinrel\rhook}

\def\dddot#1{{\mathop{#1}\limits^{\vbox to-1.4\ex@{\kern-\tw@\ex@
 \hbox{\rm...}\vss}}}}
\def\ddddot#1{{\mathop{#1}\limits^{\vbox to-1.4\ex@{\kern-\tw@\ex@
 \hbox{\rm....}\vss}}}}
\def\sphat{^{\mathchoice{}{}%
 {\,\,\botsmash{\hbox{\lower4\ex@\hbox{$\m@th\widehat{\null}$}}}}%
 {\,\botsmash{\hbox{\lower3\ex@\hbox{$\m@th\hat{\null}$}}}}}}

\def\spacute{^{\!\botsmash{\hbox{\lower\@ne ex\hbox{\'{}}}}}}
\def\spgrave{^{\mathchoice{}{}{}{\!}%
 \botsmash{\hbox{\lower\@ne ex\hbox{\`{}}}}}}
\def\spdot{^{\hbox{\raise\ex@\hbox{\rm.}}}}
\def\spddot{^{\hbox{\raise\ex@\hbox{\rm..}}}}
\def\spdddot{^{\hbox{\raise\ex@\hbox{\rm...}}}}
\def\spddddot{^{\hbox{\raise\ex@\hbox{\rm....}}}}
\def\spbreve{^{\!\botsmash{\hbox{\lower4\ex@\hbox{\u{}}}}}}

\def\textonlyfont@#1#2{\def#1{\RIfM@
 \Err@{Use \string#1\space only in text}\else#2\fi}}
\textonlyfont@\rm\tenrm
\textonlyfont@\it\tenit
\textonlyfont@\sl\tensl
\textonlyfont@\bf\tenbf
\def\oldnos#1{\RIfM@{\mathcode`\,="013B \fam\@ne#1}\else
 \leavevmode\hbox{$\m@th\mathcode`\,="013B \fam\@ne#1$}\fi}
\def\text{\RIfM@\expandafter\text@\else\expandafter\text@@\fi}
\def\text@@#1{\leavevmode\hbox{#1}}
\def\mathhexbox@#1#2#3{\text{$\m@th\mathchar"#1#2#3$}}
\def\dag{{\mathhexbox@279}}
\def\ddag{{\mathhexbox@27A}}
\def\S{{\mathhexbox@278}}
\def\P{{\mathhexbox@27B}}
\newif\iffirstchoice@
\firstchoice@true
\def\text@#1{\mathchoice
 {\hbox{\everymath{\displaystyle}\def\textfonti{\the\textfont\@ne}%
  \def\textfontii{\the\textfont\tw@}\textdef@@ T#1}}
 {\hbox{\firstchoice@false
  \everymath{\textstyle}\def\textfonti{\the\textfont\@ne}%
  \def\textfontii{\the\textfont\tw@}\textdef@@ T#1}}
 {\hbox{\firstchoice@false
  \everymath{\scriptstyle}\def\textfonti{\the\scriptfont\@ne}%
  \def\textfontii{\the\scriptfont\tw@}\textdef@@ S\rm#1}}
 {\hbox{\firstchoice@false
  \everymath{\scriptscriptstyle}\def\textfonti
  {\the\scriptscriptfont\@ne}%
  \def\textfontii{\the\scriptscriptfont\tw@}\textdef@@ s\rm#1}}}
\def\textdef@@#1{\textdef@#1\rm\textdef@#1\bf\textdef@#1\sl\textdef@#1\it}
\def\rmfam{0}
\def\textdef@#1#2{%
 \DN@{\csname\expandafter\eat@\string#2fam\endcsname}%
 \if S#1\edef#2{\the\scriptfont\next@\relax}%
 \else\if s#1\edef#2{\the\scriptscriptfont\next@\relax}%
 \else\edef#2{\the\textfont\next@\relax}\fi\fi}
\scriptfont\itfam\tenit \scriptscriptfont\itfam\tenit
\scriptfont\slfam\tensl \scriptscriptfont\slfam\tensl
\newif\iftopfolded@
\newif\ifbotfolded@
\def\topfoldedtext{\topfolded@true\botfolded@false\foldedtext@}
\def\botfoldedtext{\botfolded@true\topfolded@false\foldedtext@}
\def\foldedtext{\topfolded@false\botfolded@false\foldedtext@}
\Invalid@\foldedwidth
\def\foldedtext@{\relaxnext@
 \DN@{\ifx\next\foldedwidth\let\next@\nextii@\else
  \DN@{\nextii@\foldedwidth{.3\hsize}}\fi\next@}%
 \DNii@\foldedwidth##1##2{\setbox\z@\vbox
  {\normalbaselines\hsize##1\relax
  \tolerance1600 \noindent\ignorespaces##2}\ifbotfolded@\boxz@\else
  \iftopfolded@\vtop{\unvbox\z@}\else\vcenter{\boxz@}\fi\fi}%
 \FN@\next@}
\def\bold{\RIfM@\expandafter\bold@\else
 \expandafter\nonmatherr@\expandafter\bold\fi}
\def\bold@#1{{\bold@@{#1}}}
\def\bold@@#1{\fam\bffam\relax#1}
\def\slanted{\RIfM@\expandafter\slanted@\else
 \expandafter\nonmatherr@\expandafter\slanted\fi}
\def\slanted@#1{{\slanted@@{#1}}}
\def\slanted@@#1{\fam\slfam\relax#1}
\def\roman{\RIfM@\expandafter\roman@\else
 \expandafter\nonmatherr@\expandafter\roman\fi}
\def\roman@#1{{\roman@@{#1}}}
\def\roman@@#1{\fam\rmfam\relax#1}
\def\italic{\RIfM@\expandafter\italic@\else
 \expandafter\nonmatherr@\expandafter\italic\fi}
\def\italic@#1{{\italic@@{#1}}}
\def\italic@@#1{\fam\itfam\relax#1}
\def\Cal{\RIfM@\expandafter\Cal@\else
 \expandafter\nonmatherr@\expandafter\Cal\fi}
\def\Cal@#1{{\Cal@@{#1}}}
\def\Cal@@#1{\noaccents@\fam\tw@#1}
\mathchardef\Gamma="0000
\mathchardef\Delta="0001
\mathchardef\Theta="0002
\mathchardef\Lambda="0003
\mathchardef\Xi="0004
\mathchardef\Pi="0005
\mathchardef\Sigma="0006
\mathchardef\Upsilon="0007
\mathchardef\Phi="0008
\mathchardef\Psi="0009
\mathchardef\Omega="000A
\mathchardef\varGamma="0100
\mathchardef\varDelta="0101
\mathchardef\varTheta="0102
\mathchardef\varLambda="0103
\mathchardef\varXi="0104
\mathchardef\varPi="0105
\mathchardef\varSigma="0106
\mathchardef\varUpsilon="0107
\mathchardef\varPhi="0108
\mathchardef\varPsi="0109
\mathchardef\varOmega="010A
\newif\ifmsamloaded@
\newif\ifmsbmloaded@
\newif\ifeufmloaded@
\let\alloc@@\alloc@
\def\hexnumber@#1{\ifcase#1 0\or 1\or 2\or 3\or 4\or 5\or 6\or 7\or 8\or
 9\or A\or B\or C\or D\or E\or F\fi}
\edef\bffam@{\hexnumber@\bffam}
\def\loadmsam{\msamloaded@true
 \font@\tenmsa=msam10
 \font@\sevenmsa=msam7
 \font@\fivemsa=msam5
 \alloc@@8\fam\chardef\sixt@@n\msafam
 \textfont\msafam=\tenmsa
 \scriptfont\msafam=\sevenmsa
 \scriptscriptfont\msafam=\fivemsa
 \edef\msafam@{\hexnumber@\msafam}%
 \mathchardef\dabar@"0\msafam@39
 \def\dashrightarrow{\mathrel{\dabar@\dabar@\mathchar"0\msafam@4B}}%
 \def\dashleftarrow{\mathrel{\mathchar"0\msafam@4C\dabar@\dabar@}}%
 \let\dasharrow\dashrightarrow
 \def\ulcorner{\delimiter"4\msafam@70\msafam@70 }
 \def\urcorner{\delimiter"5\msafam@71\msafam@71 }
 \def\llcorner{\delimiter"4\msafam@78\msafam@78 }
 \def\lrcorner{\delimiter"5\msafam@79\msafam@79 }
 \def\yen{{\mathhexbox@\msafam@55 }}
 \def\checkmark{{\mathhexbox@\msafam@58 }}
 \def\circledR{{\mathhexbox@\msafam@72 }}
 \def\maltese{{\mathhexbox@\msafam@7A }}}
\def\loadmsbm{\msbmloaded@true
 \font@\tenmsb=msbm10
 \font@\sevenmsb=msbm7
 \font@\fivemsb=msbm5
 \alloc@@8\fam\chardef\sixt@@n\msbfam
 \textfont\msbfam=\tenmsb
 \scriptfont\msbfam=\sevenmsb
 \scriptscriptfont\msbfam=\fivemsb
 \edef\msbfam@{\hexnumber@\msbfam}%
 }
\def\widehat#1{\ifmsbmloaded@
  \setboxz@h{$\m@th#1$}\ifdim\wdz@>\tw@ em\mathaccent"0\msbfam@5B{#1}\else
  \mathaccent"0362{#1}\fi
 \else\mathaccent"0362{#1}\fi}
\def\widetilde#1{\ifmsbmloaded@
  \setboxz@h{$\m@th#1$}\ifdim\wdz@>\tw@ em\mathaccent"0\msbfam@5D{#1}\else
  \mathaccent"0365{#1}\fi
 \else\mathaccent"0365{#1}\fi}
\def\newsymbol#1#2#3#4#5{\define#1{}\let\next@\relax
 \ifnum#2=\@ne\ifmsamloaded@\let\next@\msafam@\fi\else
 \ifnum#2=\tw@\ifmsbmloaded@\let\next@\msbfam@\fi\fi\fi
 \ifx\next@\relax
  \ifnum#2>\tw@\Err@{\Invalid@@\string\newsymbol}\else
  \ifnum#2=\@ne\Err@{You must first \string\loadmsam}\else
   \Err@{You must first \string\loadmsbm}\fi\fi
 \else
  \mathchardef#1="#3\next@#4#5
 \fi}

\def\Bbb{\RIfM@\expandafter\Bbb@\else
 \expandafter\nonmatherr@\expandafter\Bbb\fi}
\def\Bbb@#1{{\Bbb@@{#1}}}
\def\Bbb@@#1{\noaccents@\fam\msbfam\relax#1}
\def\loadeufm{\eufmloaded@true
 \font@\teneufm=eufm10
 \font@\seveneufm=eufm7
 \font@\fiveeufm=eufm5
 \alloc@@8\fam\chardef\sixt@@n\eufmfam
 \textfont\eufmfam=\teneufm
 \scriptfont\eufmfam=\seveneufm
 \scriptscriptfont\eufmfam=\fiveeufm}
\def\frak{\RIfM@\expandafter\frak@\else
 \expandafter\nonmatherr@\expandafter\frak\fi}
\def\frak@#1{{\frak@@{#1}}}
\def\frak@@#1{\fam\eufmfam\relax#1}

\newif\ifcmmibloaded@
\newif\ifcmbsyloaded@
\def\loadbold{\cmmibloaded@true\cmbsyloaded@true
 \font@\tencmmib=cmmib10 \font@\sevencmmib=cmmib7 \font@\fivecmmib=cmmib5
 \skewchar\tencmmib='177 \skewchar\sevencmmib='177 \skewchar\fivecmmib='177
 \alloc@@8\fam\chardef\sixt@@n\cmmibfam
 \textfont\cmmibfam=\tencmmib
 \scriptfont\cmmibfam=\sevencmmib
 \scriptscriptfont\cmmibfam=\fivecmmib
 \edef\cmmibfam@{\hexnumber@\cmmibfam}%
 \font@\tencmbsy=cmbsy10 \font@\sevencmbsy=cmbsy7 \font@\fivecmbsy=cmbsy5
 \skewchar\tencmbsy='60 \skewchar\sevencmbsy='60 \skewchar\fivecmbsy='60
 \alloc@@8\fam\chardef\sixt@@n\cmbsyfam
 \textfont\cmbsyfam=\tencmbsy
 \scriptfont\cmbsyfam=\sevencmbsy
 \scriptscriptfont\cmbsyfam=\fivecmbsy
 \edef\cmbsyfam@{\hexnumber@\cmbsyfam}}
\def\mathchari@#1#2#3{\ifcmmibloaded@\mathchar"#1\cmmibfam@#2#3 \else
 \Err@{First bold symbol font not loaded}\fi}
\def\mathcharii@#1#2#3{\ifcmbsyloaded@\mathchar"#1\cmbsyfam@#2#3 \else
 \Err@{Second bold symbol font not loaded}\fi}
\def\boldkey#1{\ifcat\noexpand#1A%
  \ifcmmibloaded@{\fam\cmmibfam#1}\else
   \Err@{First bold symbol font not loaded}\fi
 \else
 \ifx#1!\mathchar"5\bffam@21 \else
 \ifx#1(\mathchar"4\bffam@28 \else\ifx#1)\mathchar"5\bffam@29 \else
 \ifx#1+\mathchar"2\bffam@2B \else\ifx#1:\mathchar"3\bffam@3A \else
 \ifx#1;\mathchar"6\bffam@3B \else\ifx#1=\mathchar"3\bffam@3D \else
 \ifx#1?\mathchar"5\bffam@3F \else\ifx#1[\mathchar"4\bffam@5B \else
 \ifx#1]\mathchar"5\bffam@5D \else
 \ifx#1,\mathchari@63B \else
 \ifx#1-\mathcharii@200 \else
 \ifx#1.\mathchari@03A \else
 \ifx#1/\mathchari@03D \else
 \ifx#1<\mathchari@33C \else
 \ifx#1>\mathchari@33E \else
 \ifx#1*\mathcharii@203 \else
 \ifx#1|\mathcharii@06A \else
 \ifx#10\bold0\else\ifx#11\bold1\else\ifx#12\bold2\else\ifx#13\bold3\else
 \ifx#14\bold4\else\ifx#15\bold5\else\ifx#16\bold6\else\ifx#17\bold7\else
 \ifx#18\bold8\else\ifx#19\bold9\else
  \Err@{\string\boldkey\space can't be used with #1}%
 \fi\fi\fi\fi\fi\fi\fi\fi\fi\fi\fi\fi\fi\fi\fi
 \fi\fi\fi\fi\fi\fi\fi\fi\fi\fi\fi\fi\fi\fi}
\def\boldsymbol#1{%
 \DN@{\Err@{You can't use \string\boldsymbol\space with \string#1}#1}%
 \ifcat\noexpand#1A%
   \let\next@\relax
  \ifcmmibloaded@{\fam\cmmibfam#1}\else\Err@{First bold symbol
   font not loaded}\fi
 \else
  \xdef\meaning@{\meaning#1.........}%
  \expandafter\math@\meaning@\math@
  \ifmath@
   \expandafter\mathch@\meaning@\mathch@
   \ifmathch@
    \expandafter\boldsymbol@@\meaning@\boldsymbol@@
   \fi
  \else
   \expandafter\macro@\meaning@\macro@
   \expandafter\delim@\meaning@\delim@
   \ifdelim@
    \expandafter\delim@@\meaning@\delim@@
   \else
    \boldsymbol@{#1}%
   \fi
  \fi
 \fi
 \next@}
\def\mathhexboxii@#1#2{\ifcmbsyloaded@\mathhexbox@{\cmbsyfam@}{#1}{#2}\else
  \Err@{Second bold symbol font not loaded}\fi}
\def\boldsymbol@#1{\let\next@\relax\let\next#1%
 \ifx\next\cdot\mathcharii@201 \else
 \ifx\next\prime{{\null\mathcharii@030 \null}}\else
 \ifx\next\lbrack\mathchar"4\bffam@5B \else
 \ifx\next\rbrack\mathchar"5\bffam@5D \else
 \ifx\next\{\mathcharii@466 \else
 \ifx\next\lbrace\mathcharii@466 \else
 \ifx\next\}\mathcharii@567 \else
 \ifx\next\rbrace\mathcharii@567 \else
 \ifx\next\surd{{\mathcharii@170}}\else
 \ifx\next\S{{\mathhexboxii@78}}\else
 \ifx\next\P{{\mathhexboxii@7B}}\else
 \ifx\next\dag{{\mathhexboxii@79}}\else
 \ifx\next\ddag{{\mathhexboxii@7A}}\else
 \DN@{\Err@{You can't use \string\boldsymbol\space with \string#1}#1}%
 \fi\fi\fi\fi\fi\fi\fi\fi\fi\fi\fi\fi\fi}
\def\boldsymbol@@#1.#2\boldsymbol@@{\classnum@#1 \count@@@\classnum@        
 \divide\classnum@4096 \count@\classnum@                                    
 \multiply\count@4096 \advance\count@@@-\count@ \count@@\count@@@           
 \divide\count@@@\@cclvi \count@\count@@                                    
 \multiply\count@@@\@cclvi \advance\count@@-\count@@@                       
 \divide\count@@@\@cclvi                                                    
 \multiply\classnum@4096 \advance\classnum@\count@@                         
 \ifnum\count@@@=\z@                                                        
  \count@"\bffam@ \multiply\count@\@cclvi
  \advance\classnum@\count@
  \DN@{\mathchar\number\classnum@}%
 \else
  \ifnum\count@@@=\@ne                                                      
   \ifcmmibloaded@
   \count@"\cmmibfam@ \multiply\count@\@cclvi
   \advance\classnum@\count@
   \DN@{\mathchar\number\classnum@}%
   \else\DN@{\Err@{First bold symbol font not loaded}}\fi
  \else
   \ifnum\count@@@=\tw@                                                    
  \ifcmbsyloaded@
    \count@"\cmbsyfam@ \multiply\count@\@cclvi
    \advance\classnum@\count@
    \DN@{\mathchar\number\classnum@}%
  \else\DN@{\Err@{Second bold symbol font not loaded}}\fi
  \fi
 \fi
\fi}
\newif\ifdelim@
\newcount\delimcount@
{\uccode`6=`\\ \uccode`7=`d \uccode`8=`e \uccode`9=`l
 \uppercase{\gdef\delim@#1#2#3#4#5\delim@
  {\delim@false\ifx 6#1\ifx 7#2\ifx 8#3\ifx 9#4\delim@true
   \xdef\meaning@{#5}\fi\fi\fi\fi}}}
\def\delim@@#1"#2#3#4#5#6\delim@@{\if#32%
\let\next@\relax
 \ifcmbsyloaded@
 \mathcharii@#2#4#5 \else\Err@{Second bold family not loaded}\fi\fi}
\def\vert{\delimiter"026A30C }
\def\Vert{\delimiter"026B30D }
\let\|\Vert

\def\boldkeydots@#1{\bold@true\let\next=#1\let\delayed@=#1\mdots@@
 \boldkey#1\bold@false}  
\def\boldsymboldots@#1{\bold@true\let\next#1\let\delayed@#1\mdots@@
 \boldsymbol#1\bold@false}
\newif\ifeufbloaded@
\def\loadeufb{\eufbloaded@true
 \font@\teneufb=eufb10
 \font@\seveneufb=eufb7
 \font@\fiveeufb=eufb5
 \alloc@@8\fam\chardef\sixt@@n\eufbfam
 \textfont\eufbfam=\teneufb
 \scriptfont\eufbfam=\seveneufb
 \scriptscriptfont\eufbfam=\fiveeufb
 \edef\eufbfam@{\hexnumber@\eufbfam}}
\newif\ifeusmloaded@
\def\loadeusm{\eusmloaded@true
 \font@\teneusm=eusm10
 \font@\seveneusm=eusm7
 \font@\fiveeusm=eusm5
 \alloc@@8\fam\chardef\sixt@@n\eusmfam
 \textfont\eusmfam=\teneusm
 \scriptfont\eusmfam=\seveneusm
 \scriptscriptfont\eusmfam=\fiveeusm
 \edef\eusmfam@{\hexnumber@\eusmfam}}
\newif\ifeusbloaded@
\def\loadeusb{\eusbloaded@true
 \font@\teneusb=eusb10
 \font@\seveneusb=eusb7
 \font@\fiveeusb=eusb5
 \alloc@@8\fam\chardef\sixt@@n\eusbfam
 \textfont\eusbfam=\teneusb
 \scriptfont\eusbfam=\seveneusb
 \scriptscriptfont\eusbfam=\fiveeusb
 \edef\eusbfam@{\hexnumber@\eusbfam}}
\newif\ifeurmloaded@
\def\loadeurm{\eurmloaded@true
 \font@\teneurm=eurm10
 \font@\seveneurm=eurm7
 \font@\fiveeurm=eurm5
 \alloc@@8\fam\chardef\sixt@@n\eurmfam
 \textfont\eurmfam=\teneurm
 \scriptfont\eurmfam=\seveneurm
 \scriptscriptfont\eurmfam=\fiveeurm
 \edef\eurmfam@{\hexnumber@\eurmfam}}
\newif\ifeurbloaded@
\def\loadeurb{\eurbloaded@true
 \font@\teneurb=eurb10
 \font@\seveneurb=eurb7
 \font@\fiveeurb=eurb5
 \alloc@@8\fam\chardef\sixt@@n\eurbfam
 \textfont\eurbfam=\teneurb
 \scriptfont\eurbfam=\seveneurb
 \scriptscriptfont\eurbfam=\fiveeurb
 \edef\eurbfam@{\hexnumber@\eurbfam}}
\def\accentclass@{7}
\def\noaccents@{\def\accentclass@{0}}
\def\makeacc@#1#2{\def#1{\mathaccent"\accentclass@#2 }}
\makeacc@\hat{05E}
\makeacc@\check{014}
\makeacc@\tilde{07E}
\makeacc@\acute{013}
\makeacc@\grave{012}
\makeacc@\dot{05F}
\makeacc@\ddot{07F}
\makeacc@\breve{015}
\makeacc@\bar{016}

\newcount\skewcharcount@
\newcount\familycount@
\def\theskewchar@{\familycount@\@ne
 \global\skewcharcount@\the\skewchar\textfont\@ne                           
 \ifnum\fam>\m@ne\ifnum\fam<16
  \global\familycount@\the\fam\relax
  \global\skewcharcount@\the\skewchar\textfont\the\fam\relax\fi\fi          
 \ifnum\skewcharcount@>\m@ne
  \ifnum\skewcharcount@<128
  \multiply\familycount@256
  \global\advance\skewcharcount@\familycount@
  \global\advance\skewcharcount@28672
  \mathchar\skewcharcount@\else
  \global\skewcharcount@\m@ne\fi\else
 \global\skewcharcount@\m@ne\fi}                                            
\newcount\pointcount@
\def\getpoints@#1.#2\getpoints@{\pointcount@#1 }
\newdimen\accentdimen@
\newcount\accentmu@
\def\dimentomu@{\multiply\accentdimen@ 100
 \expandafter\getpoints@\the\accentdimen@\getpoints@
 \multiply\pointcount@18
 \divide\pointcount@\@m
 \global\accentmu@\pointcount@}
\def\Makeacc@#1#2{\def#1{\RIfM@\DN@{\mathaccent@
 {"\accentclass@#2 }}\else\DN@{\nonmatherr@{#1}}\fi\next@}}
\def\unbracefonts@{\let\Cal@\Cal@@\let\roman@\roman@@\let\bold@\bold@@
 \let\slanted@\slanted@@}
\def\mathaccent@#1#2{\ifnum\fam=\m@ne\xdef\thefam@{1}\else
 \xdef\thefam@{\the\fam}\fi                                                 
 \accentdimen@\z@                                                           
 \setboxz@h{\unbracefonts@$\m@th\fam\thefam@\relax#2$}
 \ifdim\accentdimen@=\z@\DN@{\mathaccent#1{#2}}
  \setbox@ne\hbox{\unbracefonts@$\m@th\fam\thefam@\relax#2\theskewchar@$}
  \setbox\tw@\hbox{$\m@th\ifnum\skewcharcount@=\m@ne\else
   \mathchar\skewcharcount@\fi$}
  \global\accentdimen@\wd@ne\global\advance\accentdimen@-\wdz@
  \global\advance\accentdimen@-\wd\tw@                                     
  \global\multiply\accentdimen@\tw@
  \dimentomu@\global\advance\accentmu@\@ne                                 
 \else\DN@{{\mathaccent#1{#2\mkern\accentmu@ mu}%
    \mkern-\accentmu@ mu}{}}\fi                                             
 \next@}\Makeacc@\Hat{05E}
\Makeacc@\Check{014}
\Makeacc@\Tilde{07E}
\Makeacc@\Acute{013}
\Makeacc@\Grave{012}
\Makeacc@\Dot{05F}
\Makeacc@\Ddot{07F}
\Makeacc@\Breve{015}
\Makeacc@\Bar{016}
\def\Vec{\RIfM@\DN@{\mathaccent@{"017E }}\else
 \DN@{\nonmatherr@\Vec}\fi\next@}
\def\newbox@{\alloc@4\box\chardef\insc@unt}
\def\accentedsymbol#1#2{\expandafter\newbox@\csname\expandafter
  \eat@\string#1@box\endcsname
 \expandafter\setbox\csname\expandafter\eat@
  \string#1@box\endcsname\hbox{$\m@th#2$}\define
  #1{\expandafter\copy\csname\expandafter\eat@\string#1@box\endcsname{}}}
\def\sqrt#1{\radical"270370 {#1}}
\let\underline@\underline
\let\overline@\overline
\def\underline#1{\underline@{#1}}
\def\overline#1{\overline@{#1}}
\Invalid@\leftroot
\Invalid@\uproot
\newcount\uproot@
\newcount\leftroot@
\def\root{\relaxnext@
  \DN@{\ifx\next\uproot\let\next@\nextii@\else
   \ifx\next\leftroot\let\next@\nextiii@\else
   \let\next@\plainroot@\fi\fi\next@}%
  \DNii@\uproot##1{\uproot@##1\relax\FN@\nextiv@}%
  \def\nextiv@{\ifx\next\space@\DN@. {\FN@\nextv@}\else
   \DN@.{\FN@\nextv@}\fi\next@.}%
  \def\nextv@{\ifx\next\leftroot\let\next@\nextvi@\else
   \let\next@\plainroot@\fi\next@}%
  \def\nextvi@\leftroot##1{\leftroot@##1\relax\plainroot@}%
   \def\nextiii@\leftroot##1{\leftroot@##1\relax\FN@\nextvii@}%
  \def\nextvii@{\ifx\next\space@
   \DN@. {\FN@\nextviii@}\else
   \DN@.{\FN@\nextviii@}\fi\next@.}%
  \def\nextviii@{\ifx\next\uproot\let\next@\nextix@\else
   \let\next@\plainroot@\fi\next@}%
  \def\nextix@\uproot##1{\uproot@##1\relax\plainroot@}%
  \bgroup\uproot@\z@\leftroot@\z@\FN@\next@}
\def\plainroot@#1\of#2{\setbox\rootbox\hbox{$\m@th\scriptscriptstyle{#1}$}%
 \mathchoice{\r@@t\displaystyle{#2}}{\r@@t\textstyle{#2}}
 {\r@@t\scriptstyle{#2}}{\r@@t\scriptscriptstyle{#2}}\egroup}
\def\r@@t#1#2{\setboxz@h{$\m@th#1\sqrt{#2}$}%
 \dimen@\ht\z@\advance\dimen@-\dp\z@
 \setbox@ne\hbox{$\m@th#1\mskip\uproot@ mu$}\advance\dimen@ by1.667\wd@ne
 \mkern-\leftroot@ mu\mkern5mu\raise.6\dimen@\copy\rootbox
 \mkern-10mu\mkern\leftroot@ mu\boxz@}
\def\boxed#1{\setboxz@h{$\m@th\displaystyle{#1}$}\dimen@.4\ex@
 \advance\dimen@3\ex@\advance\dimen@\dp\z@
 \hbox{\lower\dimen@\hbox{%
 \vbox{\hrule height.4\ex@
 \hbox{\vrule width.4\ex@\hskip3\ex@\vbox{\vskip3\ex@\boxz@\vskip3\ex@}%
 \hskip3\ex@\vrule width.4\ex@}\hrule height.4\ex@}%
 }}}
\let\ampersand@\relax
\newdimen\minaw@
\minaw@11.11128\ex@
\newdimen\minCDaw@
\minCDaw@2.5pc
\def\minCDarrowwidth#1{\RIfMIfI@\onlydmatherr@\minCDarrowwidth
 \else\minCDaw@#1\relax\fi\else\onlydmatherr@\minCDarrowwidth\fi}
\newif\ifCD@
\def\CD{\bgroup\vspace@\relax\let\ampersand@&\iffalse}\fi
 \CD@true\vcenter\bgroup\Let@\tabskip\z@skip\baselineskip20\ex@
 \lineskip3\ex@\lineskiplimit3\ex@\halign\bgroup
 &\hfill$\m@th##$\hfill\crcr}
\def\endCD{\crcr\egroup\egroup\egroup}
\newdimen\bigaw@
\atdef@>#1>#2>{\ampersand@                                                  
 \setboxz@h{$\m@th\ssize\;{#1}\;\;$}
 \setbox@ne\hbox{$\m@th\ssize\;{#2}\;\;$}
 \setbox\tw@\hbox{$\m@th#2$}
 \ifCD@\global\bigaw@\minCDaw@\else\global\bigaw@\minaw@\fi                 
 \ifdim\wdz@>\bigaw@\global\bigaw@\wdz@\fi
 \ifdim\wd@ne>\bigaw@\global\bigaw@\wd@ne\fi                                
 \ifCD@\hskip.5em\fi                                                        
 \ifdim\wd\tw@>\z@
  \mathrel{\mathop{\hbox to\bigaw@{\rightarrowfill}}\limits^{#1}_{#2}}
 \else\mathrel{\mathop{\hbox to\bigaw@{\rightarrowfill}}\limits^{#1}}\fi    
 \ifCD@\hskip.5em\fi                                                       
 \ampersand@}                                                              
\atdef@<#1<#2<{\ampersand@\setboxz@h{$\m@th\ssize\;\;{#1}\;$}%
 \setbox@ne\hbox{$\m@th\ssize\;\;{#2}\;$}\setbox\tw@\hbox{$\m@th#2$}%
 \ifCD@\global\bigaw@\minCDaw@\else\global\bigaw@\minaw@\fi
 \ifdim\wdz@>\bigaw@\global\bigaw@\wdz@\fi
 \ifdim\wd@ne>\bigaw@\global\bigaw@\wd@ne\fi
 \ifCD@\hskip.5em\fi
 \ifdim\wd\tw@>\z@
  \mathrel{\mathop{\hbox to\bigaw@{\leftarrowfill}}\limits^{#1}_{#2}}\else
  \mathrel{\mathop{\hbox to\bigaw@{\leftarrowfill}}\limits^{#1}}\fi
 \ifCD@\hskip.5em\fi\ampersand@}
\atdef@)#1)#2){\ampersand@
 \setboxz@h{$\m@th\ssize\;{#1}\;\;$}%
 \setbox@ne\hbox{$\m@th\ssize\;{#2}\;\;$}%
 \setbox\tw@\hbox{$\m@th#2$}%
 \ifCD@
 \global\bigaw@\minCDaw@\else\global\bigaw@\minaw@\fi
 \ifdim\wdz@>\bigaw@\global\bigaw@\wdz@\fi
 \ifdim\wd@ne>\bigaw@\global\bigaw@\wd@ne\fi
 \ifCD@\hskip.5em\fi
 \ifdim\wd\tw@>\z@
  \mathrel{\mathop{\hbox to\bigaw@{\rightarrowfill}}\limits^{#1}_{#2}}%
 \else\mathrel{\mathop{\hbox to\bigaw@{\rightarrowfill}}\limits^{#1}}\fi
 \ifCD@\hskip.5em\fi
 \ampersand@}
\atdef@(#1(#2({\ampersand@\setboxz@h{$\m@th\ssize\;\;{#1}\;$}%
 \setbox@ne\hbox{$\m@th\ssize\;\;{#2}\;$}\setbox\tw@\hbox{$\m@th#2$}%
 \ifCD@\global\bigaw@\minCDaw@\else\global\bigaw@\minaw@\fi
 \ifdim\wdz@>\bigaw@\global\bigaw@\wdz@\fi
 \ifdim\wd@ne>\bigaw@\global\bigaw@\wd@ne\fi
 \ifCD@\hskip.5em\fi
 \ifdim\wd\tw@>\z@
  \mathrel{\mathop{\hbox to\bigaw@{\leftarrowfill}}\limits^{#1}_{#2}}\else
  \mathrel{\mathop{\hbox to\bigaw@{\leftarrowfill}}\limits^{#1}}\fi
 \ifCD@\hskip.5em\fi\ampersand@}
\atdef@ A#1A#2A{\llap{$\m@th\vcenter{\hbox
 {$\ssize#1$}}$}\Big\uparrow\rlap{$\m@th\vcenter{\hbox{$\ssize#2$}}$}&&}
\atdef@ V#1V#2V{\llap{$\m@th\vcenter{\hbox
 {$\ssize#1$}}$}\Big\downarrow\rlap{$\m@th\vcenter{\hbox{$\ssize#2$}}$}&&}
\atdef@={&\hskip.5em\mathrel
 {\vbox{\hrule width\minCDaw@\vskip3\ex@\hrule width
 \minCDaw@}}\hskip.5em&}
\atdef@|{\Big\Vert&&}
\atdef@@\vert{\Big\Vert&&}
\def\pretend#1\haswidth#2{\setboxz@h{$\m@th\scriptstyle{#2}$}\hbox
 to\wdz@{\hfill$\m@th\scriptstyle{#1}$\hfill}}
\def\pmb{\RIfM@\expandafter\mathpalette\expandafter\pmb@\else
 \expandafter\pmb@@\fi}
\def\pmb@@#1{\leavevmode\setboxz@h{#1}\kern-.025em\copy\z@\kern-\wdz@
 \kern-.05em\copy\z@\kern-\wdz@\kern-.025em\raise.0433em\boxz@}
\def\binrel@@#1{\ifdim\wd2<\z@\mathbin{#1}\else\ifdim\wd\tw@>\z@
 \mathrel{#1}\else{#1}\fi\fi}
\newdimen\pmbraise@
\def\pmb@#1#2{\setbox\thr@@\hbox{$\m@th#1{#2}$}%
 \setbox4 \hbox{$\m@th#1\mkern.7794mu$}\pmbraise@\wd4
 \binrel@{#2}\binrel@@{\mkern-.45mu\copy\thr@@\kern-\wd\thr@@
 \mkern-.9mu\copy\thr@@\kern-\wd\thr@@\mkern-.45mu\raise\pmbraise@\box\thr@@}}
\def\documentstyle#1{\input #1.sty\relax}
\font\dummyft@=dummy
\fontdimen1 \dummyft@=\z@
\fontdimen2 \dummyft@=\z@
\fontdimen3 \dummyft@=\z@
\fontdimen4 \dummyft@=\z@
\fontdimen5 \dummyft@=\z@
\fontdimen6 \dummyft@=\z@
\fontdimen7 \dummyft@=\z@
\fontdimen8 \dummyft@=\z@
\fontdimen9 \dummyft@=\z@
\fontdimen10 \dummyft@=\z@
\fontdimen11 \dummyft@=\z@
\fontdimen12 \dummyft@=\z@
\fontdimen13 \dummyft@=\z@
\fontdimen14 \dummyft@=\z@
\fontdimen15 \dummyft@=\z@
\fontdimen16 \dummyft@=\z@
\fontdimen17 \dummyft@=\z@
\fontdimen18 \dummyft@=\z@
\fontdimen19 \dummyft@=\z@
\fontdimen20 \dummyft@=\z@
\fontdimen21 \dummyft@=\z@
\fontdimen22 \dummyft@=\z@
\def\fontlist@{\\{\tenrm}\\{\sevenrm}\\{\fiverm}\\{\teni}\\{\seveni}%
 \\{\fivei}\\{\tensy}\\{\sevensy}\\{\fivesy}\\{\tenex}\\{\tenbf}\\{\sevenbf}%
 \\{\fivebf}\\{\tensl}\\{\tenit}}
\def\font@#1=#2 {\rightappend@#1\to\fontlist@\font#1=#2 }
\def\dodummy@{{\def\\##1{\global\let##1\dummyft@}\fontlist@}}
\def\nopages@{\output={\setbox\z@\box255 \deadcycles\z@}%
 \alloc@5\toks\toksdef\@cclvi\output}
\let\galleys\nopages@
\newif\ifsyntax@
\newcount\countxviii@
\def\syntax{\syntax@true\dodummy@\countxviii@\count18
 \loop\ifnum\countxviii@>\m@ne\textfont\countxviii@=\dummyft@
 \scriptfont\countxviii@=\dummyft@\scriptscriptfont\countxviii@=\dummyft@
 \advance\countxviii@\m@ne\repeat                                           
 \dummyft@\tracinglostchars\z@\nopages@\frenchspacing\hbadness\@M}
\def\S@{S } \def\G@{G } \def\P@{P }
\newif\ifbadans@
\def\printoptions{\W@{Do you want S(yntax check),
  G(alleys) or P(ages)?^^JType S, G or P, follow by <return>: }\loop
 \read\m@ne to\ans@
 \xdef\next@{\def\noexpand\Ans@{\ans@}}\uppercase\expandafter{\next@}
 \ifx\Ans@\S@\badans@false\syntax\else
 \ifx\Ans@\G@\badans@false\galleys\else
 \ifx\Ans@\P@\badans@false\else
 \badans@true\fi\fi\fi
 \ifbadans@\W@{Type S, G or P, follow by <return>: }%
 \repeat}
\def\alloc@#1#2#3#4#5{\global\advance\count1#1by\@ne
 \ch@ck#1#4#2\allocationnumber=\count1#1
 \global#3#5=\allocationnumber
 \ifalloc@\wlog{\string#5=\string#2\the\allocationnumber}\fi}
\def\document{\def\alloclist@{}\def\fontlist@{}}
\let\enddocument\bye

\let\proclaim\undefined
\let\footnote\undefined
\let\=\undefined
\let\>\undefined

\catcode`\@=\active

%
%
\def\next{AMSPPT}\ifx\styname\next \endinput\fi
\catcode`\@=11
\def\styname{AMSPPT}
\def\styversion{2.0}
{\W@{\styname.STY - Version \styversion}\W@{}}
\hyphenation{acad-e-my acad-e-mies af-ter-thought anom-aly anom-alies
an-ti-deriv-a-tive an-tin-o-my an-tin-o-mies apoth-e-o-ses apoth-e-o-sis
ap-pen-dix ar-che-typ-al as-sign-a-ble as-sist-ant-ship as-ymp-tot-ic
asyn-chro-nous at-trib-uted at-trib-ut-able bank-rupt bank-rupt-cy
bi-dif-fer-en-tial blue-print busier busiest cat-a-stroph-ic
cat-a-stroph-i-cally con-gress cross-hatched data-base de-fin-i-tive
de-riv-a-tive dis-trib-ute dri-ver dri-vers eco-nom-ics econ-o-mist
elit-ist equi-vari-ant ex-quis-ite ex-tra-or-di-nary flow-chart
for-mi-da-ble forth-right friv-o-lous ge-o-des-ic ge-o-det-ic geo-met-ric
griev-ance griev-ous griev-ous-ly hexa-dec-i-mal ho-lo-no-my ho-mo-thetic
ideals idio-syn-crasy in-fin-ite-ly in-fin-i-tes-i-mal ir-rev-o-ca-ble
key-stroke lam-en-ta-ble light-weight mal-a-prop-ism man-u-script
mar-gin-al meta-bol-ic me-tab-o-lism meta-lan-guage me-trop-o-lis
met-ro-pol-i-tan mi-nut-est mol-e-cule mono-chrome mono-pole mo-nop-oly
mono-spline mo-not-o-nous mul-ti-fac-eted mul-ti-plic-able non-euclid-ean
non-iso-mor-phic non-smooth par-a-digm par-a-bol-ic pa-rab-o-loid
pa-ram-e-trize para-mount pen-ta-gon phe-nom-e-non post-script pre-am-ble
pro-ce-dur-al pro-hib-i-tive pro-hib-i-tive-ly pseu-do-dif-fer-en-tial
pseu-do-fi-nite pseu-do-nym qua-drat-ics quad-ra-ture qua-si-smooth
qua-si-sta-tion-ary qua-si-tri-an-gu-lar quin-tes-sence quin-tes-sen-tial
re-arrange-ment rec-tan-gle ret-ri-bu-tion retro-fit retro-fit-ted
right-eous right-eous-ness ro-bot ro-bot-ics sched-ul-ing se-mes-ter
semi-def-i-nite semi-ho-mo-thet-ic set-up se-vere-ly side-step sov-er-eign
spe-cious spher-oid spher-oid-al star-tling star-tling-ly
sta-tis-tics sto-chas-tic straight-est strange-ness strat-a-gem strong-hold
sum-ma-ble symp-to-matic syn-chro-nous topo-graph-i-cal tra-vers-a-ble
tra-ver-sal tra-ver-sals treach-ery turn-around un-at-tached un-err-ing-ly
white-space wide-spread wing-spread wretch-ed wretch-ed-ly Brown-ian
Eng-lish Euler-ian Feb-ru-ary Gauss-ian Grothen-dieck Hamil-ton-ian
Her-mit-ian Jan-u-ary Japan-ese Kor-te-weg Le-gendre Lip-schitz
Lip-schitz-ian Mar-kov-ian Noe-ther-ian No-vem-ber Rie-mann-ian
Schwarz-schild Sep-tem-ber}
\Invalid@\nofrills
\Invalid@\usualspace
\newif\ifnofrills@
\def\nofrills@#1#2{\relaxnext@
  \DN@{\ifx\next\nofrills
    \nofrills@true\let#2\relax\DN@\nofrills{\nextii@}%
  \else
    \nofrills@false\def#2{#1}\let\next@\nextii@\fi
\next@}}
\def\usualspace@#1{\ifnofrills@\def\usualspace{#1}\fi}
\def\addto#1#2{\csname \expandafter\eat@\string#1@\endcsname
  \expandafter{\the\csname \expandafter\eat@\string#1@\endcsname#2}}
\newdimen\bigsize@
\def\big@#1#2{{\hbox{$\left#2\vcenter to#1\bigsize@{}%
  \right.\nulldelimiterspace\z@\m@th$}}}
\def\big{\big@\@ne}
\def\Big{\big@{1.5}}
\def\bigg{\big@\tw@}
\def\Bigg{\big@{2.5}}
\def\raggedcenter@{\leftskip\z@ plus.4\hsize \rightskip\leftskip
 \parfillskip\z@ \parindent\z@ \spaceskip.3333em \xspaceskip.5em
 \pretolerance9999\tolerance9999 \exhyphenpenalty\@M
 \hyphenpenalty\@M \let\\\linebreak}
\def\upperspecialchars{\def\ss{SS}\let\i=I\let\j=J\let\ae\AE\let\oe\OE
  \let\o\O\let\aa\AA\let\l\L}
\def\uppercasetext@#1{%
  {\spaceskip1.2\fontdimen2\the\font plus1.2\fontdimen3\the\font
   \upperspecialchars\uctext@#1$\m@th\aftergroup\eat@$}}
\def\uctext@#1$#2${\endash@#1-\endash@$#2$\uctext@}
\def\endash@#1-#2\endash@{\uppercase{#1}\if\notempty{#2}--\endash@#2\endash@\fi}
\def\runaway@#1{\DN@{#1}\ifx\envir@\next@
  \Err@{You seem to have a missing or misspelled \string\end#1 ...}%
  \let\envir@\empty\fi}
\newif\iftemp@
\def\notempty#1{TT\fi\def\test@{#1}\ifx\test@\empty\temp@false
  \else\temp@true\fi \iftemp@}
\font@\tensmc=cmcsc10
\font@\sevenex=cmex7
\font@\sevenit=cmti7
\font@\eightrm=cmr8 
\font@\sixrm=cmr6 
\font@\eighti=cmmi8     \skewchar\eighti='177 
\font@\sixi=cmmi6       \skewchar\sixi='177   
\font@\eightsy=cmsy8    \skewchar\eightsy='60 
\font@\sixsy=cmsy6      \skewchar\sixsy='60   
\font@\eightex=cmex8
\font@\eightbf=cmbx8 
\font@\sixbf=cmbx6   
\font@\eightit=cmti8 
\font@\eightsl=cmsl8 
\font@\eightsmc=cmcsc8
\font@\eighttt=cmtt8 
\loadmsam
\loadmsbm
\loadeufm
\input amssym.tex\relax
\newtoks\tenpoint@
\def\tenpoint{\normalbaselineskip12\p@
 \abovedisplayskip12\p@ plus3\p@ minus9\p@
 \belowdisplayskip\abovedisplayskip
 \abovedisplayshortskip\z@ plus3\p@
 \belowdisplayshortskip7\p@ plus3\p@ minus4\p@
 \textonlyfont@\rm\tenrm \textonlyfont@\it\tenit
 \textonlyfont@\sl\tensl \textonlyfont@\bf\tenbf
 \textonlyfont@\smc\tensmc \textonlyfont@\tt\tentt
 \ifsyntax@ \def\big##1{{\hbox{$\left##1\right.$}}}%
  \let\Big\big \let\bigg\big \let\Bigg\big
 \else
  \textfont\z@=\tenrm  \scriptfont\z@=\sevenrm  \scriptscriptfont\z@=\fiverm
  \textfont\@ne=\teni  \scriptfont\@ne=\seveni  \scriptscriptfont\@ne=\fivei
  \textfont\tw@=\tensy \scriptfont\tw@=\sevensy \scriptscriptfont\tw@=\fivesy
  \textfont\thr@@=\tenex \scriptfont\thr@@=\sevenex
        \scriptscriptfont\thr@@=\sevenex
  \textfont\itfam=\tenit \scriptfont\itfam=\sevenit
        \scriptscriptfont\itfam=\sevenit
  \textfont\bffam=\tenbf \scriptfont\bffam=\sevenbf
        \scriptscriptfont\bffam=\fivebf
  \setbox\strutbox\hbox{\vrule height8.5\p@ depth3.5\p@ width\z@}%
  \setbox\strutbox@\hbox{\lower.5\normallineskiplimit\vbox{%
        \kern-\normallineskiplimit\copy\strutbox}}%
 \setbox\z@\vbox{\hbox{$($}\kern\z@}\bigsize@=1.2\ht\z@
 \fi
 \normalbaselines\rm\ex@.2326ex\jot3\ex@\the\tenpoint@}
\newtoks\eightpoint@
\def\eightpoint{\normalbaselineskip10\p@
 \abovedisplayskip10\p@ plus2.4\p@ minus7.2\p@
 \belowdisplayskip\abovedisplayskip
 \abovedisplayshortskip\z@ plus2.4\p@
 \belowdisplayshortskip5.6\p@ plus2.4\p@ minus3.2\p@
 \textonlyfont@\rm\eightrm \textonlyfont@\it\eightit
 \textonlyfont@\sl\eightsl \textonlyfont@\bf\eightbf
 \textonlyfont@\smc\eightsmc \textonlyfont@\tt\eighttt
 \ifsyntax@\def\big##1{{\hbox{$\left##1\right.$}}}%
  \let\Big\big \let\bigg\big \let\Bigg\big
 \else
  \textfont\z@=\eightrm \scriptfont\z@=\sixrm \scriptscriptfont\z@=\fiverm
  \textfont\@ne=\eighti \scriptfont\@ne=\sixi \scriptscriptfont\@ne=\fivei
  \textfont\tw@=\eightsy \scriptfont\tw@=\sixsy \scriptscriptfont\tw@=\fivesy
  \textfont\thr@@=\eightex \scriptfont\thr@@=\sevenex
   \scriptscriptfont\thr@@=\sevenex
  \textfont\itfam=\eightit \scriptfont\itfam=\sevenit
   \scriptscriptfont\itfam=\sevenit
  \textfont\bffam=\eightbf \scriptfont\bffam=\sixbf
   \scriptscriptfont\bffam=\fivebf
 \setbox\strutbox\hbox{\vrule height7\p@ depth3\p@ width\z@}%
 \setbox\strutbox@\hbox{\raise.5\normallineskiplimit\vbox{%
   \kern-\normallineskiplimit\copy\strutbox}}%
 \setbox\z@\vbox{\hbox{$($}\kern\z@}\bigsize@=1.2\ht\z@
 \fi
 \normalbaselines\eightrm\ex@.2326ex\jot3\ex@\the\eightpoint@}
\parindent1pc
\normallineskiplimit\p@
\newdimen\indenti \indenti=2pc
\def\pageheight#1{\vsize#1}
\def\pagewidth#1{\hsize#1%
   \captionwidth@\hsize \advance\captionwidth@-2\indenti}
\pagewidth{30pc} \pageheight{47pc}
\def\topmatter{%
 \ifx\undefined\msafam
 \else\font@\eightmsa=msam8 \font@\sixmsa=msam6
   \ifsyntax@\else \addto\tenpoint{\textfont\msafam=\tenmsa
              \scriptfont\msafam=\sevenmsa \scriptscriptfont\msafam=\fivemsa}%
     \addto\eightpoint{\textfont\msafam=\eightmsa \scriptfont\msafam=\sixmsa
              \scriptscriptfont\msafam=\fivemsa}%
   \fi
 \fi
 \ifx\undefined\msbfam
 \else\font@\eightmsb=msbm8 \font@\sixmsb=msbm6
   \ifsyntax@\else \addto\tenpoint{\textfont\msbfam=\tenmsb
         \scriptfont\msbfam=\sevenmsb \scriptscriptfont\msbfam=\fivemsb}%
     \addto\eightpoint{\textfont\msbfam=\eightmsb \scriptfont\msbfam=\sixmsb
         \scriptscriptfont\msbfam=\fivemsb}%
   \fi
 \fi
 \ifx\undefined\eufmfam
 \else \font@\eighteufm=eufm8 \font@\sixeufm=eufm6
   \ifsyntax@\else \addto\tenpoint{\textfont\eufmfam=\teneufm
       \scriptfont\eufmfam=\seveneufm \scriptscriptfont\eufmfam=\fiveeufm}%
     \addto\eightpoint{\textfont\eufmfam=\eighteufm
       \scriptfont\eufmfam=\sixeufm \scriptscriptfont\eufmfam=\fiveeufm}%
   \fi
 \fi
 \ifx\undefined\eufbfam
 \else \font@\eighteufb=eufb8 \font@\sixeufb=eufb6
   \ifsyntax@\else \addto\tenpoint{\textfont\eufbfam=\teneufb
      \scriptfont\eufbfam=\seveneufb \scriptscriptfont\eufbfam=\fiveeufb}%
    \addto\eightpoint{\textfont\eufbfam=\eighteufb
      \scriptfont\eufbfam=\sixeufb \scriptscriptfont\eufbfam=\fiveeufb}%
   \fi
 \fi
 \ifx\undefined\eusmfam
 \else \font@\eighteusm=eusm8 \font@\sixeusm=eusm6
   \ifsyntax@\else \addto\tenpoint{\textfont\eusmfam=\teneusm
       \scriptfont\eusmfam=\seveneusm \scriptscriptfont\eusmfam=\fiveeusm}%
     \addto\eightpoint{\textfont\eusmfam=\eighteusm
       \scriptfont\eusmfam=\sixeusm \scriptscriptfont\eusmfam=\fiveeusm}%
   \fi
 \fi
 \ifx\undefined\eusbfam
 \else \font@\eighteusb=eusb8 \font@\sixeusb=eusb6
   \ifsyntax@\else \addto\tenpoint{\textfont\eusbfam=\teneusb
       \scriptfont\eusbfam=\seveneusb \scriptscriptfont\eusbfam=\fiveeusb}%
     \addto\eightpoint{\textfont\eusbfam=\eighteusb
       \scriptfont\eusbfam=\sixeusb \scriptscriptfont\eusbfam=\fiveeusb}%
   \fi
 \fi
 \ifx\undefined\eurmfam
 \else \font@\eighteurm=eurm8 \font@\sixeurm=eurm6
   \ifsyntax@\else \addto\tenpoint{\textfont\eurmfam=\teneurm
       \scriptfont\eurmfam=\seveneurm \scriptscriptfont\eurmfam=\fiveeurm}%
     \addto\eightpoint{\textfont\eurmfam=\eighteurm
       \scriptfont\eurmfam=\sixeurm \scriptscriptfont\eurmfam=\fiveeurm}%
   \fi
 \fi
 \ifx\undefined\eurbfam
 \else \font@\eighteurb=eurb8 \font@\sixeurb=eurb6
   \ifsyntax@\else \addto\tenpoint{\textfont\eurbfam=\teneurb
       \scriptfont\eurbfam=\seveneurb \scriptscriptfont\eurbfam=\fiveeurb}%
    \addto\eightpoint{\textfont\eurbfam=\eighteurb
       \scriptfont\eurbfam=\sixeurb \scriptscriptfont\eurbfam=\fiveeurb}%
   \fi
 \fi
 \ifx\undefined\cmmibfam
 \else \font@\eightcmmib=cmmib8 \font@\sixcmmib=cmmib6
   \ifsyntax@\else \addto\tenpoint{\textfont\cmmibfam=\tencmmib
       \scriptfont\cmmibfam=\sevencmmib \scriptscriptfont\cmmibfam=\fivecmmib}%
    \addto\eightpoint{\textfont\cmmibfam=\eightcmmib
       \scriptfont\cmmibfam=\sixcmmib \scriptscriptfont\cmmibfam=\fivecmmib}%
   \fi
 \fi
 \ifx\undefined\cmbsyfam
 \else \font@\eightcmbsy=cmbsy8 \font@\sixcmbsy=cmbsy6
   \ifsyntax@\else \addto\tenpoint{\textfont\cmbsyfam=\tencmbsy
      \scriptfont\cmbsyfam=\sevencmbsy \scriptscriptfont\cmbsyfam=\fivecmbsy}%
    \addto\eightpoint{\textfont\cmbsyfam=\eightcmbsy
      \scriptfont\cmbsyfam=\sixcmbsy \scriptscriptfont\cmbsyfam=\fivecmbsy}%
   \fi
 \fi
 \let\topmatter\relax}
\def\chapterno@{\uppercase\expandafter{\romannumeral\chaptercount@}}
\newcount\chaptercount@
\def\chapter{\nofrills@{\afterassignment\chapterno@
                        CHAPTER \global\chaptercount@=}\chapter@
 \DNii@##1{\leavevmode\hskip-\leftskip
   \rlap{\vbox to\z@{\vss\centerline{\eightpoint
   \chapter@##1\unskip}\baselineskip2pc\null}}\hskip\leftskip
   \nofrills@false}%
 \FN@\next@}
\newbox\titlebox@
\def\title{\nofrills@{\uppercasetext@}\title@%
 \DNii@##1\endtitle{\global\setbox\titlebox@\vtop{\tenpoint\bf
 \raggedcenter@\ignorespaces
 \baselineskip1.3\baselineskip\title@{##1}\endgraf}%
 \ifmonograph@ \edef\next{\the\leftheadtoks}\ifx\next\empty
    \leftheadtext{##1}\fi
 \fi
 \edef\next{\the\rightheadtoks}\ifx\next\empty \rightheadtext{##1}\fi
 }\FN@\next@}
\newbox\authorbox@
\def\author#1\endauthor{\global\setbox\authorbox@
 \vbox{\tenpoint\smc\raggedcenter@\ignorespaces
 #1\endgraf}\relaxnext@ \edef\next{\the\leftheadtoks}%
 \ifx\next\empty\leftheadtext{#1}\fi}
\newbox\affilbox@
\def\affil#1\endaffil{\global\setbox\affilbox@
 \vbox{\tenpoint\raggedcenter@\ignorespaces#1\endgraf}}
\newcount\addresscount@
\addresscount@\z@
\def\address#1\endaddress{\global\advance\addresscount@\@ne
  \expandafter\gdef\csname address\number\addresscount@\endcsname
  {\vskip12\p@ minus6\p@\noindent\eightpoint\smc\ignorespaces#1\par}}
\def\email{\nofrills@{\eightpoint{\it E-mail\/}:\enspace}\email@
  \DNii@##1\endemail{%
  \expandafter\gdef\csname email\number\addresscount@\endcsname
  {\def\usualspace{{\it\enspace}}\smallskip\noindent\eightpoint\email@
  \ignorespaces##1\par}}%
 \FN@\next@}
\def\thedate@{}
\def\date#1\enddate{\gdef\thedate@{\tenpoint\ignorespaces#1\unskip}}
\def\thethanks@{}
\def\thanks#1\endthanks{\gdef\thethanks@{\eightpoint\ignorespaces#1.\unskip}}
\def\thekeywords@{}
\def\keywords{\nofrills@{{\it Key words and phrases.\enspace}}\keywords@
 \DNii@##1\endkeywords{\def\thekeywords@{\def\usualspace{{\it\enspace}}%
 \eightpoint\keywords@\ignorespaces##1\unskip.}}%
 \FN@\next@}
\def\thesubjclass@{}
\def\subjclass{\nofrills@{{\rm1980 {\it Mathematics Subject
   Classification\/} (1985 {\it Revision\/}).\enspace}}\subjclass@
 \DNii@##1\endsubjclass{\def\thesubjclass@{\def\usualspace
  {{\rm\enspace}}\eightpoint\subjclass@\ignorespaces##1\unskip.}}%
 \FN@\next@}
\newbox\abstractbox@
\def\abstract{\nofrills@{{\smc Abstract.\enspace}}\abstract@
 \DNii@{\setbox\abstractbox@\vbox\bgroup\noindent$$\vbox\bgroup
  \def\envir@{abstract}\advance\hsize-2\indenti
  \usualspace@{{\enspace}}\eightpoint \noindent\abstract@\ignorespaces}%
 \FN@\next@}
\def\endabstract{\par\unskip\egroup$$\egroup}
\def\widestnumber#1#2{\begingroup\let\head\null\let\subhead\empty
   \let\subsubhead\subhead
   \ifx#1\head\global\setbox\tocheadbox@\hbox{#2.\enspace}%
   \else\ifx#1\subhead\global\setbox\tocsubheadbox@\hbox{#2.\enspace}%
   \else\ifx#1\key\bgroup\let\endrefitem@\egroup
        \key#2\endrefitem@\global\refindentwd\wd\keybox@
   \else\ifx#1\no\bgroup\let\endrefitem@\egroup
        \no#2\endrefitem@\global\refindentwd\wd\nobox@
   \else\ifx#1\page\global\setbox\pagesbox@\hbox{\quad\bf#2}%
   \else\ifx#1\item\setboxz@h{#2}\global\rosteritemwd\wdz@
        \global\advance\rosteritemwd by.5\parindent
   \else\message{\string\widestnumber is not defined for this option
   (\string#1)}%
\fi\fi\fi\fi\fi\fi\endgroup}
\newif\ifmonograph@
\def\Monograph{\monograph@true \let\headmark\rightheadtext
  \let\varindent@\indent \def\headfont@{\bf}\def\proclaimfont@{\smc}%
  \def\demofont@{\smc}}
\let\varindent@\noindent
\newbox\tocheadbox@    \newbox\tocsubheadbox@
\newbox\tocbox@
\def\toc{\toc@{Contents}}
\def\newtocdefs{%
   \def \title##1\endtitle
       {\penaltyandskip@\z@\smallskipamount
        \hangindent\wd\tocheadbox@\noindent{\bf##1}}%
   \def \chapter##1{%
        Chapter \uppercase\expandafter{\romannumeral##1.\unskip}\enspace}%
   \def \specialhead##1\endspecialhead
       {\par\hangindent\wd\tocheadbox@ \noindent##1\par}%
   \def \head##1 ##2\endhead
       {\par\hangindent\wd\tocheadbox@ \noindent
        \if\notempty{##1}\hbox to\wd\tocheadbox@{\hfil##1\enspace}\fi
        ##2\par}%
   \def \subhead##1 ##2\endsubhead
       {\par\vskip-\parskip {\normalbaselines
        \advance\leftskip\wd\tocheadbox@
        \hangindent\wd\tocsubheadbox@ \noindent
        \if\notempty{##1}\hbox to\wd\tocsubheadbox@{##1\unskip\hfil}\fi
         ##2\par}}%
   \def \subsubhead##1 ##2\endsubsubhead
       {\par\vskip-\parskip {\normalbaselines
        \advance\leftskip\wd\tocheadbox@
        \hangindent\wd\tocsubheadbox@ \noindent
        \if\notempty{##1}\hbox to\wd\tocsubheadbox@{##1\unskip\hfil}\fi
        ##2\par}}}
\def\toc@#1{\relaxnext@
   \def\page##1%
       {\unskip\penalty0\null\hfil
        \rlap{\hbox to\wd\pagesbox@{\quad\hfil##1}}\hfilneg\penalty\@M}%
 \DN@{\ifx\next\nofrills\DN@\nofrills{\nextii@}%
      \else\DN@{\nextii@{{#1}}}\fi
      \next@}%
 \DNii@##1{%
\ifmonograph@\bgroup\else\setbox\tocbox@\vbox\bgroup
   \centerline{\headfont@\ignorespaces##1\unskip}\nobreak
   \vskip\belowheadskip \fi
   \setbox\tocheadbox@\hbox{0.\enspace}%
   \setbox\tocsubheadbox@\hbox{0.0.\enspace}%
   \leftskip\indenti \rightskip\leftskip
   \setbox\pagesbox@\hbox{\bf\quad000}\advance\rightskip\wd\pagesbox@
   \newtocdefs
 }%
 \FN@\next@}
\def\endtoc{\par\egroup}
\let\pretitle\relax
\let\preauthor\relax
\let\preaffil\relax
\let\predate\relax
\let\preabstract\relax
\let\prepaper\relax
\def\dedicatory #1\enddedicatory{\def\preabstract{{\medskip
  \eightpoint\it \raggedcenter@#1\endgraf}}}
\def\thetranslator@{}
\def\translator#1\endtranslator{\def\thetranslator@{\nobreak\medskip
 \line{\eightpoint\hfil Translated by \uppercase{#1}\qquad\qquad}\nobreak}}
\outer\def\endtopmatter{\runaway@{abstract}%
 \edef\next{\the\leftheadtoks}\ifx\next\empty
  \expandafter\leftheadtext\expandafter{\the\rightheadtoks}\fi
 \ifmonograph@\else
   \ifx\thesubjclass@\empty\else \makefootnote@{}{\thesubjclass@}\fi
   \ifx\thekeywords@\empty\else \makefootnote@{}{\thekeywords@}\fi
   \ifx\thethanks@\empty\else \makefootnote@{}{\thethanks@}\fi
 \fi
  \pretitle
  \ifmonograph@ \topskip7pc \else \topskip4pc \fi
  \box\titlebox@
  \topskip10pt
  \preauthor
  \ifvoid\authorbox@\else \vskip2.5pc plus1pc \unvbox\authorbox@\fi
  \preaffil
  \ifvoid\affilbox@\else \vskip1pc plus.5pc \unvbox\affilbox@\fi
  \predate
  \ifx\thedate@\empty\else \vskip1pc plus.5pc \line{\hfil\thedate@\hfil}\fi
  \preabstract
  \ifvoid\abstractbox@\else \vskip1.5pc plus.5pc \unvbox\abstractbox@ \fi
  \ifvoid\tocbox@\else\vskip1.5pc plus.5pc \unvbox\tocbox@\fi
  \prepaper
  \vskip2pc plus1pc
}
\def\document{\let\fontlist@\relax\let\alloclist@\relax
  \tenpoint}
\newskip\aboveheadskip       \aboveheadskip\bigskipamount
\newdimen\belowheadskip      \belowheadskip6\p@
\def\headfont@{\smc}
\def\penaltyandskip@#1#2{\relax\ifdim\lastskip<#2\relax\removelastskip
      \ifnum#1=\z@\else\penalty@#1\relax\fi\vskip#2%
  \else\ifnum#1=\z@\else\penalty@#1\relax\fi\fi}
\def\nobreak{\penalty\@M
  \ifvmode\def\penalty@{\let\penalty@\penalty\count@@@}%
  \everypar{\let\penalty@\penalty\everypar{}}\fi}
\let\penalty@\penalty
\def\heading#1\endheading{\head#1\endhead}
\def\subheading#1{\subhead#1\endsubhead}
\def\specialheadfont@{\bf}
\outer\def\specialhead{\par\penaltyandskip@{-200}\aboveheadskip
  \begingroup\interlinepenalty\@M\rightskip\z@ plus\hsize \let\\\linebreak
  \specialheadfont@\noindent\ignorespaces}
\def\endspecialhead{\par\endgroup\nobreak\vskip\belowheadskip}
\outer\def\head#1\endhead{\par\penaltyandskip@{-200}\aboveheadskip
 {\headfont@\raggedcenter@\interlinepenalty\@M
 \ignorespaces#1\endgraf}\nobreak
 \vskip\belowheadskip
 \headmark{#1}}
\let\headmark\eat@
\newskip\subheadskip       \subheadskip\medskipamount
\def\subheadfont@{\bf}
\outer\def\subhead{\nofrills@{.\enspace}\subhead@
 \DNii@##1\endsubhead{\par\penaltyandskip@{-100}\subheadskip
  \varindent@{\usualspace@{{\subheadfont@\enspace}}%
 \subheadfont@\ignorespaces##1\unskip\subhead@}\ignorespaces}%
 \FN@\next@}
\outer\def\subsubhead{\nofrills@{.\enspace}\subsubhead@
 \DNii@##1\endsubsubhead{\par\penaltyandskip@{-50}\medskipamount
      {\usualspace@{{\it\enspace}}%
  \it\ignorespaces##1\unskip\subsubhead@}\ignorespaces}%
 \FN@\next@}
\def\proclaimheadfont@{\bf}
\outer\def\proclaim{\runaway@{proclaim}\def\envir@{proclaim}%
  \nofrills@{.\enspace}\proclaim@
 \DNii@##1{\penaltyandskip@{-100}\medskipamount\varindent@
   \usualspace@{{\proclaimheadfont@\enspace}}\proclaimheadfont@
   \ignorespaces##1\unskip\proclaim@
  \sl\ignorespaces}%
 \FN@\next@}
\outer\def\endproclaim{\let\envir@\relax\par\rm
  \penaltyandskip@{55}\medskipamount}
\def\demoheadfont@{\it}
\def\demo{\runaway@{proclaim}\nofrills@{.\enspace}\demo@
     \DNii@##1{\par\penaltyandskip@\z@\medskipamount
  {\usualspace@{{\demoheadfont@\enspace}}%
  \varindent@\demoheadfont@\ignorespaces##1\unskip\demo@}\rm
  \ignorespaces}\FN@\next@}

\def\qed{\ifhmode\unskip\nobreak\fi\quad\ifmmode\square\else$\m@th\square$\fi}

\def\definition{\runaway@{proclaim}%
  \nofrills@{.\proclaimheadfont@\enspace}\definition@
        \DNii@##1{\penaltyandskip@{-100}\medskipamount
        {\usualspace@{{\proclaimheadfont@\enspace}}%
        \varindent@\proclaimheadfont@\ignorespaces##1\unskip\definition@}%
        \rm \ignorespaces}\FN@\next@}

\newdimen\rosteritemwd
\newcount\rostercount@
\newif\iffirstitem@
\let\plainitem@\item
\newtoks\everypartoks@
\def\par@{\everypartoks@\expandafter{\the\everypar}\everypar{}}
\def\roster{\edef\leftskip@{\leftskip\the\leftskip}%
 \relaxnext@
 \rostercount@\z@  
 \def\item{\FN@\rosteritem@}%
 \DN@{\ifx\next\runinitem\let\next@\nextii@\else
  \let\next@\nextiii@\fi\next@}%
 \DNii@\runinitem  
  {\unskip  
   \DN@{\ifx\next[\let\next@\nextii@\else
    \ifx\next"\let\next@\nextiii@\else\let\next@\nextiv@\fi\fi\next@}%
   \DNii@[####1]{\rostercount@####1\relax
    \enspace{\rm(\number\rostercount@)}~\ignorespaces}%
   \def\nextiii@"####1"{\enspace{\rm####1}~\ignorespaces}%
   \def\nextiv@{\enspace{\rm(1)}\rostercount@\@ne~}%
   \par@\firstitem@false  
   \FN@\next@}%
 \def\nextiii@{\par\par@  
  \penalty\@m\smallskip\vskip-\parskip
  \firstitem@true}%
 \FN@\next@}
\def\rosteritem@{\iffirstitem@\firstitem@false\else\par\vskip-\parskip\fi
 \leftskip3\parindent\noindent  
 \DNii@[##1]{\rostercount@##1\relax
  \llap{\hbox to2.5\parindent{\hss\rm(\number\rostercount@)}%
   \hskip.5\parindent}\ignorespaces}%
 \def\nextiii@"##1"{%
  \llap{\hbox to2.5\parindent{\hss\rm##1}\hskip.5\parindent}\ignorespaces}%
 \def\nextiv@{\advance\rostercount@\@ne
  \llap{\hbox to2.5\parindent{\hss\rm(\number\rostercount@)}%
   \hskip.5\parindent}}%
 \ifx\next[\let\next@\nextii@\else\ifx\next"\let\next@\nextiii@\else
  \let\next@\nextiv@\fi\fi\next@}

\newif\ifnextRunin@
\def\endroster{\relaxnext@
 \par\leftskip@  
 \penalty-50 \vskip-\parskip\smallskip  
 \DN@{\ifx\next\Runinitem\let\next@\relax
  \else\nextRunin@false\let\item\plainitem@  
   \ifx\next\par 
    \DN@\par{\everypar\expandafter{\the\everypartoks@}}%
   \else  
    \DN@{\noindent\everypar\expandafter{\the\everypartoks@}}%
  \fi\fi\next@}%
 \FN@\next@}
\newcount\rosterhangafter@
\def\Runinitem#1\roster\runinitem{\relaxnext@
 \rostercount@\z@ 
 \def\item{\FN@\rosteritem@}%
 \def\runinitem@{#1}%
 \DN@{\ifx\next[\let\next\nextii@\else\ifx\next"\let\next\nextiii@
  \else\let\next\nextiv@\fi\fi\next}%
 \DNii@[##1]{\rostercount@##1\relax
  \def\item@{{\rm(\number\rostercount@)}}\nextv@}%
 \def\nextiii@"##1"{\def\item@{{\rm##1}}\nextv@}%
 \def\nextiv@{\advance\rostercount@\@ne
  \def\item@{{\rm(\number\rostercount@)}}\nextv@}%
 \def\nextv@{\setbox\z@\vbox  
  {\ifnextRunin@\noindent\fi  
  \runinitem@\unskip\enspace\item@~\par  
  \global\rosterhangafter@\prevgraf}%
  \firstitem@false  
  \ifnextRunin@\else\par\fi  
  \hangafter\rosterhangafter@\hangindent3\parindent
  \ifnextRunin@\noindent\fi  
  \runinitem@\unskip\enspace 
  \item@~\ifnextRunin@\else\par@\fi  
  \nextRunin@true\ignorespaces}%
 \FN@\next@}
\def\footmarkform@#1{$\m@th^{#1}$}
\let\thefootnotemark\footmarkform@
\def\makefootnote@#1#2{\insert\footins
 {\interlinepenalty\interfootnotelinepenalty
 \eightpoint\splittopskip\ht\strutbox\splitmaxdepth\dp\strutbox
 \floatingpenalty\@MM\leftskip\z@\rightskip\z@\spaceskip\z@\xspaceskip\z@
 \leavevmode{#1}\footstrut\ignorespaces#2\unskip\lower\dp\strutbox
 \vbox to\dp\strutbox{}}}
\newcount\footmarkcount@
\footmarkcount@\z@
\def\footnotemark{\let\@sf\empty\relaxnext@
 \ifhmode\edef\@sf{\spacefactor\the\spacefactor}\/\fi
 \DN@{\ifx[\next\let\next@\nextii@\else
  \ifx"\next\let\next@\nextiii@\else
  \let\next@\nextiv@\fi\fi\next@}%
 \DNii@[##1]{\footmarkform@{##1}\@sf}%
 \def\nextiii@"##1"{{##1}\@sf}%
 \def\nextiv@{\iffirstchoice@\global\advance\footmarkcount@\@ne\fi
  \footmarkform@{\number\footmarkcount@}\@sf}%
 \FN@\next@}
\def\footnotetext{\relaxnext@
 \DN@{\ifx[\next\let\next@\nextii@\else
  \ifx"\next\let\next@\nextiii@\else
  \let\next@\nextiv@\fi\fi\next@}%
 \DNii@[##1]##2{\makefootnote@{\footmarkform@{##1}}{##2}}%
 \def\nextiii@"##1"##2{\makefootnote@{##1}{##2}}%
 \def\nextiv@##1{\makefootnote@{\footmarkform@{\number\footmarkcount@}}{##1}}%
 \FN@\next@}
\def\footnote{\let\@sf\empty\relaxnext@
 \ifhmode\edef\@sf{\spacefactor\the\spacefactor}\/\fi
 \DN@{\ifx[\next\let\next@\nextii@\else
  \ifx"\next\let\next@\nextiii@\else
  \let\next@\nextiv@\fi\fi\next@}%
 \DNii@[##1]##2{\footnotemark[##1]\footnotetext[##1]{##2}}%
 \def\nextiii@"##1"##2{\footnotemark"##1"\footnotetext"##1"{##2}}%
 \def\nextiv@##1{\footnotemark\footnotetext{##1}}%
 \FN@\next@}
\def\adjustfootnotemark#1{\advance\footmarkcount@#1\relax}
\def\footnoterule{\kern-3\p@
  \hrule width 5pc\kern 2.6\p@} 
\def\captionfont@{\smc}
\def\topcaption#1#2\endcaption{%
  {\dimen@\hsize \advance\dimen@-\captionwidth@
   \rm\raggedcenter@ \advance\leftskip.5\dimen@ \rightskip\leftskip
  {\captionfont@#1}%
  \if\notempty{#2}.\enspace\ignorespaces#2\fi
  \endgraf}\nobreak\bigskip}
\def\botcaption#1#2\endcaption{%
  \nobreak\bigskip
  \setboxz@h{\captionfont@#1\if\notempty{#2}.\enspace\rm#2\fi}%
  {\dimen@\hsize \advance\dimen@-\captionwidth@
   \leftskip.5\dimen@ \rightskip\leftskip
   \noindent \ifdim\wdz@>\captionwidth@ 
   \else\hfil\fi 
  {\captionfont@#1}\if\notempty{#2}.\enspace\rm#2\fi\endgraf}}
\def\@ins{\par\begingroup\def\vspace##1{\vskip##1\relax}%
  \def\captionwidth##1{\captionwidth@##1\relax}%
  \setbox\z@\vbox\bgroup} 
\def\block{\RIfMIfI@\nondmatherr@\block\fi
       \else\ifvmode\vskip\abovedisplayskip\noindent\fi
        $$\def\endblock{\par\egroup$$}\fi
  \vbox\bgroup\advance\hsize-2\indenti\noindent}
\def\endblock{\par\egroup}
\def\cite#1{{\rm[{\citefont@\m@th#1}]}}
\def\citefont@{\rm}
\def\refsfont@{\eightpoint}
\outer\def\Refs{\runaway@{proclaim}%
 \relaxnext@ \DN@{\ifx\next\nofrills\DN@\nofrills{\nextii@}\else
  \DN@{\nextii@{References}}\fi\next@}%
 \DNii@##1{\penaltyandskip@{-200}\aboveheadskip
  \line{\hfil\headfont@\ignorespaces##1\unskip\hfil}\nobreak
  \vskip\belowheadskip
  \begingroup\refsfont@\sfcode`.=\@m}%
 \FN@\next@}

\newbox\nobox@            \newbox\keybox@           \newbox\bybox@
\newbox\paperbox@         \newbox\paperinfobox@     \newbox\jourbox@
\newbox\volbox@           \newbox\issuebox@         \newbox\yrbox@
\newbox\pagesbox@         \newbox\bookbox@          \newbox\bookinfobox@
\newbox\publbox@          \newbox\publaddrbox@      \newbox\finalinfobox@
\newbox\edsbox@           \newbox\langbox@
\newif\iffirstref@        \newif\iflastref@
\newif\ifprevjour@        \newif\ifbook@            \newif\ifprevinbook@
\newif\ifquotes@          \newif\ifbookquotes@      \newif\ifpaperquotes@
\newdimen\bysamerulewd@
\setboxz@h{\refsfont@\kern3em}
\bysamerulewd@\wdz@
\newdimen\refindentwd
\setboxz@h{\refsfont@ 00. }
\refindentwd\wdz@
\outer\def\ref{\begingroup \noindent\hangindent\refindentwd
 \firstref@true \def\nofrills{\def\refkern@{\kern3sp}}%
 \ref@}
\def\ref@{\book@false \bgroup\let\endrefitem@\egroup \ignorespaces}
\def\moreref{\endrefitem@\endref@\firstref@false\ref@}%
\def\transl{\endrefitem@\endref@\firstref@false
  \book@false
  \prepunct@
  \setboxz@h\bgroup \aftergroup\unhbox\aftergroup\z@
    \def\endrefitem@{\unskip\refkern@\egroup}\ignorespaces}%
\def\emptyifempty@{\dimen@\wd\currbox@
  \advance\dimen@-\wd\z@ \advance\dimen@-.1\p@
  \ifdim\dimen@<\z@ \setbox\currbox@\copy\voidb@x \fi}
\let\refkern@\relax
\def\endrefitem@{\unskip\refkern@\egroup
  \setboxz@h{\refkern@}\emptyifempty@}\ignorespaces
\def\refdef@#1#2#3{\edef\next@{\noexpand\endrefitem@
  \let\noexpand\currbox@\csname\expandafter\eat@\string#1box@\endcsname
    \noexpand\setbox\noexpand\currbox@\hbox\bgroup}%
  \toks@\expandafter{\next@}%
  \if\notempty{#2#3}\toks@\expandafter{\the\toks@
  \def\endrefitem@{\unskip#3\refkern@\egroup
  \setboxz@h{#2#3\refkern@}\emptyifempty@}#2}\fi
  \toks@\expandafter{\the\toks@\ignorespaces}%
  \edef#1{\the\toks@}}
\refdef@\no{}{. }
\refdef@\key{[\m@th}{] }
\refdef@\by{}{}
\def\bysame{\by\hbox to\bysamerulewd@{\hrulefill}\thinspace
   \kern0sp}
\def\manyby{\message{\string\manyby is no longer necessary; \string\by
  can be used instead, starting with version 2.0 of \styname.STY}\by}
\refdef@\paper{\ifpaperquotes@``\fi\it}{}
\refdef@\paperinfo{}{}
\def\jour{\endrefitem@\let\currbox@\jourbox@
  \setbox\currbox@\hbox\bgroup
  \def\endrefitem@{\unskip\refkern@\egroup
    \setboxz@h{\refkern@}\emptyifempty@
    \ifvoid\jourbox@\else\prevjour@true\fi}%
\ignorespaces}
\refdef@\vol{\ifbook@\else\bf\fi}{}
\refdef@\issue{no. }{}
\refdef@\yr{}{}
\refdef@\pages{}{}
\def\page{\endrefitem@\def\pp@{\def\pp@{pp.~}p.~}\let\currbox@\pagesbox@
  \setbox\currbox@\hbox\bgroup\ignorespaces}
\def\pp@{pp.~}
\def\book{\endrefitem@ \let\currbox@\bookbox@
 \setbox\currbox@\hbox\bgroup\def\endrefitem@{\unskip\refkern@\egroup
  \setboxz@h{\ifbookquotes@``\fi}\emptyifempty@
  \ifvoid\bookbox@\else\book@true\fi}%
  \ifbookquotes@``\fi\it\ignorespaces}
\def\inbook{\endrefitem@
  \let\currbox@\bookbox@\setbox\currbox@\hbox\bgroup
  \def\endrefitem@{\unskip\refkern@\egroup
  \setboxz@h{\ifbookquotes@``\fi}\emptyifempty@
  \ifvoid\bookbox@\else\book@true\previnbook@true\fi}%
  \ifbookquotes@``\fi\ignorespaces}
\refdef@\eds{(}{, eds.)}
\def\ed{\endrefitem@\let\currbox@\edsbox@
 \setbox\currbox@\hbox\bgroup
 \def\endrefitem@{\unskip, ed.)\refkern@\egroup
  \setboxz@h{(, ed.)}\emptyifempty@}(\ignorespaces}
\refdef@\bookinfo{}{}
\refdef@\publ{}{}
\refdef@\publaddr{}{}
\refdef@\finalinfo{}{}
\refdef@\lang{(}{)}

\let\refdef@\relax 
\def\ppunbox@#1{\ifvoid#1\else\prepunct@\unhbox#1\fi}
\def\nocomma@#1{\ifvoid#1\else\changepunct@3\prepunct@\unhbox#1\fi}
\def\changepunct@#1{\ifnum\lastkern<3 \unkern\kern#1sp\fi}
\def\prepunct@{\count@\lastkern\unkern
  \ifnum\lastpenalty=0
    \let\penalty@\relax
  \else
    \edef\penalty@{\penalty\the\lastpenalty\relax}%
  \fi
  \unpenalty
  \let\refspace@\ \ifcase\count@,
\or;\or.\or 
  \or\let\refspace@\relax
  \else,\fi
  \ifquotes@''\quotes@false\fi \penalty@ \refspace@
}
\def\transferpenalty@#1{\dimen@\lastkern\unkern
  \ifnum\lastpenalty=0\unpenalty\let\penalty@\relax
  \else\edef\penalty@{\penalty\the\lastpenalty\relax}\unpenalty\fi
  #1\penalty@\kern\dimen@}
\def\endref{\endrefitem@\lastref@true\endref@
  \par\endgroup \prevjour@false \previnbook@false }
\def\endref@{%
\iffirstref@
  \ifvoid\nobox@\ifvoid\keybox@\indent\fi
  \else\hbox to\refindentwd{\hss\unhbox\nobox@}\fi
  \ifvoid\keybox@
  \else\ifdim\wd\keybox@>\refindentwd
         \box\keybox@
       \else\hbox to\refindentwd{\unhbox\keybox@\hfil}\fi\fi
  \kern4sp\ppunbox@\bybox@
\fi 
  \ifvoid\paperbox@
  \else\prepunct@\unhbox\paperbox@
    \ifpaperquotes@\quotes@true\fi\fi
  \ppunbox@\paperinfobox@
  \ifvoid\jourbox@
    \ifprevjour@ \nocomma@\volbox@
      \nocomma@\issuebox@
      \ifvoid\yrbox@\else\changepunct@3\prepunct@(\unhbox\yrbox@
        \transferpenalty@)\fi
      \ppunbox@\pagesbox@
    \fi 
  \else \prepunct@\unhbox\jourbox@
    \nocomma@\volbox@
    \nocomma@\issuebox@
    \ifvoid\yrbox@\else\changepunct@3\prepunct@(\unhbox\yrbox@
      \transferpenalty@)\fi
    \ppunbox@\pagesbox@
  \fi 
  \ifbook@\prepunct@\unhbox\bookbox@ \ifbookquotes@\quotes@true\fi \fi
  \nocomma@\edsbox@
  \ppunbox@\bookinfobox@
  \ifbook@\ifvoid\volbox@\else\prepunct@ vol.~\unhbox\volbox@
  \fi\fi
  \ppunbox@\publbox@ \ppunbox@\publaddrbox@
  \ifbook@ \ppunbox@\yrbox@
    \ifvoid\pagesbox@
    \else\prepunct@\pp@\unhbox\pagesbox@\fi
  \else
    \ifprevinbook@ \ppunbox@\yrbox@
      \ifvoid\pagesbox@\else\prepunct@\pp@\unhbox\pagesbox@\fi
    \fi \fi
  \ppunbox@\finalinfobox@
  \iflastref@
    \ifvoid\langbox@.\ifquotes@''\fi
    \else\changepunct@2\prepunct@\unhbox\langbox@\fi
  \else
    \ifvoid\langbox@\changepunct@1%
    \else\changepunct@3\prepunct@\unhbox\langbox@
      \changepunct@1\fi
  \fi
}
\outer\def\enddocument{%
 \runaway@{proclaim}%
\ifmonograph@ 
\else
 \nobreak
 \thetranslator@
 \count@\z@ \loop\ifnum\count@<\addresscount@\advance\count@\@ne
 \csname address\number\count@\endcsname
 \csname email\number\count@\endcsname
 \repeat
\fi
 \vfill\supereject\end}
\def\folio{{\foliofont@\ifnum\pageno<\z@ \romannumeral-\pageno
 \else\number\pageno \fi}}
\def\foliofont@{\eightrm}
\def\headlinefont@{\eightpoint}
\def\leftheadline{\rlap{\folio}\hfill \iftrue\topmark\fi \hfill}
\def\rightheadline{\hfill \expandafter
  \hfill \llap{\folio}}
\newtoks\leftheadtoks
\newtoks\rightheadtoks
\def\leftheadtext{\nofrills@{\uppercasetext@}\lht@
  \DNii@##1{\leftheadtoks\expandafter{\lht@{##1}}%
    \mark{\the\leftheadtoks\noexpand\else\the\rightheadtoks}
    \ifsyntax@\setboxz@h{\def\\{\unskip\space\ignorespaces}%
        \headlinefont@##1}\fi}%
  \FN@\next@}
\def\rightheadtext{\nofrills@{\uppercasetext@}\rht@
  \DNii@##1{\rightheadtoks\expandafter{\rht@{##1}}%
    \mark{\the\leftheadtoks\noexpand\else\the\rightheadtoks}%
    \ifsyntax@\setboxz@h{\def\\{\unskip\space\ignorespaces}%
        \headlinefont@##1}\fi}%
  \FN@\next@}
\headline={\def\chapter#1{\chapterno@. }%
  \def\\{\unskip\space\ignorespaces}\headlinefont@
  \ifodd\pageno \rightheadline \else \leftheadline\fi}
\def\NoRunningHeads{\global\runheads@false\global\let\headmark\eat@}

\def\logo@{\baselineskip2pc \hbox to\hsize{\hfil\eightpoint Typeset by
 \AmSTeX}}
\newif\iffirstpage@     \firstpage@true
\newif\ifrunheads@      \runheads@true
\output={\output@}
\def\output@{\shipout\vbox{%
 \iffirstpage@ \global\firstpage@false
  \pagebody \logo@ \makefootline%
 \else \ifrunheads@ \makeheadline \pagebody
       \else \pagebody \makefootline \fi
 \fi}%
 \advancepageno \ifnum\outputpenalty>-\@MM\else\dosupereject\fi}
\tenpoint
\catcode`\@=\active


 \NoBlackBoxes





 \define\zz{ {{\bold{Z}_2}}}

 
 \define\calb{\Cal B}
 \define\calc{\Cal C}

 \define\calq{\Cal Q}
 \define\calr{\Cal R}
 \define\calz{\Cal Z}



 \define\cycd#1#2{{\calz}_{#1}(#2)}
 \define\cyc#1#2{{\calz}^{#1}(#2)}
 \define\cych#1{{{\calz}^{#1}}}
 \define\cycp#1#2{{\calz}^{#1}(\bbp(#2))}
 \define\cyf#1#2{{\cyc{#1}{#2}}^{fix}}
 \define\cyfd#1#2{{\cycd{#1}{#2}}^{fix}}

 \define\crl#1#2{{\calz}_{\bbr}^{#1}{(#2)}}

 \define\crd#1#2{\widetilde{\calz}_{\bbr}^{#1}{(#2)}}

 \define\crld#1#2{{\calz}_{#1,\bbr}{(#2)}}
 \define\crdd#1#2{\widetilde{\calz}_{#1,{\bbr}}{(#2)}}

 \define\crlh#1{{\calz}_{\bbr}^{#1}}
 \define\crdh#1{{\widetilde{\calz}_{\bbr}^{#1}}}
 \define\cyav#1#2{{{\cyc{#1}{#2}}^{av}}}
 \define\cyavd#1#2{{\cycd{#1}{#2}}^{av}}

 \define\cyaa#1#2{{\cyc{#1}{#2}}^{-}}
 \define\cyaad#1#2{{\cycd{#1}{#2}}^{-}}

 \define\cyq#1#2{{\calq}^{#1}(#2)}
 \define\cyqd#1#2{{\calq}_{#1}(#2)}

 \define\cqt#1#2{{\calz}_{\bbh}^{#1}{(#2)}}
 \define\cqtav#1#2{{\calz}^{#1}{(#2)}^{av}}
 \define\cqtrd#1#2{\widetilde{\calz}_{\bbh}^{#1}{(#2)}}

 \define\cyct#1#2{{\calz}^{#1}(#2)_\zz}
 \define\cyft#1#2{{\cyc{#1}{#2}}^{fix}_\zz}
\define\cxg#1#2{G^{#1}_{\bbc}(#2)}
 \define\reg#1#2{G^{#1}_{\bbr}(#2)}

 \define\cyaat#1#2{{\cyc{#1}{#2}}^{-}_\zz}


 \define\fflag#1#2{{#1}={#1}_{#2} \supset {#1}_{{#2}-1} \supset
 \ldots \supset {#1}_{0} }
 
 \define\vect#1{ {\Cal{V}ect}_{#1}}


 \define\chv#1#2#3{{\calc}^{#1}_{#2}(#3)}
 \define\chvd#1#2#3{{\calc}_{#1,#2}(#3)}
 \define\chm#1#2{{\calc}_{#1}(#2)}



 \define\Claim#1{\subheading{Claim #1}}

 \define\xrightarrow#1{\overset{#1}\to{\rightarrow}}


\hfuzz1pc 


\define\bbz{\Bbb Z}

\define\bbr{\Bbb R}
\define\bbc{\Bbb C}
\define\bbh{\Bbb H}
\define\bbp{\Bbb P}

\define\bbs{\Bbb S}

\define\bm{\bold M}

\define\cu{\Cal U}

\define\cc{\Cal C}
\define\cd{\Cal D}
\define\ce{\Cal E}

\define\cf{\Cal F}

\define\cs{\Cal S}
\define\cz{\Cal Z}
\define\co#1{\Cal O_{#1}}
\define\ct{\Cal T}
\define\ci{\Cal I}
\define\cR{\Cal R}


\define\a{\alpha}
\redefine\b{\beta}
\define\g{\gamma}
\redefine\d{\delta}
\define\r{\rho}
\define\s{\sigma}
\define\z{\zeta}
\define\x{\xi}

\define\e{\epsilon}
\redefine\D{\Delta}
\define\G{\Gamma}

\define\p#1{{\bbp}^{#1}}

\define\equdef{\overset\text{def}\to=}

\define\blbx{\hfill  $\square$}
\redefine\qed{\blbx}

\define\pf{\subheading{Proof}}
\define\Lemma#1{\subheading{Lemma #1}}
\define\Theorem#1{\subheading{Theorem #1}}
\define\Prop#1{\subheading{Proposition #1}}
\define\Cor#1{\subheading{Corollary #1}}
\define\Note#1{\subheading{Note #1}}
\define\Def#1{\subheading{Definition #1}}
\define\Remark#1{\subheading{Remark #1}}
\define\Ex#1{\subheading{Example #1}}
\define\arr{\longrightarrow}

\define\Hom{\operatorname{Hom}}

\redefine\Xi{X_{\infty}}

\define\jac#1#2{\left(\!\!\!\left(
\frac{\partial #1}{\partial #2}
\right)\!\!\!\right)}
\define\restrict#1{\left. #1 \right|_{t_{p+1} = \dots = t_n = 0}}

\define\SP#1#2{{\roman SP}^{#1}(#2)}

\define\coc#1#2#3{\cc^{#1}(#2;\, #3)}
\define\zoc#1#2#3{\cz^{#1}(#2;\, #3)}
\define\zyc#1#2#3{\cz^{#1}(#2 \times #3)}

\define\Div{\roman{ Div}}
\define\ar#1{\overset{#1}\to{\longrightarrow}}

\define\th#1#2{{\Bbb H}^{#1}(#2)}
\define\hth#1#2{\widehat{\Bbb H}^{#1}(#2)}


\define\bad#1#2{\cf_{#2}(#1)}

\define\pch#1{\bbp_{\bbc}(\bbh^{#1})}

\def\l{\ell}

\def\BP{{\Bbb P}}

\def\BR{{\Bbb R}}

\def\BZ{{\Bbb Z}}

\def\I#1#2{I_{#1, #2}}

\def\D{D}

\def\L{L}

\def\z2t{\text{$\bbz_2\ct$}}

\def\hH{\widehat{I\!\!H}}
\def\hh{\widehat{H}}
\def\oper{\operatorname}

\def\hHc{\hH_{\oper{cpt}}}

\def\Lloc{L^1_{\oper{loc}}}

\def\cpt{\oper{cpt}}
\def\bbth{T\!\!\! I}
\def\bbthh{\widehat{\bbth}}

\def\mini{\oper{min}}

\def\<{\left<}
\def\>{\right>}
\def\[{\left[}
\def\]{\right]}

\def\surrightarrow{\longrightarrow\!\!\!\!\!\!\!\!\longrightarrow}
\def\surleftarrow{\longleftarrow\!\!\!\!\!\!\!\!\longleftarrow}

\def\wt{\widetilde}
\def\vf{\varphi}

\def\DD{I\!\! D}
\def\supp{\operatorname{supp \ }}

\redefine\and{\qquad\text{and}\qquad}

\document


\define\Q{Q}               
\redefine\bbr{{\bold R}}
\redefine\bbz{{\bold Z}}

\centerline{\bf \titfont  The de Rham-Federer Theory of
Differential Characters}
\bigskip
\centerline{\bf \titfont and  Character  Duality }

\bigskip
\centerline{by}
\bigskip

\centerline{\bf  \aufont   Reese Harvey, \ 
 Blaine Lawson,\ and\ 
John Zweck\footnote{\footfont Research
of all authors was partially supported by the NSF}}
           
\vskip .3in
 
\centerline{\bf Abstract} \medskip
  \font\abstractfont=cmr10 at 9 pt

{{\parindent=.4in\narrower\abstractfont \noindent

In the first part of this paper the theory of differential characters
is  developed completely from a de Rham - Federer viewpoint. 
Characters are defined as equivalence classes of special currents,
called sparks, which appear naturally in the theory of singular
connections.   As in de Rham - Federer cohomology, there are many
different spaces of currents which yield the character groups.  The
fundamental exact sequences in the theory are easily derived from
methods of geometric measure theory. A multiplication of de
Rham-Federer characters is defined using transversality results for
flat and rectifiable currents established in the appendix.    It is
shown that there is a natural equivalence of ring functors from  de
Rham - Federer characters to the classical Cheeger-Simons characters
given, as in de Rham cohomology, via integration. This discussion 
rounds out the approach to differential character theory introduced
by Gillet-Soul\'e and Harris.

The groups of differential characters have an obvious topology and
natural {\sl smooth Pontrjagin duals} (introduced here). It is shown
that the dual  groups sit in two  exact sequences which resemble
the fundamental exact sequences for the character groups themselves.
They are essentially the smooth duals of the fundamental sequences 
with roles interchanged.

A principal result here is the formulation and proof of
 duality for characters on oriented
manifolds. It is shown that the pairing 
$$
(a,b)\ \mapsto\  a*b([X])
$$
given by multiplication and evaluation on the fundamental cycle,
gives an isomorphism of the group of differential characters of
degree $k$ with the dual  to characters in degree $n-k-1$ where n =
dim($X$).  A number of examples are examined in detail.

It is also shown that there are natural Thom homomorphisms
for differential characters.  In fact, given an orthogonal
connection on an oriented vector bundle $\pi:E\to X$, there is
a canonical {\sl Thom spark} $\gamma$ on $E$ which generates the Thom
image as a free module over the characters on $X$.  One has $d\g =
\tau-[X]$ where $\tau$ is a canonical Thom form on $E$. Under the
two basic differentials in the theory, the Thom map becomes the usual
ones  in de Rham theory and in integral cohomology. 
Gysin maps are also defined for differential characters.

Many examples including Morse sparks, Hodge sparks and
characteristic sparks are examined in detail.
   
}}
   
\vfill\eject

\vskip .3in               
               
\centerline{Table of Contents}
\medskip   
       
\hskip 1in  0.  Introduction
 
\hskip 1in  1.  De Rham - Federer Characters

\hskip 1in  2.  Isomorphisms of De Rham - Federer Type  
  
\hskip 1in  3.  Ring Structure  

\hskip 1in  4.  Cheeger-Simons Characters

\hskip 1in  5.  Dual Sequences

\hskip 1in  6.  Duality for Differential Characters 
 
\hskip 1in  7.  Examples of Duality
 
\hskip 1in  8.  Characters with Compact Support

\hskip 1in  9.  The Thom Homomorphism 

\hskip 1in   \!\! 10.  Gysin Maps

\hskip 1in   \!\! 11.   Morse Sparks

\hskip 1in   \!\! 12.   Hodge Sparks and Riemannian Abel-Jacobi
Mappings

\hskip 1in   \!\! 13.    Characteristic Sparks
and Degeneracy Sparks

\hskip 1in   \!\! 14.   Characters in Degree Zero.

\hskip 1in   \!\! 15.   Characters in Degree One.

\hskip 1in   \!\! 16.   Characters in Degree Two (Gerbes).

\hskip 1in   \!\! 17.   Some Historical Remarks on the Complex
Analogue

\medskip
              
\hskip 1.2in   Appendix A.  Slicing Currents by Sections of a 
Vector Bundle.

\hskip 1.2in   Appendix B.  Intersecting Currents

\hskip 1.2in   Appendix C.  The de Rham-Federer approach to integral
cohomology

\hskip 1.2in   Appendix D.  The Spark Product

\vskip .4in

\font\footfont = cmr10 at 8pt

\subheading{\S0. Introduction}   In 1973 Jim Simons and Jeff
Cheeger introduced the rings of differential characters and used
them, among other things, to give an important generalization of
the Chern-Weil homomorphism [S], [CS].  Their theory constituted a
vast generalization of the Chern-Simons invariant, and was applied
to prove striking global non-conformal-immersion theorems in
riemannian geometry.  The theory also gave rise to new invariants
for flat connections, and provided an insightful derivation of the
whole family of foliation invariants due to Godbillon-Vey,
Heitsch, Bott and others.

The first part of this paper presents a self-contained
formulation of the theory of differential characters based on special
currents, called {\bf sparks}, which appear in the study of singular
connections.   This approach can be thought of as a {\bf de Rham -
Federer theory of differential characters} in analogy with de Rham -
Federer cohomology theories. Such formulations are originally due
to Gillet-Soul\'e [GS$_1$] and B. Harris [H].  Our discussion 
presents, rounds out, and unifies their work.

It is of course always useful to have distinct representations
of the same theory, and there are many advantages of the de Rham -
Federer viewpoint.
One sees, for example, that the differential character groups carry 
a natural topology.  One also finds that in many ways differential
characters behave like a homology-cohomology theory.  A striking
instance is that differential characters obey a Poincar\'e-type
{\bf  Duality} (a main result of this paper) ---
that is, on an oriented $n$-manifold, taking   $*$-product and
evaluating on the fundamental cycle identifies the 
compactly supported  characters of degree $k$ with 
the smooth Pontrjagin dual of the group of characters of degree
$n-k-1$.  
That one should consider Pontrjagin duals when dealing with
differential characters is a key point in this result.
The theory of differential characters also carries 
Thom homomorphisms, Gysin mappings, and certain relative long
exact sequences.  There is also an {\sl Alexander-Lefschetz-type 
Duality Theorem}. (This and the long exact sequences will be treated 
in a separate paper.)

A {\bf spark} on a manifold $X$ is a de Rham current $a$ whose
exterior derivative can be decomposed as  
$$
da=\phi-R
\tag{0.1}
$$
where $\phi$ is a smooth form  and $R$ is  integrally
flat.\footnote{{\ftfont A current is integrally flat if it can be
written as R+dS where R and S are rectifiable.}}  Such objects arise
quite naturally.  Suppose for example that $E\to X$ is  a smooth
oriented vector bundle with connection, and $\a$ is a section with
non-degenerate zeros $\Div(\a)$. Pulling down the classical Chern
transgression gives a canonical $L^1_{\text{loc}}$-form $\tau$ with 
$d\tau = \chi(\Omega)- \Div(\a)$ where $\chi(\Omega)$ is the  
Weil polynomial in the curvature which represents the Euler class of
$E$. Quite general formulas of this type are derived systematically 
 in [HL$_{1-4}$].

  A basic fact is that the decomposition (0.1)
is {\bf unique}.  Consequently,   every spark has two
differentials: $d_1a = \phi$ and $d_2a = R$. Two sparks $a_1$
and $a_2$ are defined to be equivalent  if $a_1 =
a_2+db+S$ where $b$ is any current and $S$ is integrally flat. The
spark equivalence classes define the space of characters $\hH^*(X)$,
and the  differentials $d_1$ and $d_2$ yield the two fundamental
exact sequences   $$
0\ \arr \ \hH^k_\infty(X)\arr\hH^k(X)
\ @>{\d_2}>> H^{k+1}(X,\bbz) \ \arr\ 0.
$$
$$
0\ \arr \ H^k(X,S^1)\arr \hH^k(X) \
@>{\d_1}>> Z_0^{k+1}(X)\ \arr\ 0.
$$ 
in the theory.  (See 1.11 and 1.12.)

In the classical case of de Rham - Federer cohomology there are many
spaces of currents which represent the same
functor.\footnote{{\ftfont  For example,  H*(X;\, $\bbr$) can be
represented as the cohomology of all currents, or  of normal currents,
or of flat currents, or of smooth currents. Similarly, 
 H*(X;\, $\bbz$) can be represented by integral currents, or by
integrally flat currents, or by the currents induced by smooth
singular chains, etc..}} Analogously we show in \S 2 that there are
many spaces of sparks which under suitable equivalence yield
differential characters. For example, one can restrict $a$ to be
flat and $R$ to be rectifiable, or one can restrict $a$ to be a form
with $L^1_{\text{loc}}$-coefficients and $R$ to be a smooth singular
chain.  This flexibility is quite useful in practice.

In \S 3  we introduce a product, denoted $*$,  on differential
characters using a basic ``exterior/intersection''  product on
currents.  The generic existence of this exterior/intersection
product itself constitutes a new result which is essentially a
{\bf transversality theorem for flat and rectifiable currents}. The
proof, which combines certain geometric arguments with Federer
slicing theory, is given in the appendix.  

We remark that defining products at  the chain/cochain-level
in topology is always problematical, 
and  the intersection theory for currents in Appendix B should have
considerable  independent interest.

In \S 4 we exhibit an equivalence of  functors from deRham-Federer
characters to the classical Cheeger-Simons theory.  As in the
cohomology analogue, this equivalence is essentially given by
integration over cycles (or ``spark
holonomy'').   It is then shown that our intersection $*$-product 
agrees with the  product introduced by Cheeger on differential
characters [C].  Thus the equivalence  is an equivalence of
ring functors.

From our point of view it is natural to consider the {\bf smooth
characters} namely those which are representable by  sparks which
are smooth forms.  The space of smooth characters is exactly the
kernel of the differential $\d_2$. The characters also naturally 
inherit a topology, and this ideal of smooth characters is exactly
the connected component of 0 in this topology.

The form $\phi$ in equation (0.1) could be considered the 
``curvature'' of the spark $a$, and in this sense the kernel 
of $\d_1$ consists of the ``curvature-flat'' sparks.

One of the main new results here is the establishment of a
{\bf Poincar\'e-Pontrjagin duality} for differential characters.  We
preface this result by showing that taking  Pontrjagin duals carries
the two fundamental exact sequences above into strictly analogous
sequences with roles reversed.  We then prove that if $X$ is compact
and oriented, the pairing $$
\hH^{n-k-1}(X) \times \hH^k(X) \longrightarrow S^1
$$
given by 
$$
(a,b) \mapsto  a*b([X])
$$
is non-degenerate.  In fact, the resulting mapping
$$
\hH^{n-k-1}(X) \arr \Hom\left(\hH^k(X),\, S^1\right)
$$
into continuous   group homomorphisms has dense range consisting
exactly of the {\sl smooth} homomorphisms (See \S 6).  In \S 7 we
exhibit a number of examples of this duality on standard spaces.

In \S 8 we introduce characters with compact support.  This
allows us to extend duality to non-compact manifolds. 
There are in fact further extensions.  On
manifolds with boundary there is a Lefschetz-Pontrjagin Duality
Theorem [HL$_6$]. There is also a general result of Alexander
duality type.

The sparks studied here  appear often in topology and
geometry. They play a fundamental role in Morse Theory (See \S 11.)
They also appear in the Hodge decomposition of rectifiable
cycles on a riemannian manifold. From this we introduce  a notion
of  linear equivalence on rectifiable cycles and define  riemannian
 ``Abel-Jacobi mappings'' which strictly generalize the 
well-known maps from  K\"ahler geometry. (See \S 12.)
On 3-manifolds these maps appear in the work of Chatterjee [Cha] and 
in Hitchin's discussion of  special Lagrangian cycles and the mirror
symmetry conjecture [Hi].

Sparks are also central in the theory of singular
connections developed in [HL$_{1-4}$].  In this theory they appear 
as canonical transgression forms $T$ with $dT=\Phi-\Sigma$ where
$\Phi$ is a smooth ``geometric'' form (typically a Chern-Weil
characteristic form) and $\Sigma$ is an integrally flat cycle
(typically defined by the singularities of a mapping). 

 For example, on the total space of an oriented
bundle $E\to X$ with orthogonal connection, there is a canonical
{\bf Thom spark} $\gamma$ with $d\gamma = \tau- [X]$ where
$\tau$ is the canonical Thom form [HL$_{1}$, IV] and $[X]$ is
the 0-section in $E$. Taking $*$-product with $\gamma$ defines
a Thom homomorphism for characters  with the
property that $\d_1$ and $\d_2$ transform it to the classical Thom
mappings. (See \S 9.)  This leads to  general Gysin maps in the
theory. (See \S 10.)

For a complex bundle $F\to X$ with unitary connection there are
{\sl Chern sparks} $C_k$ with  
$$
dC_k=c_k(\Omega) - \DD_k
$$
on $X$ where $c_k(\Omega)$ is the $k$th Chern characteristic form
and  $\DD_k$ is the linear dependency current of a set of
cross-sections of $F$. For a real bundle $F$ with orthogonal
connection there are  similar {\sl Pontrjagin-Ronga sparks}
associated to rank-2 degeneracies, and if $F$ is oriented, there are
{\sl Euler sparks} $X_{\a}$ associated to any atomic section
$\a$ and satisfying: $d X_{\a} = \chi(\Omega)-\Div(\a)$.  More
generally to a bundle map $\a:E\to F$ there are {\sl Thom-Porteous
sparks} $T_k$ with $dT_k=
TP_k(\Omega^F,\Omega^E)-[\Sigma_k(\a)]$ where $\Sigma_k(\a)$ is the 
$k$th degeneracy locus of $\a$ and $TP$ is the $k$th Thom-Porteous
class expressed canonically in the curvatures of $E$ and $F$.
In each of these cases the  character represented by the spark
is independent  of the sections (or bundle maps) involved. It is
exactly the ``characteristic character'' defined by Cheeger-Simons
in  [CS]. This is elaborated in \S 13. 

In low degrees differential characters have special
interpretations, and so do the underlying  sparks. 
A character of degree 0 corresponds to a smooth mapping
$f:X\to S^1$ and its associated sparks are the
functions $\tilde f:X\to \bbr$ in $\Lloc(X)$ with  $\tilde f \equiv
f$ (mod $\bbz$).  A character of
degree 1 corresponds to a complex line bundle $E$ with
unitary connection up to gauge equivalence, and underlying sparks
are given by pairs  $(E,\s)$ where 
$\s$ is an atomic section of $E$. In degree 2 sparks are related to
``gerbes with connection''.  This is discussed in \S16
and will be treated in greater detail in a forthcoming paper
[HL$_7$].

In this sequel [HL$_7$] we shall present an expanded approach
to differential characters which involves a double 
\v Cech-deRhamCurrent complex  of {\bf hypersparks}. The theory is
strikingly parallel to the one presented here (and, in fact,
contains it).  There is an equation analogous to the spark equation
(0.1), and there is an equivalence relation, leading to characters,
which generalizes the one given here. The sparks are contained in
the hypersparks as a  ``full-subtheory'' yielding an isomorphism at
the level of characters.  However, there are other interesting
full-subtheories.  For example, the subgroup of {\bf smooth
hypersparks} which are directly related to gerbes with connection
([Hi], [Br]). 
 
Throughout this paper all manifolds are assumed to be oriented.

The authors would like to thank the referee for a number of
helpful comments and in particular for suggesting that we give an
expanded  discussion of gerbes. We are also grateful to Chris Earles
for pointing out a much better version of Proposition 1.16.

The second author would like to express gratitude to IHES and the
Clay Mathematics Institute for partial support during the 
preparation of this work.


     \vskip.4in


\subheading{\S1. De Rham - Federer characters} In this
section we introduce  and examine certain equivalence classes of
special  currents called sparks.   We shall see later that our set
of classes is naturally isomorphic to the Cheeger-Simons space of
differential characters.

Let $X$ be a smooth oriented manifold of dimension $n$ (not
necessarily compact).  We adopt the following notation as in 
[deR].  Let $\ce^k(X)$ denote the space of smooth differential 
$k$-forms on $X$ with the $C^{\infty}$-topology.  By $\cd^k(X)$ we
mean the smooth $k$-forms with compact support on $X$. By the space
of {\bf currents of degree $k$ (and dimension  $n-k$)} on $X$ we
mean the topological dual space $$
{\cd'}^k(X) \  \equiv \
{\cd'}_{n-k}(X)\ \equdef  \left\{{\cd}^{n-k}(X)\right\}'.
$$
Note that ${\cd'}^*(X)$ is a complex under $d'$, the dual
of $d$.  The forms $\ce^*(X)$ embed in the currents 
${\cd'}^*(X)$ by associating to a $k$-form $\varphi$ the
degree-$k$ current, also denoted by $\varphi$, given by
$\varphi(\psi)=\int\varphi\wedge\psi$ for $\psi \in
\cd^{n-k}(X)$.  On this subspace $d' = (-1)^{n-k+1}d$, and
thereby $d$ extends to all of ${\cd'}^k(X)$. The following
subspaces of currents are defined in Federer [F].
$$ 
\aligned
\ce^k_{\Lloc}(X) \ &= \ \text{the $k$-forms with
$\Lloc$-coefficients on $X$}\\
\cR^k(X) \ &= \ \text{the
locally {\bf rectifiable} currents of degree $k$  (dimension
$(n-k)$) on $X$}\\
\cf^k(X) \ &= \ \text{the locally  {\bf
flat} currents of degree $k$ on $X$}\\
\ci\cf^k(X) \ &= \
\text{the locally {\bf integrally flat} currents of degree
$k$ on $X$}
\endaligned
$$
In what follows the word  {\sl locally} will often be
supressed.  We recall that:
$$
\aligned \cf^k(X) \ &= \ \{ \phi + d\psi \,:\,  \phi \in
\ce^k_{\Lloc}(X) \ \ \text{and}\ \ \psi \in
\ce^{k-1}_{\Lloc}(X)\}\\
\ci\cf^k(X) \ &= \ \{ R + dS
\,:\,  R \in \cR^k(X) \ \ \text{and}\ \ S \in \cR^{k-1}(X)\}
\endaligned \tag{1.1}
$$
We now introduce a space of currents which appears
naturally in geometry and analysis (cf.  [BC$_{1,2}$],
[GS$_{1-4}$], [BGS], [HL$_{1-4}$], [HZ$_{1,2}$]).

\Def{1.2} The space of {\bf sparks} of degree~$k$ on $X$ is
defined to be
$$
\cs^k(X) \ \equiv \ \{a\in {\cd'}^k(X) \ :\  da =
\phi - R \  \text{where}  \ \phi \in
\ce^{k+1}(X) \ \text{and}  \ R \in \ci\cf^{k+1}(X)\}
$$

\Lemma{1.3} {\sl For $a\in\cs^k(X)$ the decomposition of  $da$
into smooth plus integrally flat is unique, i.e., if $da =
\phi-R =\phi'-R'$ where $\phi,\phi'$ are smooth and $R,R'$
are integrally flat, then $\phi =\phi'$ and $R=R'$. }

\pf It suffices to show that
$$
\ce^{k+1}(X) \cap \ci\cf^{k+1}(X) \ = \ \{0\}, \tag{1.4}
$$
i.e., if a current $T$ is both smooth and integrally flat,
then $T=0$. If $T$ is smooth, then it has finite mass.  By
[F, 4.2.16(3), p\. 413] any integrally flat current of finite mass
is rectifiable.  Thus $T$ is both smooth and rectifiable.  Since
$T$ is smooth, its $(n-k-1)$-density $\Theta(\|T\|, x) = 0$
for all $x\in X$.   However, for a rectifiable current $T$,
one has  $\Theta(\|T\|, x) \geq 1$ for  $\|T\|$-almost all
$x$ (cf. [F, 4.1.28, p\. 384]). Thus $\|T\|$ and hence $T$ vanish.
\qed

For any spark  $a\in\cs^k(X)$  with decomposition $da=\phi-R$ as
in Lemma 1.3 we define  
$$
d_1a = \phi
\qquad\text{and}\qquad
d_2a = R.
$$

\Lemma{1.5}  {\sl  Let $a\in\cs^k(X)$ be a spark.  Then $d_1a=\phi$ 
is a smooth $d$-closed differential form with integral periods and
 $d_2a=R$ is an integrally flat cycle. }

\pf  Since $da=\phi-R$, the current $d\phi = dR$ is both smooth and
integrally flat. Hence,
$$
d\phi = dR=0
$$
by (1.4).  Since $R$ has integral periods and $\phi$ is
cohomologous to $R$, $\phi$ also has integral periods.  \qed

Motivated by the comparison formula for euler sparks [HZ$_2$]
(and, of course, Cheeger-Simons) we define an equivalence
relation on the space of sparks.

\medskip 
\Def{1.6} For each integer $k$, $0\leq k\leq n$ we define
the  group of {\bf de Rham - Federer characters} of degree
$k$ to be the  quotient
$$
\hH^k(X) \ \equdef\ \cs^k(X)/ \left\{ d{\cd'}^{k-1}(X)+
\ci\cf^k(X)\right\}
$$
The equivalence class  in $\hH^k(X)$ of a spark $a\in\cs^k(X)$
will  be denoted by $\<a\>$ or $\widehat a$. We also set 
$$
\hH^{-1}(X)\ \equdef\ \bbz
$$
with generator ${\bold I}$ satisfying $d_1{\bold I}=1$ and 
$d_2{\bold I}= [X]=1$. Note that (1.4) fails for $k=-1$.

\medskip 
Suppose $a$ and $\bar a$ are equivalent sparks; i.e. $\bar
a= a+db-S$, with $b\in{\cd'}^{k-1}(X)$ and $S\in\ci\cf^k(X)$. If
$da=\phi-R$ then
$$
d\bar a = \phi - (R + dS).\tag{1.7}
$$
Let $\cz_0^{k+1}(X)$ denote the lattice of smooth $d$-closed,
degree $k+1$ forms on $X$ with integral periods. Also, given a
cycle $R\in\cR^{k+1}(X)$, let $[R]$ denote its class in
$H^{k+1}(X,\bbz)$.

\Lemma{1.8}  {\sl The maps
$$
\align
\d_1:  &\hH^k(X)\to\cz_0^{k+1}(X)\quad\text{and}\\
\d_2:  &\hH^k(X)\to H^{k+1}(X,\bbz)
\endalign
$$
given by $\d_1(\<a\>)\equdef\phi$ and $\d_2(\<a\>)\equdef[R]$
are well defined surjections.}

\pf  Lemma 1.5 and (1.7) imply that $\d_1$ and $\d_2$ are
well defined.  Surjectivity is a direct consequence of the fundamental
isomorphisms
$$
H^*({\cd'}^*(X))\ \cong\ H^*(X;\,\bbr)
\qquad\text{and}\qquad
H^*({\ci\cf}^*(X))\ \cong\ H^*(X;\,\bbz)
\tag{1.9}
$$
due to de Rham [deR] and Federer  (cf. [HZ$_1$]) respectively.
\qed

\medskip 
\Def{1.10}  A differential character is said to be {\bf
smooth} if it can be represented by a smooth differential
form.  That is, the space $\hH^k_\infty(X)$ of {\bf smooth
characters} of degree $k$ on $X$  is the
image of the smooth forms $\ce^k(X)$ in $\hH^k(X)$. Note the natural
isomorphism
$$
\hH^k_\infty(X) \ \cong\ \ce^k(X)/\cz_0^{k}(X)
$$
To see this suppose $a \in \ce^k(X)$ represents the 0-character.
Then $a=db+S$ where $S$ is integrally flat. Thus $b$ is a spark,
and so $a\in\cz_0^{k}(X)$ by Lemma 1.5.

\Prop{1.11}  {\sl There is a  short exact sequence}
$$
0\ \arr \ \hH^k_\infty(X)\arr\hH^k(X)
\ @>{\d_2}>> H^{k+1}(X,\bbz) \ \arr\ 0.
\tag {A}
$$

\pf  If $\<a\>\in\hH^k(X)$ is represented by a smooth spark
$a\in\ce^k(X)$, then $da$ is smooth and so 
$\d_2(\<a\>)=0$. 
Conversely,
suppose $da=\phi-R$ as in 1.5 and $[R]=0$ in $H^k(X,\bbz)$. Then
$R=dS$ with $S\in\ci\cf^k(X)$ and hence $d(a+S)=\phi$. This
proves that $[\phi]\in H^{k+1}(X,\bbr)$ vanishes. Therefore,
by de Rham there exist a smooth $\psi\in\ce^k(X)$ with $d\psi=\phi$.
Since $a+S-\psi$ is $d$-closed, the class of $a+S-\psi$ in
$H^k(X,\bbr)$ is represented by a smooth $d$-closed form
$\eta\in\ce^k(X)$. That is, $a+S-\psi=\eta+db$ for some
$b\in{\cd'}^{k-1}(X)$. The equivalent spark $\bar
a=a+S-db=\psi+\eta$ representing $\<a\>$ is smooth. This
proves that the sequence is exact. \qed

\smallskip\ 
\Prop{1.12}  {\sl There is a  short exact sequence}
$$
0\ \arr \ H^k(X,S^1)\arr \hH^k(X) \
@>{\d_1}>> Z_0^{k+1}(X)\ \arr\ 0.
\tag {B}
$$

\pf In keeping with our de Rham - Federer approach we give the
following representation of $H^*(X;\, S^1)$ in terms of currents:
$$
H^k(X, S^1) \ \cong\ \frac{\{a\in{\cd'}^k(X):da\in \ci\cf^{k+1}(X)\}}
{d{\cd'}^{k-1}(X)+\ci\cf^k(X)}.\tag1.13
$$
To prove this,  consider the commutative diagram of
sheaves on $X$:
$$
\align
0\arr  &\bbz\arr \ci\cf^0\ @>{d }>>\ \ci\cf^1\ @>{d}>>\ \ci\cf^2\
@>{d }>>\ \cdots\ .\\
&\downarrow  \hskip25pt\downarrow\hskip35pt\downarrow
\hskip35pt\downarrow\\
0\arr &\bbr\arr{\cd'}^0\ @>{d }>>\ {\cd'}^1\ @>{d }>>\ 
{\cd'}^2\ @>{d }>>\ \cdots\ .
\endalign
$$
The rows are acyclic resolutions of the constant sheaves $\bbz$ and
$\bbr$ respectively (cf. [HZ$_1$]).  (This fact establishes the
isomorphisms (1.9) above.)  \ \ Taking the quotient we deduce from
general arguments  that 
$$
0\arr\bbr/\bbz\arr{\cd'}^0/\ci\cf^0@>{d }>>\ \cdots\ 
$$
is an acyclic resolution of the constant sheaf $S^1=\bbr/\bbz$,
and therefore (1.13) holds.

To prove Proposition 1.12 it suffices by Lemma 1.8 to show
that ker$(\d_1) \cong H^1(X,S^1)$. 
Since $\d_1(\<a\>)=\phi=0$ if and only if $da=-R\in
\ci\cf^{k+1}(X)$,  equation (1.13) shows that kernel 
$\d_1=H^k(X,S^1)$ as desired. \qed

\medskip
Consider the direct sum of the maps $\d_1$ and $\d_2$ as a
single map $\d$. Define 
$$
\Q^k(X)\  = \ \{(\phi,u)\in Z^k_0(X)
\times H^k(X,\bbz): [\phi] = u\otimes\bbr\ \text{ in }\ 
H^k(X,\bbr)\}.
$$
Following Cheeger-Simons this should be denoted $R^k(X)$. We have
changed the notation  to avoid confusion with the space of rectifiable
currents.

\Prop{1.14}  {\sl There is a  short exact sequence}
$$
0\ \arr \ \frac{H^k(X;\,\bbr)}{H^k_{\oper{free}}
(X;\, \bbz)} \ \arr\ \hH^k(X) \
@>{\d}>> \Q^{k+1}(X) \ \arr\ 0,
$$

\pf  To see that $\d$ is surjective consider a pair $(\phi, u)
\in\Q^{k+1}(X)$ and choose an
integrally flat cycle $R$ with $[R]=u$. Then $\phi - R$ is
a  $d$-closed current whose de Rham cohomology class is
zero.  By [deR] the cohomology of currents agrees with the
real standard cohomology of $X$. Hence, there exists a current
$a$ with $da = \phi -R$, and surjectivity is proved.  Suppose
now that $a\in\cs^k(X)$ and $\d(\<a\>) = 0$.  Then  $da =
dS$ for some integrally flat current $S$. Let $a$ denote the
equivalent $d$-closed spark $\widetilde a=a-S$. Suppose $a'$
is another $d$-closed spark in the spark class $\<\widetilde
a\>\in\hH^k(X)$. Then $a'=\widetilde a+db+S_0$ with $S_0$
integrally flat. Since $a$ and $a'$ are $d$-closed, $S_0$
must be a cycle. Thus we see that $[\widetilde a]\in
H^k(X,\bbr)$ is determined by the class $\<a\>\in\hH^k(X)$
exactly up to translations by the subgroup
$H^k_{\text{free}}(X;\, \bbz)$ generated by the integrally
flat cycles. \qed

\Remark{1.15}  Via 1.14  the group of
differential  characters carries  a  {\bf natural topology} 
coming from  the  $C^{\infty}$-topology on forms in
$\cz_0^{k+1}(X)$  and the standard topology on the torus 
${H^k(X;\, \bbr)}/{H^k_{\oper{free}}(X;\, \bbz)}$.  In this
topology we see that $\hH^k_{\infty}(X)$ is the {\sl connected
component of 0}, and the sequence (1.11) can be reinterpreted as
the sequence obtained by dividing by this component.

\medskip

The decompositions (A) and (B) given above can be further refined 
by restricting $\d_1$ to kernel $\d_2=\hH^k_\infty(X)$ and 
 $\d_2$ to kernel $\d_1\equiv H^k(X,S^1)$ and then also
decomposing the images. This leads to a large diagram where
(A) and (B) are the central column and row respectively.

\Prop{1.16} {\sl The following diagram commutes, and each row and
column is exact} 
$$
\CD
\ @.  0 @.  0 @.  0 @.  \ \\
@.    @VVV    @VVV    @VVV @. \\
0  @>>>  \frac{H^k(X,\bbr)}{H^k_{\oper{free}}
(X, \bbz)}   @>>>  \hH^k_\infty(X) 
@>{\d_1}>> d\ce^k(X)   @>>>   0  \\
@.    @VVV    @VVV    @VVV @. \\
0  @>>> H^k(X,S^1)  @>>>  \hH^k(X) 
@>{\d_1}>> \cz^{k+1}_0(X)   @>>>   0  \\
@.    @V{\d_2}VV    @V{\d_2}VV    @VVV @. \\
0  @>>>  H^{k+1}_{\text{tor}}(X,\bbz)  @>>> H^{k+1}(X,\bbz)
@>{}>>  H^{k+1}_{\text{free}}(X,\bbz)  @>>>   0  \\
@.    @VVV    @VVV    @VVV @. \\
\ @.  0 @.  0 @.  0 @.  \ 
\endCD
$$

\pf  In the top row and  in the left column  the
kernel clearly coincides with the kernel of $\d$,  given by the
torus in Proposition~1.14.  The surjectivity of $\d_1$  in the
top row is immediate. The surjectivity of $\d_2$ in the left
column follows directly  from the long exact sequence based on
$0\to \bbz\to\bbr\to S^1\to0$. However, we include a proof using
the de Rham-Federer representation (1.13).
 Note first that if $\<a\>\in H^k(X,S^1) =\oper{kernel}(\d_1)$, then
$da=-R\in \ci\cf^{k+1}(X)$ and hence $[R]\in
H^{k+1}_{\text{tor}}(X,\bbz)$. Conversely, if  $u\in 
H^{k+1}_{\text{tor}}(X,\bbz)$ has representative
 $R\in \ci\cf^{k+1}(X)$,
then the equation $da=-R$ can be solved with $a\in{\cd'}^k(X)$. 
The exactness of the bottom row and right column are obvious. 
The  commutativity of the diagram is straightforward to
verify.   \qed

   \vskip.3in

 

\subheading{\S2. Isomorphisms of de Rham - Federer type} 
Among the deepest and most useful theorems of Federer and de Rham
are those which assert that ordinary cohomology groups on a
manifold can be computed in many ways as the cohomology of 
distinct complexes of currents.  
In this section we prove analogous results for differential
characters, namely, we show that $\hH^*(X)$ has several
equivalent definitions whose interchangeablity is often useful in
practice.

The definition of $\hH^k(X)$ given in \S 1 used a ``maximal'' 
set of representatives.  We begin by examining  a nearly 
``minimal'' set.  Let ${\widetilde {\cc}}^k(X) \subset
{\cd'}^k(X)$ denote the currents defined by locally finite sums
of $C^{\infty}$ singular $(n-k)$-chains on $X$
(the {\sl current chains} of de Rham [deR]).  Note that 
$$
{\widetilde {\cc}}^k(X)\subset \cR^k(X) \subset \ci\cf^k(X).
$$
\medskip

\Def{2.1}
$$
\cs^k_{\mini}(X)\ \equiv\ \left\{a\in \ce^k_{\Lloc}(X)\ :\ 
da = \phi-R\ \text{where}\ \phi \in \ce^{k+1}(X)\ \text{and}\ 
R \in {\widetilde {\cc}}^{k+1}(X)
\right\}
$$
$$
\hH^k(X)_{\mini}\ \equiv\ \cs^k_{\mini}(X) / 
\left[d\ce^{k-1}_{\Lloc}(X) + {\widetilde {\cc}}^k(X)
\right]_0
$$
where by definition $V/[W]_0 = V/(W\cap V)$ for subspaces $V$
and $W$ of ${\cd'}^k(X)$.
\medskip

\Prop{2.2}  {\sl  The inclusion $\cs^k_{\mini}(X) \subset
\cs^k(X)$ induces an isomorphism of quotients}
$$
\hH^k(X)_{\mini}\ \cong\ \hH^k(X)
$$

\pf
We first prove surjectivity.  Suppose $\s \in{\cd'}^k(X)$,
with $d\s = \phi-R$ with $dR=d\phi=0$ as in 1.5.  Since 
$H^*({\widetilde {\cc}}^*(X)) \cong H^*(\ci\cf(X))$, the class
of $R$ is represented by a chain $C\in 
{\widetilde {\cc}}^{k+1}(X)$, i.e., $R=dS+C$, with $S\in
\ci\cf^k(X)$. Note that  $\s+S \cong S$ in $\hH^k(X)$ and
$d(\s+S)=\phi-C$. Since $H^*(\cf^*(X))\cong 
H^*({\cd'}^*(X))$ and $\phi-C\in \cf^{k+1}(X)$ represents the
zero class in $H^*({\cd'}^{k+1}(X))$, we can solve the equation
$df = \phi-C$ with $f\in \cf^k(X)$. Now the current $\s+S-f$
is $d$-closed and therefore modulo $d {\cd'}^{k-1}(X)$ it is
represented by a flat current.  That is, we have 
$$
\s+S-f = g+d\a
$$
for some $g\in \cf^k(X)$ and some $\a\in {\cd'}^{k-1}(X)$. 
Hence,  the current $\widetilde{\s} =\s+S-d\a$ is flat
and represents the same class as $\s$ in $\hH^k(X)$.

Now by (1.1) the flat current $\widetilde{\s}$ can be written
as $\widetilde{\s}=h+d\ell$ where $h$ and $\ell$ are forms with
$\Lloc$-coefficients. Thus $h\in \ce^k_{\Lloc}(X)$ represents
the same class as $\widetilde{\s}$ in $\hH^k(X)$.  Since $dh =
\phi-C$ with $C\in {\widetilde {\cc}}^{k+1}(X)$, we have 
$[h]\in \cs^k_{\mini}(X)$ and surjectivity is proved.

\def\cct{{\widetilde {\cc}}}

We now prove injectivity.  Consider $h\in \ce^k_{\Lloc}(X)$
with $dh=\phi-C$ where $\phi$ is smooth and  
$C\in {\widetilde {\cc}}^{k+1}(X)$.  Suppose that $h = d\a +R$
with $\a\in {\cd'}^{k-1}(X)$ and $R\in\ci\cf^k(X)$, i.e.,
suppose that $h$ represents the zero class in $\hH^k(X)$.
Then $dR = dh = \phi-C$.  By (1.4) this implies that $\phi = 0$
and $dR=-C$. Since 
$$
H^*({\widetilde {\cc}}^{*}(X))\cong
H^*(\ci\cf^*(X)),
\tag{2.3}
$$
 there exist a chain $K_1\in \cct^k(X)$
satisfying $dK_1=C$.  Therefore $R+K_1$  is $d$-closed and
integrally flat. Again using (2.3) we see that the class of 
$R+K_1$ is represented by a chain $K_2$, i.e., $R+K_1=dS+K_2$
with $S\in\ci\cf^{k-1}(X)$ and $K_2\in \cct^k(X)$.  Set
$C'\equiv K_1-K_2$.  Then $h+C'$ represents the same class as
$h \in \hH^k(X)_{\mini}$, and 
$$\aligned
h+C' &= d\a+R+K_1-K_2\\
&=d(\a+S),
\endaligned$$
with $\a+S\in {\cd'}^{k-1}(X)$. Since $H^*(\cf^*(X))\cong
H^*({\cd'}^*(X))$, we have $h+C'=dT$ for $T\in \cf^{k-1}(X)$.
Finally by (1.1) we have $T=f+dg$ where $f$ and $g$ are forms
with $\Lloc$-coefficients, and so $h +C'=df$ with $f\in
\ce^{k-1}_{\Lloc}(X)$.  This proves that $[h]=0$ in 
$\hH^k(X)_{\mini}$.   \qed
\medskip

In practice the most useful spaces of representatives for
differential characters are the following. (See [HL$_{1-4}$],
[HZ$_{1,2}$] for example.)

\Def{2.4}
$$\aligned
\hH^k(X)_0\ &= \ \frac{\{a \in \ce^k_{\Lloc}(X)\,:\, da = \phi-R, \ \
\ \phi\  \text{is smooth and $R$ is rectifiable}\} } 
                   { [d\ce^{k-1}_{\Lloc}(X)+\cR^k(X)]_0 }\\
\hH^k(X)_1\ &= \ \ \frac{\{a \in \cf^k(X)\,:\, da = \phi-R,
\ \ \ \phi\  \text{is smooth and $R$ is rectifiable}\} } 
{[d\cf^{k-1}(X) + \cR^k(X)]_0} \\
\hH^k(X)_2\ &= \ \frac{\{a \in \ce^k_{\Lloc}(X)\,:\, da = \phi-I, \
\ \ \phi\  \text{is smooth and $I$ is integrally flat}\} } 
{ [d\ce^{k-1}_{\Lloc}(X)+\ci\cf^k(X)]_0} \\
\hH^k(X)_3\ &= \ \frac{\{a \in \cf^k(X)\,:\, da = \phi-I,
\ \ \ \phi\  \text{is smooth and $I$ is integrally flat}\} } 
{d\cf^{k-1}(X) + \ci\cf^k(X)} 
\endaligned
$$
where by definition $V/[W]_0 = V/(W\cap V)$ for subspaces $V$
and $W$ of $\cf^k(X)$.   There is a natural commutative diagram 
of homomorphisms:
$$
\CD
 \ & \ &  \ & \ & \hH^k(X)_1 & \ & \  &  \ \\ 
 \ & \ & \ &  \nearrow & \ & \searrow & \  &  \ \\ 
  \hH^k(X)_{\mini}  @>>> \ \ \hH^k(X)_0   & \ &
 \ & \ & \hH^k(X)_3  @>>> \hH^k(X) \\ 
 \ & \ & \ &  \searrow & \ & \nearrow & \ &  \ \\ 
 \ & \ &  \ & \ & \hH^k(X)_2 & \ & \ &  \
\endCD
\tag{2.5}
$$

\Prop{2.6}  {\sl All the mappings in diagram (2.5) are 
isomorphisms.}

\pf
Surjectivity is immediate from 2.2.  Injectivity is
straightforward to check.
\qed

\medskip

\noindent
{\bf Note 2.7} \ \   An {\bf integral current} is by definition
a rectifiable current whose boundary is also rectifiable.
Thus in Definition 2.4, ``rectifiable'' can be everywhere
replaced by ``integral''  (in both the numerators
and the denominators).

\medskip

\noindent
{\bf Note 2.8 (concerning topologies).} \ \   As pointed
out in 1.15, the group $\hH^*(X)$ carries a natural topology
for which the homomorphisms $\d_1$ and $\d_2$ are continuous.
However this topology is far from being a quotient of any of the
standard topologies on the representing spaces of currents
discussed above. To see this consider the following 
basic examples of sparks.  Let $R$ be a rectifiable cycle
of degree $k$  on $X$, and let $H_{\epsilon}R$,
$0<\epsilon<1$ be Federer's family of {\bf smoothing homotopies}
[F; 4.1.18].  Then $H_{\epsilon}R$ satisfies the spark equation
$$
d( H_{\epsilon}R)\ =\ R_{\epsilon} - R
$$
where $R_{\epsilon}$ is a smooth $k$-form and 
$$
\bm(R_{\epsilon})\leq \bm(R)
\qquad\text{and}\qquad
\bm(H_{\epsilon}R) \leq \epsilon \bm(R)
$$
where $\bm(\bullet)$ denotes the mass norm.  Clearly 
$H_{\epsilon}R\arr 0$ as $\epsilon \to 0$, and so 
$dH_{\epsilon}R\arr 0$ by continuity.  However, 
$\d_1(H_{\epsilon}R) = R_{\epsilon}$ and $\d_2(H_{\epsilon}R) =R$ 
do not converge to 0.
The topology on sparks which does descend to the topology on
characters comes from the embedding 
$
i\times d_1\times d_2 :\cs^*(X) \to {\cd'}^*(X)\times
\ce^{*+1}(X)\times \ci\cf^{*+1}(X).
$

We note that Federer defines his smoothing homotopies 
only in $\bbr^N$. However, they  are easily transferred to any
manifold $X$ by embedding it in $\bbr^N$, smoothing the cycle in
a tubular neighborhood and then projecting back to $X$.

\medskip

\noindent
{\bf Remark 2.9 (functoriality).} \ \   The contravariant
functoriality of differential characters is an immediate
consequence of Theorem 4.1 below. However,  
this pull-back map can be constructed {\sl directly} on the de
Rham-Federer characters.  In fact for a projection map 
$\pi: X\times Y\to Y$ there is a well-defined pull-back mapping
on sparks which  descends to a pull-back $\pi^*$ on characters.
For a general mapping $f:X\to Y$ between manifolds we would like
to define the pull-back of a spark $a$ on $Y$ by the formula
$$
f^*a \ = \ \pi'_*\left(\pi^*a\wedge [\Gamma_f]\right)
\tag2.10
$$ 
where $\pi':X\times Y \to X$ is projection and $\Gamma_f\subset
X\times Y$ is the graph of $f$. The intersection theory of
Appendices A and B shows that given $a$ this is well
defined for almost all $f$. Furthermore, for fixed $f$,  it is
well defined   for arbitrarily small deformations of $\pi^*a$ on
the product. This is enough to give the pull-back on characters.
The details are given in the following paragraphs.
 
Consider first the projection map $\pi: X\times Y\to Y$, and
note that push-forward (integration over the fibre) takes smooth
forms with compact support on $X\times Y$ to smooth forms with
compact support on $Y$.  Hence there is a well defined
pull-back on currents.  Now let  $a\in \cs^k(X)$ be an
$\Lloc$-spark with  $da=\phi-R$ as in 1.5. Then it is
sraightforward to see that the pull-backs  $\pi^*a$, $\pi^*\phi$,
and  $\pi^*R$ are $\Lloc$, smooth and rectifiable respectively,
and satisfy the equation $d\pi^* a = \pi^*\phi- \pi^* R$. Hence
there is a well-defined mapping   
$\pi^*:\cs^k(Y)\to\cs^k(X\times Y)$
 Furthermore, if $a'=a + d b +S$ where $S$ is
rectifiable, then $\pi^*a'=\pi^*a+d\pi^* b+\pi^*S$, 
from which we see that our map on sparks determines a 
mapping of characters $\pi^*:\hH^k(Y)\to\hH^k(X\times Y)$

We now show that for any character $\a\in \hH^k(X\times Y)$, 
there are sparks $a \in\a$  (in fact a dense set of them) with
the property that $a\wedge [\Gamma_f]$ is well-defined flat
current, and $R\wedge [\Gamma_f]$ is a well-defined rectifiable
current, where $da=\phi-R$ as above.  To start fix any
$\Lloc$-spark $a_0\in\a$.  By Appendices A and B there is a
diffeomorphism $F$ of $X\times Y$, which can be taken
arbitrarily close to the identity in the $C^1$-topology, such
that $F_*a_0\wedge \Gamma_f$ is  a well-defined flat current,
and  $F_*R_0\wedge \Gamma_f$ is  a well-defined rectifiable
current where $da_0=\phi_0-R_0$ as always.  By de Rham-Federer
theory we have  $\phi_0-F_*\phi_0 = d\psi$ where $\psi $ is
smooth and  $R_0-F_*R_0 = dS$ where $S$ is rectifiable.
Since $d(a_0-F_*a_0-\psi+S) = 0$, de
Rham-Federer theory implies the existence of a smooth $k$-form
$\eta$ and a flat current $b$ such that
$a_0-F_*a_0-\psi+S=\eta+db$. The  spark $a = F_*a_0+\psi+\eta$,
which is equivalent to $a_0$, does the trick.

Taking $\pi'_*(a\wedge \Gamma_f)$, where $\pi':X\times Y\to X$
is projection, gives a spark on $X$. We claim that the
associated character is independent of the choice of ``good''
spark $a$ for $\a$ as above. Indeed if $a'$ is another good
spark, then by definition of the equivalence, $a-a'=dc +T$
where $T$ is rectifiable.  We again deform by a diffeomorphism
$F$ so that the deformed currents all meet $\Gamma_f$ properly
using Appendices A and B, and then we readjust $F_*a$ and $F_*a'$
as in the last paragraph.  Intersecting with $\Gamma_f$ and
pushing down to $X$ shows that the two sparks are equivalent. 
Thus we get a well-defined map $\Gamma_f^*:\hH^k(X\times Y)\to
\hH^k(X)$. Composition with $\pi^*$ gives a homomorphism
$f^*:\hH^k(Y)\to \hH^k(X)$.

The proof that $(f\circ g)^*=g^*\circ f^*$ is straightforward.

    \vskip .3in

 

\subheading{\S 3. Ring structure} In this section we introduce a
densely defined product on sparks which descends to characters
and makes $\hH^*(\bullet)$ into a graded-commutative ring.  We begin
with the following Proposition which is proved in  Appendix D.

\Prop{3.1}  {\sl  Given  classes 
$$
\a \in \hH^k(X) \qquad\text{and}\qquad
\b \in \hH^{\ell}(X).
$$
there exist representatives $a\in\a$ and $b\in\b$ with
$$
da \ =\ \phi -R \qquad\text{and}\qquad
db \ =\ \psi - S
$$
so that $a\wedge S$, $R\wedge b$ and $R\wedge S$ are well-defined
flat currents on $X$ and $R\wedge S$ is rectifiable.}

\bigskip

\Def{3.2}  Given $a$ and $b$ as in Proposition 3.1, define products
$$
\aligned
a*b  \ &\equdef \ a\wedge \psi \,+\, (-1)^{k+1}R\wedge b \\
a\tilde*b  \ &\equdef \ a\wedge S \,+\, (-1)^{k+1}\phi\wedge b
\endaligned
\tag{3.3}
$$
and note that by (B.3) 
$$
d(a*b) \ = \ d(a\tilde*b) \ = \ \phi\wedge \psi \,-\, R\wedge S.
\tag{3.4}
$$

\Theorem{3.5} {\sl   The classes $\langle a*b\rangle$ and
$\langle a\tilde*b\rangle$ in $\hH^{k+\ell+1}(X)$ agree and are
independent of the choice of representatives $a\in\a$ and $b\in\b$. 
Setting $$
\a * \b \ \equdef \ \langle a*b\rangle=\langle a\tilde*b\rangle  
$$
gives $\hH^*(X)$ the structure of a graded commutative ring  with
unit such that $\d$ is a ring homomorphism. 
}

\pf
To see that the two definitions agree first note that
$$
\aligned
a*b  - a\tilde*b \ &=\ a \wedge \psi +(-1)^{k+1}R\wedge b
-\left\{  a \wedge S +(-1)^{k+1}\phi\wedge b  \right\}  \\
&=  \   a\wedge (\psi - S)+(-1)^{k}(\phi-R)\wedge b  \\
&=  \  a\wedge db +(-1)^{k} da \wedge b.
\endaligned
$$
We now apply the following.

\Lemma{3.6}  {\sl  Let $a,b \in {\cd'}^*(X)$ be currents for which
the products $da\wedge b$ and $a\wedge db$ are well defined.  Then
$u = da\wedge b +(-1)^{\text{deg}(a)}a\wedge db$ is an exact current,
i.e., $u = d v$ for some current $v$ on $X$.
}
\pf
Choose smooth approximations $a_{\e}, b_{\e} \in \ce^*(X)$ with 
$a_{\e} \to a$ and $b_{\e} \to b$  as $\e\to 0$.  Then 
$da_{\e}\wedge b_{\e} +(-1)^{\text{deg}(a)}a_{\e}\wedge b_{\e}  =
d(a_{\e}\wedge b_{\e})$, and the lemma now follows from the fact
(cf. [deR]) that the range of $d$ is closed.  \qed

We now show that the product is independent of the choice of
representatives $a$ and $b$ chosen.  To begin suppose that 
$\tilde a = a + du$ for some current $u$.  Then
$$
\tilde a * b \ -\ a * b \ = \ du \wedge \psi \ = \ d(u\wedge \psi)
$$ 
and so $\tilde a *b$ and $a*b$ define the same class in $\hH^*(X)$.

Suppose now that $\tilde a = a +Q$ where $Q$ is locally rectifiable,
and write  $da = \phi - R$ and $d\tilde a = \phi - \widetilde R$ as in
1.5.
 Note that by hypothesis $dQ = d \tilde a -da = R -\tilde R$
meets $S$ properly.  Using a deformation which keeps $dQ$ fixed can find
$\tilde Q = Q + d\sigma$ where $\tilde Q$
meets $S$ properly.  Using the paragraph above and the second
definition of the product we have
$$\aligned
\tilde a * b \ -\ a * b \ &= \ (\tilde a + d\sigma) * b\ -\ a * b 
\\  &= \  \tilde Q \wedge S \ \equiv\ 0 \qquad\text{in}\ \ \hH^*(X).
\endaligned
$$
This shows  that the product $[a*b]$ is independent of the choice
of $a \in [a]\in \hH^k(X)$.  The argument for  $b$ is completely
analogous.

A direct calculation using (3.3) now shows that 
$$
\a *\b \ = \ (-1)^{(k+1)(\ell+1)} \b \tilde{*}\a
\ \equiv \ (-1)^{(k+1)(\ell+1)} \b {*}\a
\tag{3.7}
$$
and the multiplication is graded-commutative as claimed. \qed

     \vskip .3in

 

\subheading{\S 4. Cheeger-Simons characters} In their fundamental
paper [CS]   Cheeger and Simons defined the differential characters
of degree $k$ with coefficients in $\bbr/\bbz$ on a manifold $X$ to
be the group  
$$
\hh^k(X; \bbr/\bbz)\ =\ \left\{ h \in \Hom( \cz_k(X), \bbr/\bbz)\
: \   \d h \equiv \phi \ ( \text{mod}\, \bbz)\ \text{\ for some\
} \phi\in \ce^{k+1}(X)\right\}
$$
where $\cz_k(X)$ denotes the smooth singular $k$-chains on $X$ with
$\bbz$-coefficients and $\d$ the standard codifferential.
  The first result of this section is the following.

\Theorem{4.1}  {\sl  For any manifold $X$ there is a natural
isomorphism 
$$
\Psi : \hH^k(X) @>{\cong}>> \hh^k(X; \bbr/\bbz)
$$
induced by integration.} 

\medskip

To define the homomorphism $\Psi$ we will need the following.

\Prop{4.2}  {\sl  Given any class $\a\in\hH^k(X)$ and any $Z \in
\cz_k(X)$, there is a representative spark  $a \in \a$ which is
smooth on a neighborhood of supp$(Z)$.  Furthermore, given any two
such representives $a,a' \in \a$, one has that}
$$
\int_Z a \ \equiv\ \int_Z a' \ \ (\text{mod} \ \ \bbz).
$$

\Def{4.3}  Given $\a \in \hH^k(X)$ and  $Z \in \cz_k(X)$ we define
$$
\Psi(\a)(Z)\ = \ \int_Z a  \ \ (\text{mod} \ \ \bbz)
$$
where $a$ is chosen as in Proposition 4.2.  By 4.2 it is clear that
$\Psi(\a)$ is well defined and that $\Psi$ defines a homomorphism
into $S^1 = \bbr/\bbz$.

\medskip
\noindent
{\bf Proof of 4.2.} We begin with the following assertion which  is
of independent interest.

\Lemma{4.4} {\sl Given any locally rectifiable cycle $R\in \d_2(\a)
\in  H^{k+1}(X;\,\bbz)$, there is an $\Lloc$ spark 
$a\in \a$ such that $da = \d_1(\a) - R$.}

\pf
By Proposition 2.6 there exists an $\Lloc$-form         
$\widetilde a \in \a$ with $d\widetilde a = \phi - R_0$ where
$\phi = \d_1(\a)$ is a smooth form and $R_0$ is a
locally rectifiable cycle of dimension $n-k-1$.  Since the
integral cohomology class of $R_0$ is $\d_2(\a)$, there exists a
locally rectifiable current $S$ with $dS =  R_0 - R$.   The current
$\widetilde a + S$ is flat and can be written as $a +db$ where $a$
and $b$ are $\Lloc$-forms (cf. (1.1)). Then $a\in \a$  has the
asserted properties. \qed

Now the $C^\infty$-singular chains over $\bbz$  map to the
rectifiable (in fact, integral) currents and their image computes
integral cohomology. For dimensional reasons we may choose such a
cocycle $R\in \d_2(\a)$ with   
$$
\oper{supp}(R) \cap \oper{supp}(Z) \ = \ \emptyset.
\tag{4.5}
$$
We then choose $a \in \a$ as in Lemma 4.4.  In
$\Omega \equiv X-\oper{supp}(R)$ we have $da = \phi$, a smooth
form.

\Lemma {4.6} {\sl Given $a \in \ce^{k}_{\Lloc}(\Omega)$ such that
$da$ is smooth,  there exists $b \in \ce^{k-1}_{\Lloc}(X)$ such that
$a-db$ is smooth.}

\pf  This follow directly form the fact that the
cohomology of locally flat currents coincides with the 
de Rham cohomology of smooth forms [F]. \qed\medskip

 Now choose $b$ as in Lemma 4.6, and  choose $f\in
C_0^{\infty}(\Omega)$ with  $f\equiv 1$ on a neighborhood of the
compact subset $\oper{supp}(Z)\subset \Omega$. Set $b_0=f b$ and
extend it by zero to all of $X$. Then $a_0=a - d b_0$ is equivalent
to $a$ and is smooth on a neighborhood of $\oper{supp}(Z)$. 

To prove the second assertion of 4.2 
let  $a' \in \a$ be another such choice.  Then we have
$$
da = \phi - R
\qquad \text{and}\qquad
da' = \phi - R'
$$
where $\phi = \d_2(\a)$, $a$ and $a'$ are both smooth on 
$\oper{supp}(Z)$, and $R$ and $ R'$ are $C^{\infty}$-singular
chains.     Now $d(a-a') = R'-R = dS$ where $S$ is a  
$C^{\infty}$-singular chain which we may assume is transversal to
$Z$.  Thus $(a-a' -S, Z) \equiv (a-a', Z) \mod \bbz$. However, since
$[a] = [a'] = \a$ we have $a  = a' +S_0 +db$ where $S_0$ is
rectifable and $b$ is flat.  This shows that the  cycle $a-a'-S$ is
homologous to the rectifiable cycle $S_0-S$, and so  $(a-a' -S, Z)
\equiv 0 \mod \bbz$ as desired. \qed

\medskip

Note that the arguments above actually prove the following.
\Cor{4.7}  {\sl Fix $\a \in \hH^k(X)$ and  $Z \in \cz_k(X)$.
Then given any rectifiable cycle $R\in \d_2(\a)$,
whose support is disjoint from $\oper{supp}(Z)$, there is a
representative $a\in \a$ with $\d_2(a)=R$ which is smooth on a
neighborhood of  $\oper{supp}(Z)$.}
\medskip

\noindent
{\bf Proof of Theorem 4.1.}
We now recall that Theorem 1.1 of [CS] establishes a functorial
short exact sequence:
$$
0\ \arr \ \frac{H^k(X;\, \bbr)}{H^k_{\oper{free}}(X;\, \bbz)}
\ \arr\ \hh^k(X; \bbr/\bbz) \ @>{\d}>> Q^{k+1}(X) \ \arr\ 0.
\tag{4.8}
$$
One checks directly that the mapping $\Psi$ induces a commutative
diagram
$$
\CD
0 @>>> \frac{H^k(X;\, \bbr)}{H^k_{\oper{free}}(X;\, \bbz)}
  @>>> \hH^k(X)   @>{\d}>>  Q^{k+1}(X)  @>>>  0  \\
@.  @VVV   @VV{\Psi}V   \|   \\
0 @>>> \frac{H^k(X;\, \bbr)}{H^k_{\oper{free}}(X;\, \bbz)}
  @>>> \hh^k(X; \bbr/\bbz)   @>{\d}>>  Q^{k+1}(X)  @>>>  0  
\endCD
\tag{4.9}
$$
where the top  line is the sequence from 1.11, the bottom line
is (4.8), and  the  left vertical arrow is an
isomorphism.  It follows immediately that $\Psi$ is an
isomorphism. This proves the theorem. \qed

\bigskip

From the definition it is obvious that the Cheeger-Simons
characters are a contravariant functor on the category of smooth
manifolds. What is not obvious is that they carry a graded ring
structure for which $\delta_1$ and $\delta_2$ are ring
homomorphisms.  Such a product was established by Cheeger [C].
Our next result asserts that under the isomorphism $\Psi$,
our product coincides with the Cheeger product.  While our product
is rather straightforward to define, the proof of its coincidence
with the Cheeger product is not trivial.  The difficulty stems from
the fundamental difference between the intersection product of
cochains, which is local in nature, and the cup product, which is
not (cf. Note C.2).

\Theorem{4.10}  {\sl Under the isomorphism $\Psi$ given in Theorem
4.1,  the $*$-product on $\hH^*(X)$ coincides with the Cheeger
product defined on $\hh^*(X;\, \bbr/\bbz)$.}

\pf   Recall that the additive isomorphism $\Psi^{-1}:
\widehat{H}^*(X;\,\bbr/\bbz) \to\hH^*(X)$   preserves $\d_1$ and 
$\d_2$ and therefore introduces two ring structures on $\hH^*(X)$
for which these maps are ring homomorphisms.  Thus the biadditive map
$$
\rho: \hH^*(X) \times \hH^*(X) \ \arr\ \hH^*(X)
\qquad\text{given by}\qquad \rho(\a,\b)\equiv \a*\b - \a \hat{*}\b
$$
has the property that 
$$
\delta_1\rho(\a,\b)\ =\ \delta_2\rho(\a,\b)\ =\ 0
$$
for all $\a,\b$. Thus the image of $\rho$ is contained in the 
subgroup $H^{*}(X;\, \bbr)/H^{*}(X;\, \bbz)$.
Our next observation is the following.

\Prop{4.11}{\sl   If either $\a$ or $\b$ is represented by a smooth
differential form, then $\rho(\a,\b)=0$.
}

\pf  Suppose $\b$ is represented by a
smooth differential $\ell$-form $b$.  Via integration $b$
defines a real cochain which when restricted to cycles gives the
character $\b$.  Note that $\d_1(\b)= db \equiv \psi$ and
$\d_2(\b)=0$.  Let $\d_1(\a) \equiv \phi$.  Then by the definition given
in [C] we have 
$$
\a\hat*\b \ \equiv \ (-1)^{k+1}\phi\cup b + E(\phi,\psi) \qquad
\text{mod}\ \bbz 
$$
where $E$ is defined as follows.  Let $\phi_1, \phi_2$ be smooth forms
of degrees $k_1, k_2$ respectively.  Let $\Delta^i$ denote the
$i^{\text{th}}$ subdivision operator acting on $C^{\infty}$ cubical
chains.  Let $h$ be the standard chain homotopy on cubical chains
satisfying $1-\Delta=\partial \circ h + h \circ\partial$. Then 
$E(\phi_1,\phi_2)$ is a real cochain of degree $k_1+k_2+1$ defined 
on a $C^{\infty}$ cubical chain $c$ by 
$$
E(\phi_1,\phi_2)(c)\ =\ -\sum_{i=0}^{\infty} \, (\phi_1\cup\phi_2)
(h\Delta^ic).
$$
Cheeger proves that this series converges and when 
$d\phi_1=d\phi_2=0$, one has $\d E(\phi_1,\phi_2) = \phi_1\wedge\phi_2
- \phi_1\cup \phi_2$ where $\d$ denotes the coboundary. Recalling that
$d\phi = 0$ and $d b = \psi$ we compute that
$$
\aligned
\d E(\phi, b)(c) \ &= \ -\sum_{i=0}^{\infty} \, (\phi\cup b)
( h\Delta^i\partial c) \ = \ -\sum_{i=0}^{\infty} \, (\phi\cup b)
( h\partial \Delta^ic)
 \\
&= \ -\sum_{i=0}^{\infty} \, (\phi\cup b)
\{(1-\Delta -\partial  h)\Delta^ic\} \\
&= \ -\sum_{i=0}^{\infty} \, (\phi\cup b)
\{(1-\Delta)\Delta^ic\} + \sum_{i=0}^{\infty} \, \{d(\phi\cup b)\}
(h\Delta^ic)  \\
&= \ -\lim_{n\to\infty}\sum_{i=0}^{n} \, (\phi\cup b)
\{(1-\Delta)\Delta^ic\} + \sum_{i=0}^{\infty} \, \{d(\phi\cup b)\}
(h\Delta^ic)  \\
&= \ -\lim_{n\to\infty}\sum_{i=0}^{n} \, (\phi\cup b)
(1-\Delta^{n+1}c)  + (-1)^{k+1}\sum_{i=0}^{\infty} \, (\phi\cup d b)
(h\Delta^ic)  \\
&=\ (\phi\wedge b - \phi\cup b)(c) + (-1)^k E(\phi,\psi).
\endaligned
$$
Therefore,
$$
\aligned
\a\hat*\b - (-1)^{k}\d E(\phi,b) \ &=\ 
(-1)^{k+1}\phi\cup b + E(\phi,\psi) -
(-1)^{k}\{\phi\wedge b - \phi\cup b
+(-1)^kE(\phi,\psi)\} \\
&=\ (-1)^{k+1}\phi\wedge b \\ & = \ \a*\b
 \endaligned
$$
Consequently, $\a\hat*\b$ and $\a*\b$ have the same values on cycles
and therefore they define the same differential characters.   \qed
\medskip

From 1.11 and 1.14 we conclude that $\rho$ descends to a biadditive
mapping
$$
\widetilde{\rho}:H^{k+1}(X;\,\bbz)\times H^{\ell+1}(X;\,\bbz)\ \arr\ 
H^{k+\ell+1}(X;\,\bbr)/H^{k+\ell+1}_{\text{free}}(X;\,\bbz)
$$
which we will  show to be zero. To compute $\widetilde{\rho}$
fix classes $(u,v)\in H^{k+1}(X;\,\bbz)\times H^{\ell+1}(X;\,\bbz)$
and choose good cycles $R\in u$ and $S\in v$, e.g., cycles from
smooth triangulations such that all simplices in one triangulation
are transversal to all simplices in the other.  Let $H_{\epsilon}R$
and $H_{\epsilon}S$ be the smoothing homotopies from 2.8. 
We only need to evaluate $\widetilde{\rho}(R,S)$
on one cycle in each class of  $H_{k+\ell+1}(X;\,\bbz)$, so we may
assume these cycles are nice with respect to $R$ and $S$ (as above).
For such a cycle $Z$ it is an direct calculation that 
$(H_{\epsilon}R * H_{\epsilon}S - 
H_{\epsilon}R \hat* H_{\epsilon}S)(Z) \arr 0$ as $\epsilon\to 0$.
\qed

     \vskip .3in

%

\subheading{\S 5. Dual sequences}   As noted in 1.15 the group of
characters carries  a  natural topology having $\hH^*_{\infty}$
as the connected component of 0, and inducing the standard 
$C^{\infty}$-topology on $\cz^*_0(X)$ in 1.12.
It therefore makes sense to consider the {\sl Pontrjagin dual}
which we define to be the group
$$
\hH^k(X)^{*}\ \equiv\ \Hom\left(\hH^k(X), S^1\right)
$$
of {\bf continuous} group homomorphisms from $\hH^k(X)$ to the circle.  
The interesting fact is that this dual group sits in two fundamental
exact sequences which resemble the fundamental sequences 1.11 and 1.12
for  $\hH^{n-k-1}(X)$ but which are interchanged under duality.

Here we examine these sequences and lay the groundwork for
duality which will be established next. The symbol $\Hom$ shall always
denote {\bf continuous} homomorphisms. We assume throughout this
section that $X$ is compact, oriented and of dimension $n$.

\Lemma{5.1} {\sl There are  natural isomorphisms:
$$\aligned
\Hom(H^k(X;\,\bbz),\, S^1)
\ &\cong\  H^{n-k}(X;\,S^1)  
 \qquad\text{(Poincar\'e Duality)}       \\
\Hom\left(\frac{H^k(X;\,\bbr)}
{H^k_{\text{free}}(X;\,\bbz)},\, S^1\right)
\  
&\cong\  
H^{n-k}_{\text{free}}(X;\,\bbz) \\
\Hom( d\ce^{k}(X),\, S^1) \ &\cong\  
d{{\cd'}}^{n-k-1}(X)
\endaligned$$
}

\pf
By Pontrjagin duality,  classical Poincar\'e duality,  
 $H^k(X;\,\bbz) \cong H_{n-k}(X;\,\bbz)$,
is equivalent to
 $\Hom(H^k(X;\,\bbz),\, S^1) \cong 
\Hom(H_{n-k}(X;\,\bbz),\, S^1)$.
By the Universal Coefficient Theorem,
$\Hom(H_{n-k}(X;\,\bbz),\, S^1) \cong H^{n-k}(X;\,S^1)$. This
proves the first assertion.

To prove the second assertion we recall that if
$\Lambda\subset\bbr^m$ is a lattice of full rank, then
$$
\Hom(\bbr^m/\Lambda, S^1) \cong \Hom(\Lambda, \bbz)
\qquad\text{and}\qquad
\Hom(\Lambda, S^1) \cong \Hom(\Lambda, \bbr)/\Hom(\Lambda,\bbz).
\tag{5.2}$$
Now by (5.2) and Poincar\'e duality we have 
$
\Hom\left(\frac{H^k(X;\,\bbr)}{H^k_{\text{free}}(X;\,\bbz)},\, S^1\right)
\cong
\Hom(H^k_{\text{free}}(X;\,\bbz), \bbz)
\cong 
\Hom(H_{n-k}(X;\,\bbz)_{\text{free}},\, \bbz)
\cong  H^{n-k}_{\text{free}}(X;\,\bbz)$.  

For the third assertion first note that  
$\Hom( d\ce^{k}, \bbr/\bbz) = \Hom( d\ce^{k},\bbr) =
(d\ce^k)'$.
To identify this vector space dual, consider the sequence 
$
\ce^k\  \overset{d}\to{\surrightarrow}\  d\ce^k 
\ \subset \ \ce^{k+1}. 
$
On any manifold the range of $d$ is closed (cf. [deR]).
Hence, by Hahn-Banach the right hand mapping in the dual
sequence
$
(\ce^k)'\  \hookleftarrow\ ( d\ce^k)'\ \surleftarrow (\ce^{k+1})'
$
is surjective.   This composition is exactly the exterior derivative
on currents, and so we conclude that
$(d\ce^k)' = d\{({\ce}^{k+1})'\} = d\{{\cd'}^{n-k-1}\}$.   \qed

\medskip

\noindent
{\bf Observation 5.3} \ \  
There are (non-canonical) group homomorphisms:
$$
\hH^k(X) \ \cong \ d\ce^k(X)\times 
\frac{H^k(X;\,\bbr)}{H^k_{\text{free}}(X;\,\bbz)}\times
H^{k+1}(X;\,\bbz)
$$
$$
\hH^{n-k-1}(X) \ \cong \ d\ce^{n-k-1}(X)\times 
H^{n-k}_{\text{free}}(X;\,\bbz)\times
H^{n-k-1}(X;\,S^1)
$$
The first follows from 1.11 and the top row of 1.16, and the
second  from 1.12 and the right column of 1.16.  Combined with 5.1
above this justaposition is quite suggestive. To pin down the implied
duality we examine the Pontrjagin duals of the sequences in \S 1.

\medskip

Taking the duals of the sequences 1.11 and 1.16 (top row)  gives us
$$
0\to\Hom(H^{k+1}(X;\,\bbz),\, S^1)
\to\Hom(\hH^k(X),\, S^1)
\to\Hom(\hH^k_{\infty}(X),\, S^1)\to 0
$$
$$
0\to\Hom(d\ce^{k}(X),\, S^1)
\to\Hom(\hH^k_{\infty}(X),\, S^1)
\to\Hom\left(\frac{H^k(X;\,\bbr)}{H^k_{\text{free}}(X;\,\bbz)},\,
S^1\right)\to 0 
$$
Applying Lemma 5.1 gives us
\Prop{5.4} {\sl There are exact sequences}
$$
0\to H^{n-k-1}(X;\, S^1)
@>{\d_2^*}>>\Hom(\hH^k(X),\, S^1)
@>{\rho_2}>>\Hom(\hH^k_{\infty}(X),\, S^1)\to 0
\tag{5.5}
$$
$$
0\to d{\cd'}^{n-k-1}(X)
@>{\d_1^*}>>\Hom(\hH^k_{\infty}(X),\, S^1)
@>{\rho_1}>> H^{n-k}_{\text{free}}(X;\,\bbz)
\to 0 
\tag{5.6}
$$
\medskip

Similarly,  taking the duals of 1.12 and 1.16 (right column) and then
applying Lemma 5.1 and Pontrjagin duality gives
\Prop{5.7} {\sl There are exact sequences}
$$
0\to \Hom(\cz_0^{k+1}(X),\, S^1)
@>{\d_1^*}>>\Hom(\hH^k(X),\, S^1)
@>{\rho_1}>>  H^{n-k}(X;\, \bbz)\to 0
\tag{5.8}
$$
$$
0\to \frac{H^{n-k-1}(X;\,\bbr)}{H^{n-k-1}_{\text{free}}(X;\,\bbz)}
\to\Hom(\cz_0^{k+1}(X),\, S^1)
@>{}>>  d{\cd'}^{n-k-1}(X)
\to 0 
\tag{5.9}
$$

We now come to a central concept in this paper.

\Def{5.10}  A homomorphism $h:\hH^k_{\infty}(X) \to S^1$ is called {\bf
smooth }  if there exists a smooth form $\omega \in \cz_0^{n-k}(X)$ such
that 
$$
h(\a) \ \equiv\ \int_X a\wedge \omega \ \ (\text{mod} \ \bbz)
\tag{5.11}
$$
for $a\in\a$.  One checks easily that the right hand side is independent
of the choice of $a\in\a$.  Let $\Hom_{\infty}(\hH^k_{\infty}(X),\, S^1)$
denote the group of all such homomorphisms.  We then define the {\bf smooth
Pontrjagin dual} of $\hH^k(X)$ to be the group 
$$
\Hom_{\infty}(\hH^k(X),\, S^1) \ 
\equdef\ \rho_2^{-1}\Hom_{\infty}(\hH^k_{\infty}(X),\,S^1)
$$
where $\rho_2$ is the restriction homomorphism in (5.5).

\medskip

The obvious map induces an isomorphism
$
\cz_0^{n-k}(X)\ \cong\ \Hom_{\infty}(\hH^k_{\infty}(X),\, S^1)
$
and the restriction of (5.6) to the smooth homomorphisms becomes
$$
0\to d{\ce}^{n-k-1}(X)
\to\Hom_{\infty}(\hH^k_{\infty}(X),\, S^1)
\to H^{n-k}_{\text{free}}(X;\,\bbz)
\to 0 
\tag{5.12}
$$
Let $\Hom_{\infty}(\cz_0^{\ell}(X),\, S^1) = \Hom(\cz_0^{\ell}(X),\, S^1)
\cap \Hom_{\infty}(\hH^{\ell-1}(X),\, S^1)$ \ in (5.8).  Restricting
(5.9) gives a short exact sequence
$$
0\to \frac{H^{n-k-1}(X;\,\bbr)}{H^{n-k-1}_{\text{free}}(X;\,\bbz)}
\to\Hom_{\infty}(\cz_0^{k+1}(X),\, S^1)
@>{}>>  d{\ce}^{n-k-1}(X)
\to 0 
\tag{5.13}
$$

\Prop{5.14} {\sl There are natural isomorphisms}
$$
\cz_0^{n-k}(X)\ \cong\ \Hom_{\infty}(\hH^k_{\infty}(X),\, S^1)
\qquad\text{and}\qquad
\hH^{n-k-1}_{\infty}(X) \ \cong\ \Hom_{\infty}(\cz_0^{k+1}(X),\, S^1)
$$

\pf
The first isomorphism is noted above.  For the second we define a mapping
$\hH^{n-k-1}_{\infty}(X) \to \Hom_{\infty}(\cz_0^{k+1}(X),\, S^1)$ by 
formula (5.11) above.  Note again that this integral mod $\bbz$ does not
depend on the choice of representative $a\in\a\in \hH^{n-k-1}_{\infty}(X)$.  We
then observe that this mapping induces a mapping of short exact sequences 
$$\CD
0 @>>> \frac{H^{n-k-1}(X,\bbr)}{H^{n-k-1}_{\oper{free}}
(X, \bbz)}   @>>>  \hH^{n-k-1}_\infty(X)  
@>{\d_1}>> d\ce^{n-k-1}(X) @>>>  0  \\
@. @V{\cong}VV  @VV{}V  @VV{\cong}V  \\
0 @>>> \frac{H^{n-k-1}(X;\,\bbr)}{H^{n-k-1}_{\text{free}}(X;\,\bbz)}
@>>> \Hom_{\infty}(\cz_0^{k+1}(X),\, S^1)
@>{}>>  d{\ce}^{n-k-1}(X)
@>>> 0 
\endCD
$$
from 1.16 (top row) to (5.13). The right and left vertical maps are
isomorphisms and so therefore is the middle.  \qed

\medskip

This leads to the main result of this section.

\Theorem{5.15}  {\sl  If $X$ is compact and  oriented, 
there are short exact sequences:
$$
0 \  \arr \ 
H^{n-k-1}(X;\,S^1)
\  \arr \  \Hom_{\infty}\left(\hH^{k}(X),\, S^1\right)
 @>{\rho_2}>>  \cz^{n-k}_0(X)
\  \arr \  0,
$$
$$
0 \  \arr \ 
\hH_{\infty}^{n-k-1}(X)
\  \arr \  \Hom_{\infty}\left(\hH^{k}(X),\, S^1\right)
 @>{\rho_1}>>  H^{n-k}(X;\,\bbz)
\  \arr \  0,
$$
}

\pf
The first comes from (5.5) and 5.14; the second from 
(5.8) and 5.14.  \qed

      \vskip .3in


\subheading{\S 6. Duality for Differential Characters} 
Let $X$ be a compact connected
oriented  manifold of dimension $n$.  Then the top character group 
$\hH^n(X) = {\cd'}^n(X)/\{d{\cd'}^{n-1}(X)+\ci\cf^n(X)\}$ is easily
seen to satisfy
$$
\hH^n(X) \ \cong \ \Hom\left(\cz_n(X), \bbr/\bbz\right) \ \cong \ 
\Hom\left(\bbz, \bbr/\bbz\right) \ \cong \ \bbr/\bbz
$$
As noted in Remark 1.15 each group of de Rham characters has a
natural topology.  In what follows $\Hom$ shall always denote 
{\sl continuous} group homomorphisms.

\Def{6.1}  For each integer $k$, $0\leq k< n$  we define the
duality mapping
$$
\cd : \hH^{n-k-1}(X) \ \arr\ \hH^{k}(X)^*\equiv
\Hom_{\infty}\left(\hH^{k}(X), \bbr/\bbz\right)
$$
by
$$
\cd(\a)(\b) \ \equdef \ \a *\b \in \hH^n(X) = \bbr/\bbz.
$$

\medskip

\Theorem {6.2.  (Poincar\'e-Pontrjagin Duality for Characters)}
{\sl The duality mapping $\cd$ is an isomorphism.}

\medskip

\pf
Given $a \in \a \in \hH^{n-k-1}(X)$ and 
$b \in \b \in \hH^{k}(X)$ we write $da = d_1a-d_2a$
and $db = d_1b-d_2b$ where $d_1a, d_1b$ are smooth and
$d_2a, d_2b$ are rectifiable.  Then (3.3) can be rewritten
$$
 a*b \ = \ a\wedge d_1b + (-1)^{k+1} d_2a\wedge b 
\tag{6.3}
$$
$$
 a\tilde* b \ = \ a\wedge d_2b + (-1)^{k+1} d_1a\wedge b 
\tag{6.4}
$$
It is easily seen that $a$ and $b$ can be chosen so that the
terms in these equations are well defined.  Indeed, since
$\dim d_2a + \dim d_2 b = n-1$, we may choose $a$ and $b$ with 
supp$(d_2a) \cap \text{supp}(d_2b) = \emptyset$, and then by
Corollary 4.7 we may assume that $a$ is smooth on supp$(d_2b)$
and  $b$ is smooth on supp$(d_2a)$.

Consider now the composition 
$\rho_1\circ \cd : \hH^{n-k-1}(X) \arr \cz^{n-k}_0(X)$.
Note that $\rho_1$ is defined by restriction of homomorphisms to
the subgoup  $\ker \d_2 \subset \hH^k(X)$.  Now if $\d_2\b = 0$,
we may choose $b\in\b$ to be a smooth form.  In particular, $d_2b
= 0$ and by  (6.4)
$$
\cd(a)(b)\ =\ (-1)^{k+1} d_1a \wedge b
$$
which implies that
$$
\rho_1\circ \cd\ =\ (-1)^{k+1} d_1.
\tag{6.5}
$$

Similarly we have the composition
$\rho_2\circ \cd : \hH^{n-k-1}(X) \arr H^{n-k}(X;\,\bbz)$ 
where $\rho_2$ is defined by restriction to $\ker \d_1$.
If $\d_1\b = 0$, we may choose $b\in\b$ with $d_1b=0$, and by 
(6.3) we have
$$
\cd(a)(b)\ =\ (-1)^{k+1} d_2a \wedge b.
$$
This shows that
$$
\rho_2\circ \cd\ =\ (-1)^{k+1} d_2.
\tag{6.6}
$$

Combining (6.5) and (6.6) and using Theorem 5.15 we get a
commutative diagram:
$$
\CD
0 @>>>
\frac{H^{n-k-1}(X;\,\bbr)}{H^{n-k-1}_{\oper{free}}(X;\,\bbz)}
@>>>   \hH^{n-k-1}(X)  @>{\d=(\d_1,\d_2)}>>  R^{n-k}(X)  @>>> 0\\
@. @V{\e }VV @VVV @VV{(-1)^{k+1}}V  \\ 
0 @>>> 
\frac{H^{n-k-1}(X;\,\bbr)}{H^{n-k-1}_{\oper{free}}(X;\,\bbz)}
@>>>   \Hom_{\infty}\left(\hH^{k}(X),\, \bbr/\bbz\right)  
@>{\rho=(\rho_1,\rho_2)}>>  R^{n-k}(X)  @>>> 0    
\endCD 
$$
It remains to show that the mapping $\e$ is an isomorphism.
However, this map asociates to a smooth closed $(n-k-1)$-form
$\psi$ the element 
$$
[\psi] \in H^{n-k-1}(X;\,\bbr)/
H^{n-k-1}_{\oper{free}}(X;\,\bbz) \cong
\Hom\left(H_{n-k-1}(X;\,\bbz)_{\oper{free}},\, \bbr/\bbz\right)
$$
given by 
$$
[\psi](z) \ \equiv\ \int_z \psi \ \ (\text{mod}\ \bbz).
$$
This map is clearly an isomorphism.\qed

\vskip .3in

\noindent
{\bf A Second Proof of Duality.}  We present here a second proof of
Theorem 6.2 which does not use the dual sequences.  This argument
will be used in a subsequent article to establish a general  
``Alexander-Pontrjagin'' duality for characters.

We first prove that the pairing is non-degenerate.  Fix $\a \in 
\hH^{n-k-1}(X)$ and assume that 
$$
(\a*\b) \,[X] \ =  \ 0
$$
for all $\b \in\hH^k(X)$. We shall prove that $\a=0$.

To begin consider a class $\b$ which is represented by a smooth
form $b$. Then
$$
(\a*\b) [X] \ = \ (-1)^{n-k} \int_X \d_1\a \wedge b\ \equiv\ 0
\ \ (\text{mod}\ \ \bbz).
$$
Since this holds for all smooth $k$-forms $b$, we conclude that
$\d_1\a=0$. Hence, there is a representative $a\in \a$ with 
$da = R$ a rectifiable cycle whose integral homology class is
torsion.

We shall show that this homology class is zero.
Choose a class $u\in H^{k+1}(X;\,\bbz)_{\text{tor}}$ of order
$m$, and let  $T$ be a rectifiable cycle with $[T]=u$.  Since
$m[T]=0$ there exists a rectifiable current $S$ of degree $k$
with $mT=dS$.  Set $\b = [\frac 1 m S] \in \hH^{k}(X)$ and
note that $\d_1\b=0$. Now
$$\aligned
(\a*\b) [X] \ &=\ (-1)^{n-k} (d_2 a \wedge b) [X] \ =\ 
(-1)^{n-k}( R\wedge {\tsize \frac 1 m}S) [X]  \\
&=\ (-1)^{n-k}  {\tsize \frac 1 m}\{\text{intersection number of 
$R$ with $S$}\}\\
&=\ (-1)^{n-k} Lk([R], u)\\
& \equiv\ 0 \ \ (\text{mod}\ \ \bbz)
\endaligned$$
where $Lk$ denotes the Seifert-deRham linking number. Since this
holds for all $u$, the non-degeneracy of this linking pairing on
torsion cycles implies that $[R]=0$ in $H^{n-k}(X;\,\bbz)$.
Hence $\d_2\a=0$ and so after adding an exact current we may
assume that $a$ is a smooth $d$-closed differential form of
degree $n-k-1$. 

Now choose any class $v\in H^{k+1}(X;\,\bbz)\cong
H_{n-k-1}(X;\,\bbz)$ and let $S$ be a rectifiable cycle in $v$. 
Let $\phi$ be a smooth $d$-closed form representing $v\otimes\bbr$
and choose a current $b$ with $db=\phi-S$.  Let
$\b=[b]\in\hH^k(X)$.  Then $$
(\a*[b])[X]\ =\ (a\wedge d_2b)[X]\ \equiv \ \int_S\, a \ \equiv\ 0
\ \ (\text{mod}\ \ \bbz).
$$
Hence, $[a]=0$ in
$H^{n-k-1}(X;\,\bbr)/H^{n-k-1}_{\text{free}}(X;\,\bbz)
 = \ker \d$, and so $\a=0$ as claimed.

By the commutativity of the $*$-product we may interchange
$\a$ and $\b$ above and conclude that the pairing is
non-degenerate as asserted.

Finally we recall that $\hH^k_{\infty}(X) \subset \hH^k(X)$ is a
closed subgroup with discrete quotient, and $\hH^k_{\infty}(X)
\cong \ce^k(X)/\cz^k_0(X) \cong 
d\ce^k(X)\oplus\{H^k(X;\,\bbr)/H^k_{\text{free}}(X;\,\bbz)\}$. 
Recall, by Lemma 5.1 that
$$
\Hom(d\ce^k(X), S^1) \ = \ d {\cd'}^{n-k-1}(X).
$$
It follows that the image of the homomorphism $\cd$
consists exactly of the smooth homomorphisms. \qed

\bigskip

\Remark {6.7} {\bf Evaluation homomorphisms lie at the edge of
the smooth dual.}  Note that for any rectifiable cycle $Z$ of
dimension $k$ there is a natural {\sl evaluation homomorphism} $$
h_Z : \hH^k(X) \ \arr \ S^1
\qquad\text{ given by }\ \  h_Z(\b) \equiv \Psi(\b)(Z)
$$
where $\Psi$ is the map to Cheeger-Simons characters defined
in \S 4. {\sl  These homomorphisms are not in the smooth dual of
$\hH^k(X)$, however they lie at the boundary of the smooth dual.
In fact elements of the smooth dual can all be considered as
generalized smoothing homotopies associated to these
homomorphisms.}  To see this let $a\in \cs^{n-k-1}(X)$ be  an
$\Lloc$-spark which satisfies the equation
$$
da\  = \ \psi - Z
$$
for a smooth $(n-k)$-form $\psi$.  Fix $\b\in \hH^k(X)$  and
choose an $\Lloc$-spark $b\in \b$ which is smooth on
a neighborhood of supp$(Z)$ with spark equation            
$db = \phi - R$.  Then by 3.2
$$\aligned
(\a * \b) \, [X]
\ &\equiv\ \int_X a\wedge \psi  + h_{Z}(\b) \ \ \ \
\left(\text{mod}\  \bbz \right)  \\ 
&\equiv \ \int_{R} a +  \int_X \phi\wedge b \ \ \ \
\left(\text{mod}\  \bbz \right)  
\endaligned
\tag6.8
$$
where $\a=\langle a \rangle$ and where we assume $a$ to be
smooth on supp$(R)$. 

Now suppose we take  a   family of Federer
smoothing homotopies $a_\e=H_{\e}(Z)$ as in 2.8. Then
$a_\e$ is an $\Lloc$-form of degree $n-k-1$ with support in
an $\e$-neighborhood of $Z$.  It satisfies the equation
$$
da_\e \ =\ Z_\e-Z
$$
where $Z_\e$ is a smooth $(n-k)$-form which is a local
smoothing of the cycle $Z$.  This family has the property
that 
$$
\lim_{\e\to  0} a_\e\ =\ 0\qquad\text{in the flat topology}
$$
and in particular
$$
\lim_{\e\to 0} Z_\e \ =\ Z
$$
Now $\supp (Z)\cap\supp( R)=\emptyset$  since $b$ is smooth
near $Z$.  Hence $\supp (a_\e)\cap\supp( R)=\emptyset$ for all
$\e$ sufficiently small.    With $\a_\e=\langle a_\e \rangle$,
equation  (6.8) becomes 
$$\aligned
(\a_\e * \b) \, [X]
\ &\equiv\ \int_X a_e\wedge \psi  + h_{Z}(\b) \ \ \ \
\left(\text{mod}\  \bbz \right)  \\ 
&\equiv \   \int_X Z_\e\wedge b \ \ \ \
\left(\text{mod}\  \bbz \right)  
\endaligned
\tag6.9
$$ 
and we see (in two ways) that as $\e\to 0$ the smooth
homomorphisms 
 $(a_\e*\bullet)[X]$
converge to  $h_Z$. Thus the evaluation homomorphisms
 lie at the boundary of the smooth ones.     

\vskip .3in 
 
 

\subheading{\S 7. Examples}   We examine the groups of
differential  characters on some standard spaces to illustrate
the duality established above.   Given a closed subspace $V$ of 
differential forms, let $V_{\infty}' = \Hom_{\infty}(V, S^1) =
\Hom_{\infty}(V, \bbr)$ denote the smooth Pontrjagin dual of $V$.  
We begin with the following basic observation.

\Lemma{7.1} {\sl  let $X$ be a compact oriented manifold of
dimension $n$. Then for each $k$,}
$$
\{d\ce^k(X)\}_{\infty}' \ =\   d\ce^{n-k-1}(X).
$$

\pf
By Lemma 5.1 we know that $(d\ce^k)' = d  {\cd'}^{n-k-1}
\subset{\cd'}^{n-k}$.
Hence, $(d\ce^k)_{\infty}' = (\ce^k)_{\infty}'\cap d{\cd'}^{n-k-1}
= \ce^{n-k}\cap d{\cd'}^{n-k-1} = d\ce^{n-k-1}$ since by the basic
results of de Rham any smooth form which is weakly exact is smoothly
exact.  \qed

\Ex {7.2.\ \ (Surfaces)}  Let $\Sigma$ be a compact oriented
surface of genus $g$.  The duality result asserts that
$\hH^1(\Sigma)$ is the smooth Pontrjagin dual of $\hH^0(\Sigma)$.
We see this explicitly as follows.  From 1.14 we have the (split)
short exact sequences
$$\aligned
0\to \frac{H^1(\Sigma;\,\bbr)}{H^1(\Sigma;\,\bbz)}  \to
&\hH^1(\Sigma) \to d\ce^1(\Sigma)\times \bbz\to 0    \\
0\to \frac{H^0(\Sigma;\,\bbr)}{H^0(\Sigma;\,\bbz)}  \to
&\hH^0(\Sigma) \to d\ce^0(\Sigma)\times \bbz^{2g}\to 0    
 \endaligned$$
from which we deduce that
$$
\hH^1(\Sigma) \cong (S^1)^{2g}\times \bbz\times d\ce^1(\Sigma)
\qquad\text{and}\qquad
\hH^0(\Sigma) \cong \bbz^{2g}\times S^1\times d\ce^0(\Sigma)
$$

\Ex {7.3.\ \ (3-manifolds)}  Let $M$ be a compact oriented
3-manifold.  Then by Observation 5.3,  $\hH^1(M) \cong
\{H^1(M;\,\bbr)/H^1(M;\,\bbz)_{\text{free}}\}
\times d\ce^1\times H^2(M;\,\bbz)$
which using Poincar\'e duality can be rewritten
$$
\hH^1(M) \ \cong\ \Hom(H_1(M;\,\bbz)_{\text{free}}, S^1)  \times
H_1(M;\,\bbz)_{\text{free}} \times d\ce^1(M) \times
H_1(M;\,\bbz)_{\text{torsion}} 
$$
where (via 7.1) the self-duality is manifest.

\Ex {7.4.\ \ (Complex Projective Space $\bbp^n_{\bbc}$)}  
If $X=\bbp^n_{\bbc}$, we see that 
$$
\frac{H^k(X;\,\bbr)}{H^k(X;\,\bbz)}\ =\ \cases
S^1 \qquad &\text{if $k$ is even}\\
0 \qquad &\text{if $k$ is odd}.\endcases
$$
and
$$
 \cz_0^{k+1}(X) \ = \ 
\cases
d\ce^k(X) \qquad &\text{if $k$ is even}\\
\bbz\times d\ce^k(X) \qquad &\text{if $k$ is odd}.
\endcases
$$
From  1.14 we then compute the following:

\medskip
$$
\matrix
k  & \ & \hH^k(X) \\
\ & \ &\ \\
-1   & \ &   \bbz \\
0    & \ &    S^1\times d\ce^0   \\
1    & \ &    \bbz\times d\ce^1   \\
2    & \ &    S^1\times d\ce^2   \\
3    & \ &    \bbz\times d\ce^3   \\
\ & \ &\ \\
\cdot  & \ &  \cdot\\
\cdot  & \ &  \cdot\\
\cdot  & \ &  \cdot\\
\ & \ &\ \\
2n-3    & \ &    \bbz\times d\ce^{2n-3}   \\
2n-2    & \ &    S^1\times d\ce^{2n-2}   \\
2n-1    & \ &    \bbz\times d\ce^{2n-1}   \\
2n    & \ &    S^1  
\endmatrix
$$

\medskip\noindent
{\bf Example 7.5.\ \ (Real Projective Space $\bbp^{2n+1}_{\bbr}$)}
  Proceeding as above we calculate that when 
$X=\bbp^{2n+1}_{\bbr}$, one has
$$
\matrix
k  & \ & \hH^k(X) \\
\ & \ &\ \\
-1   & \ &   \bbz \\
0    & \ &    S^1\times d\ce^0   \\
1    & \ &    \bbz_2\times d\ce^1   \\
2    & \ &     d\ce^2   \\
3    & \ &    \bbz_2\times d\ce^3   \\
4    & \ &     d\ce^4   \\
5    & \ &    \bbz_2\times d\ce^5   \\
6    & \ &     d\ce^6   \\
\ & \ &\ \\
\cdot  & \ &  \cdot\\
\cdot  & \ &  \cdot\\
\cdot  & \ &  \cdot\\
\ & \ &\ \\
2n-4   & \ &     d\ce^{2n-4}   \\
2n-3    & \ &    \bbz_2\times d\ce^{2n-3}   \\
2n-2   & \ &     d\ce^{2n-2}   \\
2n-1    & \ &    \bbz_2\times d\ce^{2n-1}   \\
2n-0     & \ &    \bbz\times d\ce^{2n}   \\
2n+1   & \ &    S^1  
\endmatrix
$$
Recall that $\Hom(\bbz_2, S^1)= \bbz_2$.

\Ex {7.6.\ \ (Products of projective
spaces)}  When $X=\bbp^{2}_{\bbc}\times \bbp^{2}_{\bbc}$, we
find that:
$$
\matrix
k  & \ & \hH^k(X) \\
\ & \ &\ \\
-1   & \ &   \bbz \\
0    & \ &    S^1\times d\ce^0   \\
1    & \ &    \bbz\times \bbz\times d\ce^1   \\
2    & \ &      S^1\times S^1\times d\ce^2   \\
3    & \ &    \bbz\times \bbz\times\bbz \times d\ce^3   \\
4    & \ &     S^1\times S^1\times S^1\times d\ce^4   \\
5    & \ &    \bbz\times\bbz\times d\ce^5   \\
6    & \ &       S^1\times S^1\times d\ce^6   \\
7      & \ &    \bbz\times d\ce^{7}   \\
8   & \ &    S^1  
\endmatrix
$$

\bigskip

     \vskip .3in

 

\subheading{\S 8. Characters with compact support}  The theory
developed in the previous sections can be carried through for
characters with compact support.  This represents the
``homology version'' of the theory.  Let $X$ be a smooth
$n$-manifold as before and let ${\cd'}^k_{\cpt}(X)
\subset {\cd'}^k(X)$ denote the {\bf currents with compact
support} on $X$.  For each of the subspaces 
$W\subset {\cd'}^k(X)$ introduced in \S 1, we denote
$$
W_{\cpt}\ \equdef \ W \cap {\cd'}^k_{\cpt}(X).
$$
The analogues of  statements (1.1) hold for currents
with compact support, and one has

\Def{8.1}\ {\bf (Sparks with compact support)}
$$
\cs^k_{\cpt}(X) \ \equdef \ \{a\in {\cd'}^k_{\cpt}(X) \ :\ 
da = \phi - R \  \text{where}  \ \phi \in \ce^{k+1}_{\cpt}(X)
\ \text{and}  \ R \in \ci\cf^{k+1}_{\cpt}(X)\}
$$

\Def{8.2}  The space
$$
\hH^k_{\cpt}(X) \ \equdef \ \cs^k_{\cpt}(X) / 
\{ d{\cd'}^{k-1}_{\cpt}(X) + \ci\cf^k_{\cpt}(X) \}
$$
will be called the group of {\bf
de Rham - Federer characters with compact support} on $X$.

\bigskip

The analogues of Propositions 1.11, 1.12, 1.14 hold for
$\hH^k_{\cpt}(X)$.  In particular we have

\Prop{8.3}  {\sl There is a functorial short exact sequence
$$
0\ \arr \ \frac{H_{n-k}(X;\, \bbr)}{H_{n-k}(X;\,
\bbz)_{\oper{free}}} \ \arr\ \hH^k_{\cpt}(X) \ @>{\d}>>
\Q^{k+1}_{\cpt}(X) \ \arr\ 0
$$
where 
$$
\Q^{\ell}_{\cpt}(X) \ \equdef\  \left\{(\phi, u) \in
\cz^{\ell}_{0,\cpt}(X) \times H_{n-\ell}(X;\,\bbz) \ : \ [\phi]
= u\otimes \bbr\right\}
$$
}

The analogues of the results in \S 2 and \S 3 carry over to
these spaces.

Arguments parallel to those of \S 4 establish the following.
We shall say that a Cheeger-Simons differential character
$a:\cz_k(X)\to \bbr/\bbz$ has {\bf compact support} if there is
a compact subset $K\subset X$ such that $a(Z) = 0$ for all
cycles  $Z$ with supp$(Z) \cap K = \emptyset$.

\Theorem{8.4}  {\sl Let  $\hh^k_{\cpt}(X; \bbr/\bbz)
\subset \hh^k(X; \bbr/\bbz)$ denote the  subgroup of
differential  characters of degree $k$ with compact
support on $X$. Then there is a  natural isomorphism   
$$
\Psi : \hH^k_{\cpt}(X) @>{\cong}>> \hh^k_{\cpt}(X; \bbr/\bbz)
$$
induced by integration.} 

\bigskip

The definition of the $*$-product given in \S 3 extends
straightforwardly to give a pairing
$$
* : \hH^k(X) \times \hH^{\ell}_{\cpt}(X)
\ \arr \ \hH^{k+\ell+1}_{\cpt}(X)
\tag{8.5}
$$
which in particular makes $\hH^*_{\cpt}(X)$ a graded
commutative ring.

Suppose now that $X$ is  connected, oriented, and of
dimension $n$.  
\Def{8.6}  For
each integer $k$, $0\leq k<n$ we define  duality mappings
$$
\cd : \hH^{n-k-1}_{\cpt}(X) \ \arr\ \hH^{k}(X)^*\equiv
\Hom_{\infty}\left(\hH^{k}(X), \bbr/\bbz\right)
$$
and
$$
\cd : \hH^{n-k-1}(X) \ \arr\ \hH^{k}_{\cpt}(X)^*\equiv
\Hom_{\infty}\left(\hH^{k}_{\cpt}(X), \bbr/\bbz\right)
$$
by
$$
\cd(\a)(\b) \ \equdef \ \a *\b \in \hH^n_{\cpt}(X) = \bbr/\bbz.
$$
where $\Hom_\infty(\bullet, \bbr/\bbz)$ denotes the smooth
homomorphisms, defined as in 5.10.\medskip

\Theorem {8.7. \ (Poincar\'e-Pontrjagin Duality  on Non-compact
Manifolds)} {\sl For any connected, oriented manifold $X$ the
duality mappings above are  isomorphisms.}

\pf  The argument follows directly the second proof given in \S
6. Full details appear in [HL$_6$]. \qed

      \vskip .3in

 

\subheading{\S 9. The Thom homomorphism}  An interesting feature of
the de Rham theory is that one easily constructs Thom isomorphisms 
for characters.  

Let $X$ be a  smooth $n$-manifold and
$$
\pi : E \arr X
$$
a smooth oriented riemannian vector bundle of rank $d$.  Once and
for all we fix a  {\sl compact approximation mode} in the sense of
[HL$_1$; \S I.4].  This amounts to choosing a $C^{\infty}$-function
$\chi:[0,\infty] \to [0,1]$ with $\chi(0) = 0$, $\chi' \geq 0$ and
$\chi(t) \equiv 1$ for  $t\geq 1$.
Then from [HL$_1$; \S\S IV1-2] we have the following. \medskip\ 

\Theorem {9.1.   (Harvey-Lawson [HL$_1$; \S\S IV1-2]) } {\sl
To each orthogonal connection $D_E$ on $E$ there is a canonically
associated {\bf Thom form } $\tau \in  \ce^d(E)$ such that
\roster
\item $d\tau = 0,$
\item $\oper{supp}(\tau) \subset D(E) =  \{e\in E : \|e\|\leq 1\},$
\item $\pi_* (\tau) = 1$
\item $\tau\bigr|_X = \chi(D_E)$
\endroster
where $\chi(D_E)$ is the Chern-Euler form of $D_E$
when $d$ is even and is 0 when $d$ is odd. 

There is furthermore a canonically associated {\bf Thom spark}
$\g \in \ce^{d-1}_{\Lloc}$ which is smooth outside the zero
section $X\subset E$, has support in the unit disk bundle  $D(E)$
and satisfies the spark equation
$$
d\g \ = \ \tau - [X] \qquad \ \ \text{on}\ \ E.
$$
}

\Remark{9.2}
For $d$ even, $\tau$ is  the Chern-Euler
characteristic form of the time-1 pushforward connection on
$\pi^*E$ associated to the tautological cross-section. 

\Remark {9.3}  
Since $\pi_* : \ce^{k+d}_{\cpt}(E) \arr \ce^k_{\cpt}(X)$
(integration over the fibre), there is a pull-back map on
general currents 
$$
\pi^* : {\cd'}^k(X) \ \arr\ {\cd'}^{k+d}(E).
$$
This map preserves the properties of local integrability, local
rectifiability,  local flatness and local integral flatness.

\medskip

As a consequence we have an induced  homomorphism
$$
\pi^* : \hH^*(X) \arr \hH^{*+d}(E).
$$
We would like the corresponding mapping for the functor
$\hH^*_{\cpt}$.

\Def{9.4}  We define a {\bf Thom mapping}
$$
\ct : \cs^k_{\cpt}(X) \ \arr\ \cs^{k+d}_{\cpt}(E)
$$
as follows.  For $a\in \cs^k_{\cpt}(X)$ write $da = \phi-R$
as in 1.5.  Then
$$
\ct(a) \ \equdef \ \pi^*a \wedge \tau \,+\,
(-1)^{k+d+1}\pi^*(R)\wedge\g 
$$
where $\tau$ and $\g$ are the Thom form and the Thom spark
associated to the connection in Theorem 9.1. Under a
local orthogonal trivialization $\pi^{-1}U \cong U\times \bbr^d$ 
of $E$ we have $\pi^*R \cong R \times [\bbr^d]$ and $\g =
\text{pr}^*\g_0$ where $\g_0$ is the universal Thom spark on 
$\bbr^d$.  Thus $\pi^*(R)\wedge\g$ is well-defined.
Note that 
$$
d \ct(a) \ =\ \pi^*\phi\wedge \tau \,-\, \pi^*R\wedge [X].
\tag{9.5}
$$
\medskip\ 

\Theorem{9.6} {\sl Let $E\arr X$ be as above.  Then for
all $k\geq 0$ the Thom mapping $\ct$ induces an injective
homomorphism 
$$
\bbthh : \hH^k_{\cpt}(X) \ \arr\ \hH^{k+d}_{\cpt}(E)
$$
with the following properties:
$$\aligned
(i)&\ \ \ \ \bbthh(1) \ =\ [\g] 
\qquad\qquad\qquad\qquad \\
(ii) &\ \ \ \     \bbthh(\a) \ =\ (\pi^*\a) *  \bbthh(1) 
                    \ =\ (\pi^*\a) *  [\g]   
\qquad\qquad\qquad\qquad \\
(iii) &\ \ \ \    \d_1\bbthh(\a) \ =\ \bbth_{\oper{deR}}(\d_1\a) 
 \ \equdef\ \pi^*(\d_1\a) \wedge \omega\qquad \text{and}
\qquad\qquad\qquad\qquad \\
& \ \ \ \ 
\d_2\bbthh(\a) \ =\ \bbth(\d_2\a) \endaligned
$$
where $\bbth$ denotes the Thom isomorphism on integral cohomology.}

\Note{9.7} The unit $1 \in \hH^{-1}(X)$ is the deRham character 
$\g_0$ with 
$$
d\g_0 = 1 - [X] \qquad\qquad\text{on} \ \ X.
$$
Of course, $1 = [X]$ and $\g_0 = 0$ as currents. Think of $\g_0$
as the identity spark.  Part (i) says that  the image of  
$\g_0$ is represented by the  Thom spark $\g$ which
satisfies 
$$
d\g = \tau - [X] \qquad\qquad\text{on} \ \ E.
$$ 
The class $[\g]\in \hH^{d-1}_{\cpt}(E)$ will be called the 
{\bf Thom character}.

\Note{9.8}  Part (ii) states that the Thom homomorphism for
characters is a $\hH^{*}(X)$-module mapping.  The image
is the free rank-1 module generated by the  Thom character.

\Note{9.9} Part (iii) states that the differential $\d$ carries 
the Thom homomorphism on characters over to the classical Thom 
isomorphisms.

\bigskip\noindent
{\bf Proof of Theorem 9.6}\ 
We have directly from  3.2 and 9.4 that
$$
\ct(1) \ =\ \g
\qquad\text{and}\qquad
\ct(a) \ =\ a * \g.
$$
which proves (i) and (ii).  Part  (iii) follows immediately from
(9.5). It remains to prove that $\bbthh$ is injective.
For this we observe that $\bbthh$ induces a map of the short exact
sequences  8.3 and 1.14:
$$
\CD
0  @>>>   \frac{H_{n-k}(X;\, \bbr)}{H_{n-k}(X;\,
\bbz)_{\oper{free}}}  @>>>   \hH^k(X) \ @>{\d}>>
Q^{k+1}(X) @>>>   0  \\
    @.     @V{\cong}VV         @VV{\bbthh}V  
@VV{\bbth_{\text{deR}}\times\bbth}V    \\ 
0  @>>>  \frac{H_{n-k}(E;\, \bbr)}{H_{n-k}(E;\,
\bbz)_{\oper{free}}}  @>>>   \hH^{k+d}_{\cpt}(E) \ @>{\d}>>
Q^{k+d+1}_{\cpt}(E) @>>>   0.   \\ \endCD
$$
The map on the left is a quotient of Thom isomorphisms under 
Poincar\'e duality.  The map  on the right is a product 
of $\bbth_{\text{deR}}(\phi) = \pi^*(\phi)\wedge\tau$, which is
clearly injective, and the Thom isomorphism $\bbth$ on integral
cohomology. It follows that $\bbthh$ is injective.
\qed

\Remark{9.10. (Atomic sections)} A smooth cross-section $s:X\to E$
is called {\bf atomic} if, whenever it is written in a
local frame $e_1,...,e_d$ for $E$ as $s =  u_1e_1+\dots+u_d e_d$,
it  has the property that $du^I/|u|^{|I|} \in L^1_{\text{loc}}$ for
all $|I|<d$. Associated to an atomic section $s$ is  a
codimension-$d$ integrally flat current  $\Div(s)$ defined by the
vanishing  of $s$  [HS].  Any section $s$ with non-degenerate zeros
is atomic, and $\Div(s)$ is the manifold defined by $s=0$.

If $s:X\to E$ is atomic and $d$ is even, then $s^*\tau =
\chi(D_s)$ is the Euler form of the time-1 push-forward
connection $D_s$ associated to $s$ (cf. 9.2). Furthermore, $s^*\g$
is the associated Euler spark as defined in [HZ$_2$]. This is a
canonical $\Lloc$-form on $X$ which satisfies $$ d s^*\g \ = \
\chi(D_s) - \Div(s). $$
Note that in general 
$$
s^*\ct(a) \ = \ a * (s^*\g).
\tag{9.11}
$$

\medskip

We now examine how $\bbthh$ changes under a change of connection
on $E$.  Suppose $D_0$ and $D_1$ are two orthogonal connections on
$E$ with Thom sparks $\g_0$ and $\g_1$ satisfying
$$
 d\g_0\ = \ \tau_0 - [X]\qquad\ \ \text{and}\qquad\ \ 
 d\g_1\ = \ \tau_1 - [X]
$$
Consider the convex family of connections $D_t = (1-t)D_0+tD_1$,
and define a connection $D$ on the bundle  $E\times \bbr
\arr X\times \bbr$   by $D = D_t + dt\otimes (\partial/\partial
t)$.  Let $\g$ and $\tau$ be the Thom spark and Thom form
associated to $D$. They satisfy the equation
$$
d\g\ =\ \tau-[X\times \bbr]
$$  
Consider the flat current $\b$ and the
$C^{\infty}$-form $\g_{01}$ defined by 
$$
\b \ \equdef\ \oper{pr}_*\left( \psi \g\right)
\qquad\   \text{and}\qquad\   
\g_{10}\ \equdef\ -\oper{pr}_*\left( \psi \tau\right)
$$
where $\psi$ is the characteristic function of 
$E\times [0,1]\subset E\times \bbr$  and $\oper{pr}: E\times \bbr
\to E$ is the projection. Since $\oper{pr}_*(\psi[X\times
\bbr])=0$, we see that 
$d\b = \oper{pr}_*(d(\psi \g)) = \g_1 -\g_0
+\oper{pr}_*(\psi\tau)$, and so  
$$
\g_1-\g_0 \ = \ \g_{10} + d\b.
\tag{9.12}
$$
Taking $d$ gives the equation
$$
\tau_1-\tau_0 \ = \   d\g_{10}.
\tag{9.13}
$$

\Prop{9.14}  {\sl  Let $\bbthh_k$ be the Thom homomorphism on
differential characters associated to the connection $D_k$ on $E$
for $k=0,1$.  Then 
$$
\bbthh_1(\a) - \bbthh_0(\a) \ = \ \pi^*(\a) * [\g_{10}]
$$
where $[\g_{10}]$ denotes the character defined by the smooth
transgression form $\g_{10}$.}

\pf
From (9.12) we have
$$
\aligned
\qquad\bbthh_1(\a) - \bbthh_0(\a) \ &= \ 
\pi^*(\a) * [\g_1] - \pi^*(\a) * [\g_0]\\ 
&= \  \pi^*(\a) * ([\g_{10} +d\b])\ = \ \pi^*(\a) * ([\g_{10}]).
\qquad\qquad\square  \endaligned
$$      

\vskip .3in

 

\subheading{\S 10. Gysin maps} One advantage of the de Rham
point of view is that it is easy to define Gysin maps.  For
example if $f:X\to Y$ is a smooth fibre bundle with compact
oriented fibres of dimension $d$, then push-forward of currents
defines   ``Gysin'' homomorphisms
$$
f_! : \hH^{k+d}(X) \ \arr\ \hH^{k}(Y)
\qquad\text{and}\qquad
f_! : \hHc^{k+d}(X) \ \arr\ \hHc^{k}(Y)
\tag{10.1}
$$
This raises the question of whether there exist Gysin maps
of differential characters for more general mappings, and in
particular for smooth embeddings.  Such maps do exist --  the Thom 
homomorphism of \S 9 is a basic example.  It depends on a
choice of connection.  In general  the 
Gysin homomorphism  depends   on a choice of ``normal
geometry'' to the embedding.

\Def{10.2} Suppose $X$ and $Y$ are   manifolds and 
$f:X\hookrightarrow Y$ is a smooth embedding with oriented normal
bundle $N$.  Given an orthogonal connection on N and a smooth
embedding $F:N\hookrightarrow Y$ extending $f$, there is a Gysin
homomorphism 
$$
f_! : \hH^*_{\cpt}(X) \ \arr\ \hH^{*+d}_{\cpt}(Y),
$$
where $d= \dim(Y)-\dim(X)$,
defined by 
$$
f_! \ = \ F_*\circ \bbthh
$$
where $\bbthh$ is the Thom homomorphism for $N$ defined in \S 9.

\Ex{10.3}  Suppose $X$ and $Y$ are riemannian and  
$f:X\hookrightarrow Y$ is an isometric embedding.
Let $\rho_f$ be the injectivity radius of $f$, i.e., the sup of
the numbers $\rho$ such that the exponential mapping 
exp$:N\to Y$ is injective on 
$D_{\rho}(N)=\{v\in N : \|v\|\leq \rho\}$.
Then for any $\rho<\rho_f$ the mapping $F:N\to Y$ given by 
$F= \oper{exp}\circ\mu_{\rho}$ is an embedding on a neighborhood
of $D_1(N)$ and since the Thom spark of $N$ has support in
$D_1(N)$ the composition  $f_! \ = \ F_*\circ \bbthh$ is well
defined.

\medskip\noindent
{\bf Definition 10.4.} Let $f:X\to Y$ be {\sl any} smooth mapping
between  compact oriented manifolds.  Assume Y is 
riemannian with injectivity radius $i_Y$.  Consider the smooth
embedding  $\gamma_f : X \hookrightarrow X\times Y$ defined by
graphing   $\gamma_f(x) = (x, f(x))$ for $x\in X$.
Now the normal bundle $N$ to $\gamma_f$ can be identified with
the restriction of $TY$ to $\gamma_f(X)$, and we have  a {\bf
canonical} identification of $\{v\in N : \|v\|<i_Y\}$ with a
tubular neighborhood  of $\gamma_f(x)\subset X\times Y$ which
associates to  $v \in T_{f(x)}Y$ the point $\exp_{f(x)}(v)$. Thus
as above we have a canonically defined Gysin mapping
$(\gamma_f)_!$.   Following Grothendieck we combine this with the
projection $\oper{pr}:X\times Y \to Y$ to get a Gysin mapping
$$
f_! \ \equdef \ \oper{pr}_!\circ \left(\gamma_f \right)_!
$$

\Remark {10.5}  If $X$ and $Y$ are compact and oriented,
then Gysin maps for ordinary cohomology  can be defined using
Poincar\'e duality $D:H^k(X;\,\bbz)\to H_{n-k}(X;\,\bbz)$  by
setting $f_! = D^{-1}\circ f_*\circ D$. This does not quite work
for differential characters because the inverse $\cd^{-1}$ to the
duality map is only densely defined.   However  if one takes the
family of Gysin maps $f^t_!$, $0<t\leq 1$ obtained  via the  
family of Thom forms $\tau_t = \rho_t^*\tau$ (where $\rho_t$ is
scalar multiplication by $t$ in the normal bundle), then these maps
converge to the ``non-smooth character'' $ \cd^{-1}\circ f_*\circ
\cd$, i.e.,  
$$
\lim_{t\to 0}\cd\circ f^t_!\ =\ f_*\circ \cd
$$
on $\hH^*(X)$

      \vskip .3in


\subheading{\S 11. Morse sparks}  Sparks  represent
differential characters just as forms (or currents) 
represent cohomology, and they arise naturally in many situations.
Interestingly a rich source of sparks is provided by  Morse Theory.
In [HL$_5$] it is shown that for each Morse function $f$ on
$X$ there is a  subcomplex $\ce_{0,f}^*(X) \subset\ce^*(X)$ of finite
codimension with a surjection $P:\ce_{0,f}^*(X) \to \bbs_f^{\bbz}$ 
onto the {\sl integral}  Morse complex, which has affine fibres and
induces an isomorphism in cohomology. Furthermore, if  $\cz_{0,f}^*(X)
\subset\cz_{0}^*(X)$ denotes the closed forms in  $\ce_{0,f}^*(X)$,
then there is a continuous linear mapping $T:\cz_{0,f}^*(X) \to
S^{*-1}(X)$ into the space of sparks with $d(T\phi) = \phi-P(\phi)$.
This inverts the differential $\d_1$ and canonically splits the
characters on this set. We present below a summary  of this
construction.

Let $f:X\to \bbr$ be
a Morse function on a compact manifold and choose a gradient-like flow
$\vf_t:X\to X$ satisfying the Morse-Stokes Axioms in [HL$_5$] (such
flows always exist).  Let $\text{Cr}_f$ denote the set of critical
points of $f$ and for each $p\in \text{Cr}_f$ let  $S_p$, $U_p$ denote
the stable and unstable manifolds of $p$ respectively. These are
manifolds of finite volume in $X$ and so define integral currents
$[S_p]$ and $[U_p]$.  Consider the finite-dimensional vector space and
the integral lattice inside it given by $$
\bbs_f \ =\ \text{span}_{\bbr}\left\{[S_p]\right\}_{p\in  \text{Cr}_f}
\and
\bbs_f^{\bbz} \ =\ \text{span}_{\bbz}\left\{[S_p]\right\}_{p\in 
\text{Cr}_f}.
$$
Each of these spaces of currents is $d$ invariant and there are
natural isomorphisms 
$$
H^*(\bbs_f) \cong H^*(X;\,\bbr)
\and
H^*(\bbs_f^{\bbz}) \cong H^*(X;\,\bbz)
$$
In fact the   continuous linear projection
$$
P:\ce^*(X) \ \arr\ \bbs_f
$$
given by 
$$
P(\phi) = \sum_{p\in  \text{Cr}_f}  \left\{\int_{U_p}\phi\right\}
[S_p]
\tag11.1
 $$
is chain homotopic to the identity. That is, there exists a continuous
linear mapping $T:\ce^*(X)\to {\cd'}^*(X)$  of degree -1, such that
$$
d\circ T+T\circ d\ =\ I-P.
\tag11.2
$$
where $I:\ce^*(X)\to {\cd'}^*(X)$ is the obvious inclusion.

Let us denote by $Z\bbs_f$ and $Z\bbs_f^{\bbz}$ the cycles in  
$\bbs_f$ and $\bbs_f^{\bbz}$ respectively.

\Def{11.3} Fix an integral cycle $R=\sum n_p[S_p] \in Z\bbs_f^{\bbz}$. A
smooth form $\phi \in \cz^k_0(X)$ is called a  {\bf Thom form  for} $R$
if $P(\phi)=R$, i.e., if $\int_{U_p}\phi = n_p$ for all $p$.
Let $\cz_{0,f}^*(X) = \cz_{0}^*(X)\cap P^{-1}(Z\bbs_f^{\bbz})$ denote
the set of all such forms. We now show that Thom forms exists for every
cycle in the Morse complex.

\Lemma{11.4}{\sl   The restricted mapping 
$P:\cz_{0,f}^*(X)\arr Z\bbs_f^{\bbz}$ is surjective. Its fibres are
affine subspaces of finite codimension in $\cz_{0}^*(X)$.
}

\pf  For a Morse-Stokes flow we have that $S_p$
intersects $U_p$ transversely in one point and that 
$\overline{S_p}\cap \overline{U_{p'}}=\emptyset$ for all distinct pairs
$p,p'\in\text{Cr}_f$ of the same index.  Fix a cycle $R\in
Z\bbs_f^{\bbz}$ of degree $k$, and denote by $H_{\e}(R)$ the Federer
homotopy smoothing supported in an $\e$-neighborhood of $R$ (See 2.8).
Since supp$(R)\cap \partial U_{p}=\emptyset$ for all $p$ of index
$n-k$, for $\e$ sufficiently small we also have supp$(H_{\e}(R))\cap
\partial U_{p}=\emptyset$ for all such $p$. Now $d H_{\e}(R) =
R_{\e}-R$ where $R_{\e}$ is a smooth $d$-closed $k$-form. Now
$$\aligned
\int_{U_p}R_{\e} &= (R\wedge [U_p])[X] + (dH_{\e}(R)\wedge[U_p])[X]\\
&= \sum_{p'} n_{p'} ([S_{p'}] \wedge [U_p])[X]
+(H_{\e}(R)\wedge\partial[U_p])[X]  = n_p+0
\endaligned
$$
as desired.  Note that $P^{-1}(R)$ is defined by the affine equations
$\int_{U_p}\phi = n_p$ for all $p\in \text{Cr}_f$ whereas 
$\cz^k_0(X)$ is defined by the weaker conditions
$$
\sum_{p\in \text{Cr}_f} m_p\int_{U_p}\phi\ \in \ \bbz
\qquad \text{whenever  } d\left(\sum m_p{U_p}\right)=0
\text{ and all } m_p\in \bbz.
$$
This proves the second assertion of the lemma. \qed

\medskip  
 
Observe now that there is a commutative diagram
$$
\CD
\cz_{0,f}^*(X) @>{\subset}>>   \cz_{0}^*(X)  \\
@V{P}VV  @VV{P_0}V  \\
Z\bbs_f^{\bbz} @>{\text{pr}_0}>> H^*_{{\text{free}}}(X;\,\bbz)
\endCD
\tag11.5
$$
where $P_0:\cz_0^*(X)\to \Hom(H_*(X;\,\bbz),\, \bbz) \cong
 H^*_{{\text{free}}}(X;\,\bbz)$ is the {\bf period mapping}
and $\text{pr}_0$ is the obvious projection. Both $P$ and $P_0$ are surjective
with affine fibres.  They represent a ``fattening'' by smooth forms of
 free abelian groups  $Z\bbs_f^{\bbz}$
and $H^*_{{\text{free}}}(X;\,\bbz)$.    Note
that $\cz_{0,f}^*(X) =   \cz_{0}^*(X)$  if and only if 
${\bbs_f}=Z\bbs_f$, i.e., iff $f$ is a perfect Morse function and 
$H^*(X;\,\bbz)$ is torsion-free.

The map $T$ from (11.2) associates a Morse spark to each Thom form on
$X$. This gives the following.

\Theorem{11.6} [HL$_5$] {\sl  Let $f:X\to\bbr$ be a Morse
function as above.   Then there is a continuous linear operator of
degree -1 
$$
 T: \cz^*_{0,f}(X) \ \arr\ \cs^{*-1}(X)
$$
with the property that for each $\phi\in \cz^*_{0,f}(X)$
$$
d\{T(\phi)\}\ =\ \phi - \sum_{p\in\text{Cr}\,(f)} n_p \,[S_p]
\qquad
\text{ where }
\qquad
n_p = \int_{U_p}\,\phi.
$$
}

\noindent
{\bf Remark 11.7.} Note that $T$ is a right inverse of the mapping
$d_1:\cs^{*-1}(X)\to \cz^*_0(X)$ on the subset $\cz^*_{0,f}(X)$, and
so it splits the fundamental sequence 1.12  on this set.  To be more
explicit we define the subgroup of {\bf $f$-characters}: 
$$
\hH^*_f(X) \ \equiv \d_1^{-1}\left\{\cz^*_{0,f}(X)\right\}\ \subset\
\hH^*(X) 
$$
This is a union of subspaces of constant finite codimension in
$\hH^*(X)$. The map $T$ above gives a canonical
splitting
$$
\hH^*_f(X) \ =\ H^*(X;\, S^1) \times \cz^*_{0,f}(X)
$$
under which $\d_1$ becomes projection to the second factor and
$\d_2=\text{pr}\circ P\circ \d_1$  where $\text{pr}:Z\bbs_f^{\bbz} \to
H^*(X;\,\bbz)$  is the projection.

      \vskip .3in


\def\wh{\widehat}
\def\supth{{\text{th}}}
\def\picc{\operatorname{Pic}}
\def\har{\operatorname{Har}}

\subheading{\S 12. Hodge sparks and  riemannian Abel-Jacobi
mappings} Another important source of sparks is provided by
Hodge theory. Suppose $X$ is a compact riemannian manifold. 
Recall (cf. [HP]) that any current $R$ on $X$, not just the
$L^2$-forms,  has a {\bf Hodge decomposition}
$$
R\ =\ H(R)+dd^*G(R)+d^*dG(R)
\tag12.1
$$
where $H$ is harmonic
projection  and $G$ is the Greens operator.  
Also recall that $d$ commutes with $G$, so that
if $R$ is a cycle, then $dG(R)=0$.

\Def{12.2}  Given an integrally flat cycle $R\in
\ci\cf^{k+1}(X)$, we define its {\bf Hodge spark}
$\s(R)\in\cs^k(X)$ by
$$
\s(R)\ \equiv\ -d^*G(R)
$$
and observe that it satisfies the spark equation
$$
d\s(R)\ =\ H(R)-R
\tag12.3
$$
Note that  $d_1\circ \s =H$ and $d_2 \circ \s = $ Id on the
space of integrally flat cycles, where $d_1$ and $d_2$ are the
differentials introduced in \S 1. In particular $\s$ is a
right inverse to  $d_2$.

Let $\wh\s(R)\in\hH^k(X)$ denote the differential character
corresponding to the Hodge spark $\s(R)$. Define the 
{\bf $k^\supth$ character Jacobian} of $X$ to be the torus
\def\Jac{\operatorname{Jac}}
$$
\Jac^k(X) \ \equiv\ \ker\d\ \cong\ 
H^k(X;\,\bbr)/H^k_{\text{free}}(X;\,\bbz)
$$
contained in
$\hH^k(X)$. If we restrict $\wh\s$ to the space
$\calb^{k+1}(X)= d \ci\cf^k(X) = d\calr^k(X)$ of integrally flat
boundaries then
$$
\wh\s: \calb^{k+1}(X)\ \to\ \Jac^k(X).\tag12.4
$$
 To see that $R=d\G$ implies that $\d_1\wh\s(R)=0$ note
that $H(R)=0$ and apply (12.3). Obviously
$\d_2(\wh\s(R))=0$. The map $\wh\s$ will be called the {\bf character
Abel-Jacobi map}.

\Def{12.5}  An integrally flat boundary
$R=d\G\in\calb^{k+1}(X)$ is {\bf linearly equivalent to zero}
or a {\bf principal boundary} if $\wh\s(R)=0$. Let
$\calb^{k+1}_P(X)$ denote the space of principal boundaries, and
define  the $k^\supth$ {\bf Picard space of $X$} to be
$\picc^k(X)=\calb^{k+1}(X)/\calb^{k+1}_P(X)$.

In summary, $\wh\s$ induces an injection (also denoted
$\wh\s$) on the quotient $\picc^k(X)$. That is
$$
\picc^k(X) \ \overset\hat\s\to{\hookrightarrow} \ 
\Jac^k(X).\tag12.6
$$
The character Jacobian is independent of the Riemannian
metric on $X$. However, fixing a metric, and letting $\har^k(X)$
and $\har^k_0(X)$ denote the spaces of harmonic $k$ forms and
harmonic $k$ forms with integral periods respectively, we have
isomorphisms
$$
\Jac^k(X)\cong\har^k(X)/\har_0^k(X)\cong\har^{n-k}(X)^*/H_{n-k}(X,\BZ)\cong
\Hom(\har^{n-k}_0,\BR/\BZ).
$$

\Def{12.8} The map
$$
\widehat\s_{r} :\calb^{k+1}(X)\ \longrightarrow\ 
\Hom(\har^{n-k}_0,\BR/\BZ)
$$
induced by $\wh\s$ will be called the $k^\supth$ {\bf
Riemannian Abel-Jacobi map}.

In order to compute $\widehat\s_{r}$  more explicitly we must
find a harmonic spark which is equivalent to the spark $\s(R)$.

\Lemma{12.9}  {\sl Suppose $R=d\G$ where $\G$ is integrally
flat. Then the Hodge decomposition of $\G$  has the form
$$
\G \ =\ H(\G)+dB-\s(R).
\tag12.10
$$
where $B=d^*G(\G)$}

\medskip\noindent
{\bf Proof.}\ \ Note that
$d^*dG(\G)=d^*G(d \G  )=-\s(R)$. \qed

\Cor{12.11}  {\sl  The sparks $\s(R)$ and $H(\G)$ determine
the same differential character, i.e.,
$$
\widehat\s(R) = \widehat H(\G) \in \hH^k(X)
$$
In particular, $\widehat H(\G)$ is independent of the choice of
integrally flat $\G$ with $d\G=R$.}
\medskip

Consequently, the mapping from $\calb^{k+1}(X)$ to
$\har^k(X)/\har^k_0(X)$ induced by $\wh\s$ is given by sending
$R=d\G\in\calb^{k+1}(X)$ to $H(\G)\in\har^k(X)$. Utilizing
the duality theorem between $\har^k(X)$ and $\har^{n-k}(X)$
we see that:
$$
\widehat\s_{r}(R)(\theta)  \ \equiv\   \int H(\G)\land\theta
\ \ \bmod\BZ\tag12.12
$$
for all $\theta\in\har^{n-k}_0(X)$.

Since $\theta$ is $d$ and $d^*$ closed, equation (12.10)
implies that
$$
\widehat\s_{r}(R)(\theta) \ \equiv\ 
\int_\G\theta\mod\BZ,\quad\text{for all }\ \theta\in\har^{n-k}_0(X).
$$
Of course, $\wh\s$ and $\widehat\s_{r}$ have the same kernel, yielding a
generalization of a result of Chatterjee [Cha].

\Prop{12.14}  {\sl  Let $R=d\G$ where $\G$ is integrally
flat. Then $R$ is linearly equivalent to zero if and only if
$\int_{\G}\theta \in\bbz $ for all $\theta\in
\har^{n-k}_0(X)$.}

\Ex{12.15. (The case of dimension 0)} For each point $p\in X$,
considered as a cycle of degree $n$, equation (12.3) takes the form
$d\s(p)=\omega-p$ where $\omega$ is the riemannian volume element
nomalized to have total volume one. For any other point $q$ we have
$\s(q-p)=\s(q)-\s(p)$ and the Abel-Jacobi map in this case restricts
to a map
$$
J:X \ \arr\ \Jac^{n-1}(X)
$$
given by
$$
J(q) = \widehat{\s}(q-p) = \widehat H(\G_{pq})
$$
where $\G_{pq}$ is any smooth curve in $X$ joining $p$ to $q$.
Note that $q$ is linearly equivalent to $p$  iff
$\int_{\G_{pq}}\theta\equiv 0$ mod $\bbz$ for all harmonic 1-forms
$\theta$ with integral periods. The map $J$ extends linearly to all
0-cycles on $X$ and in particular to symmetric powers of $X$.
This strictly generalizes the classical map defined in Riemann surface
theory.  For 3-manifolds it appeared in Hitchin's discussion of gerbes
and the mirror symmetry conjecture [Hi].

\Ex{12.16}  Let $X$ be a compact K\"ahler $n$-manifold and
denote by $Z_{m,0}$ the  holomorphic $m$-chains on $X$ which
are homologous to  zero.  (A {\sl holomorphic $m$-chain} is
a  rectifiable cycle of the form $z=\sum_i n_i[V_i]$ where
$n_i\in \bbz$  and the $V_i$ are irreducible complex
analytic  subvarieties of dimension $m$.) The restriction of
$\widehat\s$ to $Z_{m,0}\subset d\cR_{2m}$ can be seen to
coincide with the classical mapping into the Griffiths $m$th
intermediate Jacobian
$$
\har^{2(n-m)-1}/\har^{2(n-m)-1}_0\ \cong\
(H^{2m+1,0}\oplus H^{2m,1}\oplus  \dots\oplus   H^{m+1,m})^*
/H_{2m+1}(X;\,\bbz)
$$
To see this isomorphism note first that the map $\phi\mapsto h_{\phi}$
defined by  $h_{\phi}(\psi)=\int_X\phi\wedge\psi$ gives an isomorphism
$\har^{2(n-m)-1}/\har^{2(n-m)-1}_0\ \cong\
(\har^{2m+1})^*/H_{2m+1}(X;\,\bbz)$ and then use the canoncial
isomorphism $\har^{2m+1} =
H^{2m+1,0}\oplus H^{2m,1}\oplus  \dots\oplus   H^{m+1,m}$.

     \vskip .3in

 

\subheading{\S 13. Characteristic sparks and degeneracy
sparks}
Sparks occur naturally in the theory of singular connections developed
by the authors in [HL$_1$].  They appear as canonical
coboundary, or ``transgression'' terms relating the  classical
Chern-Weil forms  with certain characteristic currents defined by the
singularites of a  bundle map. We present some basic examples here.

Suppose  that $E$ and $F$ are
vector bundles with connection over a manifold $X$, and let $G$ denote
the Grassmann compactification of Hom$(E,F)$. There is a   flow $\vf_t$
on $G$ engendered by the flow $\wt{\vf}_t(e,f) = (te, f)$ on $E\oplus
F$. Constructing operators like those in \S 11 from the flow $\vf_t$ and
applying them to characteristic forms of the tautological bundle
$U\to G$ leads to a generalized Chern-Weil theory relating
characteristic forms and singularities of bundle maps
$\a:E\to F$ (See [HL$_{1-4}$]).  In general, for each characteristic
form $\Phi(\Omega)$ for $E$ or $F$ one obtains a formula of the type
$$
d T \ =\ \Phi(\Omega) - \sum_k \text{Res}_{k,\Phi}\, [\Sigma_k(\a)]
\tag13.1
$$
where $T\in \ce^k_{\Lloc}(X)$, $[\Sigma_k(\a)]$ is a current defined
by the dropping by $k$ of the rank of $\a$, and  $\text{Res}_{k,\Phi}$
is a smooth residue form defined along $\Sigma_k(\a)$.
 
In many important cases the residues are constants.  Whenever this is
true, $T$ defines a spark, and so in particular a differential
character.  Here are some examples.

\Ex{13.2 (Euler sparks)}  Let $E$ be an oriented real vector bundle
of even rank with an orthogonal connection, and let $s$ be a
cross-section with non-degenerate zeros, or more generally an {\sl
atomic}  cross-section.  Then as pointed out in Remark 9.10 there is a
canonically defined  {\bf Euler spark} $\s(s)\in \ce^{n-1}_{\Lloc}(X)$
satisfying the  spark equation 
$$
d\s(s) \ =\ \chi(\Omega) - \Div(s)
$$
where $\chi(\Omega)$ is the Euler form, which depends on the connection
but not on the section $s$, and $\Div(s)$  depends on $s$ but is
independent of the connection. The comparison formula for Euler sparks
(see [HZ$_2$, Thm. 7.6]) says that the Euler differential character
$\widehat\chi(E,D) = \widehat{\s}(s)$ determined by  $\s(s)$  is
independent of the section $s$.

\Ex{13.3 (Chern sparks)} Let $E$ be a complex vector bundle of rank
$n$ with unitary connection. Let $s=(s_0,\dots,s_{n-k})$ be a set of
$n-k+1$ sections with  a well defined { linear dependency current}
$\DD(s)$. This is an integrally flat cycle essentially defined as 
the locus where $s_0,\dots,s_{n-k}$ are linearly dependent
(See [HL$_1$, p.260]). Then there is a canonically defined
{\bf k$^{\text{th}}$ Chern spark}  $C_k\in \ce^{2k-1}_{\Lloc}(X)$ with
$$
dC_k\ =\ c_{k}(\Omega) -\DD(s)
$$
where $c_{k}(\Omega)$ is the k$^{\text{th}}$ Chern form of the
connection. Here $[\DD(s)]=c_k(E)\in H^{2k}(X;\bbz)$.

\Ex{13.4 (Pontrjagin-Ronga sparks)} Let $E$ be a real vector bundle of
rank $n$ with orthogonal connection. Let $s=(s_1,\dots,s_{n-2k+2})$ 
be a set of  sections of $E$ with  a well defined 2-dependency current 
$\DD_2(s)$.  This is an integrally flat cycle essentially defined as 
the locus where  dim(span$\{s_1,\dots,s_{n-2k+2}\})\leq n-2k$  (See
[HL$_{2}$]). Then there is a canonically defined {\bf k$^{\text{th}}$
Pontrjagin-Ronga spark}  $P_k\in \ce^{4k-1}_{\Lloc}(X)$ with
$$
dP_k\ =\ p_{k}(\Omega) -\DD_2(s)
$$
where $p_{k}(\Omega)$ is the k$^{\text{th}}$ Pontrjagin form of the
connection. It was shown by Ronga [R] that
$$
[\DD_2(s)] \,\,=\,\, p_k(E) \,\, + \,\, W_{2k-1}(E)\, W_{2k+1}(E)
\in H^{4k} (X,\bbz)
$$
where $p_k$ denotes the $k$th integral Pontrjagin class 
and $W_{j}$ is the $j$th integer Stiefel-Whitney
class defined by  $W_{j} := \delta_*(w_{j-1})$ where
$w_{j-1}$ is the  mod 2 Stiefel-Whitney class and 
$\delta_*$ is the Bockstein coboundary from the
exact sequence $0\to\bbz  \overset 2 \to \to \bbz\to
\bbz_2\to 0$.  The characters determined by $C_k$ and $P_k$ are 
independent of the choice of sections (See [HZ$_2$]).

\Ex{13.5 (Thom-Porteous sparks)}  Let $\a:E\to F$ be a bundle map
between complex vector bundles with unitary connections where
rank$(E)=m$ and rank$(F)=n$.  Suppose $\a$ is $k$-atomic (See [HL$_1$,
p.260]) so the  k$^{\text{th}}$ degeneracy current $\Sigma_k(\a)$ is
well defined. This is an integrally flat cycle which measures the
dropping of $\a$ in rank by $k$.
Then there is a canonical {\bf k$^{\text{th}}$ Thom-Porteous
spark}  $TP_k\in \ce^{2(m-k)(n-k)-1}_{\Lloc}(X)$ with
$$
d TP_k\ =\ \Delta^{(m-k)}_{n-k}\left\{c(\Omega^F)
c(\Omega^E)^{-1}\right\} - \Sigma_k(\a)
$$
where $\Delta^{(m-k)}_{n-k}$ denotes the Shur polynomial in
the total Chern form of $E-F$ (See [HL$_2$]).
\medskip

There is an analogue of 13.5 for real bundles which involves 
the Pontrjagin forms.

These degeneracy sparks appear in many geometric situations.  For
example the associated degeneracy currents can be measuring: 

 (i)  The higher complex tangencies of an immersion $f:X\to Z$ of a real
manifold into

\ \ \ \ \  a complex manifold,

(ii) The higher order tangencies of a pair of foliations on $X$,

(iii)  The higher order singularities of mappings $f:X\to X'$ between
smooth manifolds,

(iv)  The higher degeneracies of $k$-frame fields,

\noindent
etc.  See [HL$_2$] for a full discussion.      \vskip .3in


\subheading{\S 14. Characters in degree zero} 
In low dimensions differential characters have particularly nice
interpretations, and so do the corresponding spaces of sparks.
These cases nicely illustrate the relationship between sparks and
characters, and also the duality we have established here. 

We begin with degree zero. Recall that $\hH^0(X) =
\cs^0(X)/\cR^0(X)$  where $\cs^0(X)$ denotes the space of
$\Lloc$-sparks  and $\cR^0(X)$ is the 
group of rectifiable currents of degree 0 on $X$. Note that
$\cR^0(X)=\Lloc(X, \bbz)$, the space of integer-valued functions in
$\Lloc(X)$.

For a manifold $X$ let 
$G(X) = C^\infty(X, S^1)$ denote the smooth maps to the
circle $S^1=\bbr/\bbz$ (the gauge group for complex line
bundles) and define  
$$
\widetilde G(X) = \{f\in \Lloc(X) \,:\, \pi\circ f \text{\ is
smooth}\}
$$
where $\pi:\bbr\to\bbr/\bbz$ is the projection.

\Prop{14.1} {\sl The spaces $\wt G(X)$ and $\cs^0(X)$ agree as
subspaces of $\Lloc(X)$.  
In particular, there is an isomorphism of short exact sequences:
$$
\CD
0 @>>>\Lloc(X;\,\bbz) @>>>\widetilde G(X) @>{\pi_*}>> G(X) @>>>0 \\
@.  @V{=}VV      @V{=}V{}V   @V{\cong}V{}V  \\
0 @>>>\Lloc(X;\,\bbz) @>>>\cs^0(X) @>{q}>>  \hH^0(X) @>>>0
\endCD
$$
where $\pi_*(f) \equiv \pi\circ f$, and $q$ is the canonical
quotient.}

\pf
For each point $p\in S^1$ choose $r\in\bbr$  with $\pi(r)=p$, and let
$\theta_p:S^1\to[r,r+1)$ be the unique function with $\pi\circ
\theta_p=$Id. This $\Lloc$-function is smooth outside $p$ and
satisfies the spark equation 
$$
d\theta_p\ =\ \Theta - [p]
\tag14.2
$$ 
where $\Theta$ is the unit volume (arc-length) form on $S^1$.
Given $g\in G(X)$, pick $p\in S^1$ a regular value of $g$.
Define $\overline f = g^*\theta_p$, $\overline \phi = g^*\Theta$,
and $\overline R = [g^{-1}(p)]$.  Then the equation
$$
d\overline f \ =\ \overline {\phi} - \overline R
\tag14.3
$$
is the pull-back of (14.2) by $g$.  Note that $\pi\circ \overline f =
g$. In particular,   $\pi_*$ is surjective.  Now given any $f\in \wt
G(X)$ such that $\pi\circ f=g$, we see that $f=\overline f +f_0$ with
$f_0 \in \Lloc(X, \bbz)$.  Hence, by (14.3), $df =        \overline
{\phi}-\overline R+df_0$ and since $\overline R+df_0$ is integrally
flat, we conclude that $\wt G(X)\subset \cs^0(X)$.

Now for any given $f\in \cs^0(X)$,  $d_1\langle f\rangle \equiv\phi$
belongs to $\cz^1_0(X)$ and hence $\phi$ determines a smooth map
$g:X\to S^1$, unique up to rotation, with $\phi = \overline{\phi} =
g^*\Theta$.  Therefore $d(f-\overline f)= R-\overline R$.  This
implies that $f-\overline f +c\in \Lloc(X,\bbz)$ for some constant $c$
since $R-\overline R$ is zero in integral cohomology and the kernel of
$d$ on ${\cd'}^0$ is the constants. Therefore, $\cs^0(X)\subset \wt
G(X)$. \qed

\Remark {14.4} For every pair 
$(g,p)\in G(x)\times S^1$ where $p$ is a regular value of $g$, the
construction above gives a $\bbz$-family of particularly nice sparks
$\overline f$ with $\pi_*(\overline f)=g$.

\medskip

\Remark {14.5} Note that in this degree regularity implies that every
generalized spark  $a\in {\cd'}^0(X)$ is actually in $\Lloc$ since
$da=\phi-R$ is a measure.

\medskip

 We now examine duality in this degree. Note that there
is a basic family of continuous homomorphisms 
$$
h_x : C^{\infty}(X,\, S^1) \ \arr \ S^1 \qquad\qquad\text{for   \ }
x\in X
\tag14.6
$$
given by $h_x(g) = g(x)$. These homomorphisms do {\bf not} lie
in the smooth dual of $\hH^0(X)$.  However, as noted in 6.7 {\sl any
Federer smoothing homotopy of the current $[x]$ does give a smooth
homomorphism.} By 2.8, or by standard constructions, there are
families of $\Lloc$ $(n-1)$-forms  $a_\e$ with supp$(a_\e)\subset
B_\e(x)$, an $\e$-ball about $x$ in local coordinates, such that  
$$
d a_\e \ =\ \Omega_\e - [x]
\and
\lim_{\e\to 0} a_{\e} = 0 
$$
where the limit is taken in the {\sl flat} topology.  
In particular $\Omega_\e \to [x]$ as $\e\to0$.
Now   for each $f\in C^{\infty}(X,\, S^1)$ we have a well-defined 
integral
$$
h_{\e,x}(f) \ \equiv\ \int_X \overline f \Omega_\e
 \ \ \ \ \left(\text{mod}\  \bbz \right) 
$$ 
where $\overline f :X\to \bbr$ is any lift of $f$ which is
continuous on $B_\e(x)\supset \supp(\Omega_\e)$.  By (6.9) we see
that  $h_{\e,x}(f) \equiv (a_\e * f)$.  Now each of the homomorphisms  
$$
h_{\e,x} :C^{\infty}(X,\, S^1)\ \arr \  S^1
$$ 
is smooth, and we have that $h_{\e,x}\to h_x$ as $\e\to 0$.
Thus the homomorphisms (14.6) lie at the boundary of the smooth ones. 

One could consider any spark $a\in \cs^{n-1}(X)$ with 
$
da = \Omega - [x],
$
to be a generalized smoothing homotopy of the point mass $[x]$. 
Note that  $\cs^{n-1}(X)$ is exactly the group generated by such
things.  The smooth dual of $C^{\infty}(X,\, S^1)$ is exactly the group
generated by these smoothing  homotopies. Notice that for 
$f\in \wt G(X) = \cs^0(X)$ with $df=\phi-R$ as above, we have by 3.2
that 
$$\aligned
\langle a \rangle * \langle f \rangle \, [X]
\ &\equiv\ \int_X a\wedge\phi  +(-1)^n f(x) \ \ \ \ \left(\text{mod}\ 
\bbz \right)  \\ 
&\equiv \ \int_R a +(-1)^n \int_X f \Omega \ \ \ \ \left(\text{mod}\ 
\bbz \right)  
\endaligned
$$
where we assume $a$ to be smooth on the support of $R$ (cf. 4.7).
This shows explicitly how the smooth homomorphism corresponding to a
character $\a\in \hH^{n-1}(X)$ differs from point evaluation on 
$\d_2(\a)$ or integration against the form $\d_1(\a)$. 

\medskip

When $H^1(X;\,\bbz)=0$, we have $C^{\infty}(X,\, S^1) =
C^{\infty}(X)/\bbz$ whose smooth dual is just  
$
\cz^n_0(X)
$
where the pairing is given by integration.  In general there is a
short exact sequence
$$
0\arr C^{\infty}(X)/\bbz \arr C^{\infty}(X,\, S^1) \arr 
H^1(X;\,\bbz)\arr0
$$
whose dual sequence 
$$
0 \longleftarrow  \cz^n_0(X)  \longleftarrow 
\Hom_{\infty}(C^{\infty}(X,\, S^1),\, S^1) \longleftarrow 
\Hom(H_{n-1}(X,\, \bbz),\, S^1) \longleftarrow  0 
$$
corresponds to 1.12.      \vskip .3in


\subheading{\S 15. Characters in degree 1} 
We now consider degree one.
Let $E\to X$ be a smooth complex line bundle with unitary connection
$D$. Then to each atomic cross-section $\s$
(for example any section which vanishes non-degenerately), there is
a canonically associated spark $T(\s)$ with the property that 
$$
d T(\s) \ =\ c_1(\Omega) - \text{Div}(\s)
\tag15.1$$
where $\Omega$ is the curvature form of $D$ ([HL$_1$, II.5]). It is
given locally by the formula 
$$
T(\s)\ =\  \text{ Re\  } \left\{{\tsize {\frac 1 {2\pi i}}}
\left( \omega -\frac {du}{u}\right)\right\}
\tag15.2
$$
where $\s=u e$, $De=-\omega\otimes e$, and
$e$ is a local section of $E$ with $\|e\|\equiv 1$.
The smooth section $\s$ is defined to be {\bf atomic} if $ du/u\in
\Lloc$, and by definition
$$
\Div(\s) \ = \ d  \text{ Re\  } \left\{{\tsize {\frac 1 {2\pi i}}}
 \frac {du}{u}\right\}
\tag15.3
$$

Denote by $L(X)$ the set of (isomorphism classes of) complex line
bundles with unitary connection up to gauge equivalence on $X$,
and set
$$
\widetilde L(X)\ =\ \{(E,D,\s) :  \text{$E$, $D$ are as above and $\s$
is an atomic cross-section of } E\} 
$$ 
\medskip\noindent
{\bf Proposition 15.4.} {\sl There is a commutative diagram
$$\CD
\widetilde L(X) @>{p}>> L(X) \\
@V{T}VV   @V{\cong}V{T_0}V  \\
\cs^1(X) @>{q}>>  \hH^1(X)
\endCD
$$
where $p(E,D,\s) \equiv (E,[D])$, $q$ is the canonical quotient, and
$T_0$ is an isomorphism.}

\pf  The comparison formula for Euler sparks proven in 
[HZ$_2$, Thm 7.6] says that the map $q\circ T$ is independent
of the choice of cross-section $\s$ and therefore descends to the 
mapping $T_0$. To see that $T_0$ is surjective, fix $\a \in \hH^1(X)$ 
with $\d\a = \Phi-e$. The class $e\in H^2(X;\,\bbz)$ corresponds
to a complex line bundle $E$ on $X$  and the de Rham class of the
smooth form $\Phi$ is the real Chern class $c_1(E)$.  Hence there
exists a unitary connection $D$ on $E$ whose associated Chern form
$c_1(\Omega)=\Phi$. To see this note that for any unitary connection
$D'$ one has  $c_1(\Omega')-\Phi = d \phi$ for a smooth real-valued
1-form $\phi$, and then set $D= D'-\phi$.   For any atomic section
$\s$ of $E$ we have  $\d\{\a-q\circ T(\s,E,D)\}=0$.  Thus by 1.14,
$\a-q\circ T(\s,E,D)=[\psi] \in
H^1(X;\,\bbr)/H^1_{\text{free}}(X;\,\bbz)$ for a smooth closed 1-form
$\psi$.  By (15.2),  $T(\s,E, D+\psi) =T(\s,E,D) +\psi$, and so
$q\circ T(\s,E, D+\psi) 
 =\a$ as desired.  

 We now prove injectivity.  Suppose
$T_0(E,D)=0$.  Then since $\d T_0(E,D)=0$, $E$ is topologically
trivial and $D$ is flat. Choose a section $e$ of $E$ with
$\|e\|\equiv 1$, i.e., a trivialization of $E$. Then $T(e,E, D)$ is a
smooth, $d$-closed 1-form whose class in 
$H^1(X;\,\bbr)/H^1_{\text{free}}(X;\,\bbz)$ represents the holonomy of
$D$.  Since this is 0, the connection is gauge equivalent to the
trivial one. \qed 
 
\Remark{15.5} Given any $(E,D)\in L(X)$ one can choose a section
$\s$ of $E$ with regular zeros, so that $\Div(\s)$ is a smooth
submanifold of codimension-2 in $X$. These give particularly nice
sparks in $\wt L(X)$ representing the character associated to
$(E,D)$.  They are analogous to the nice sparks constructed in
the last section (cf. (13.4)).

\Remark{15.6}  Essentially all sparks in $\cs^1(X)$ can be
represented by triples $(E,D,\s)$ where the $\s$'s are certain 
measurable sections of the unit circle bundle of $E$. 
In particular all sparks $a$ such that $d_2a$ has support of measure
zero can be so constructed.

\medskip

We now consider the dual to $\hH^1(X)\cong L(X)$.  As in degree 0
there is a natural family of continuous homomorphisms which are not
smooth; namely for each piecewise-$C^1$-loop $\g$ in $X$ we have the
{\bf holonomy mapping}
$$
h_\g : L(X)\ \arr\ S^1.
$$
The characters $\hH^{n-2}(X)$ represent generalized smoothing
homotopies of these holonomy homomorphisms.  They are represented by  
$\Lloc$-forms $a$ of degree $n-2$ with $$da = \Omega - \g$$
where $\g$ is such a loop and $\Omega$ is a smooth $(n-1)$-form.

 Notice that for 
$e\equiv (E,D,\s) \in \wt L(X)$ whose associated  spark $T(\s) \in
\cs^1(X)$ satisfies (15.1), we have by 3.2 that 
$$\aligned
\langle a \rangle * \langle e \rangle \, [X]
\ &\equiv\ \int_X a\wedge c_1(\Omega)  + h_{\g}(e) \ \ \ \
\left(\text{mod}\  \bbz \right)  \\ 
&\equiv \ \int_{\Div (\s)} a +  \int_X \Omega\wedge T(\s) \ \ \ \
\left(\text{mod}\  \bbz \right)  
\endaligned
$$
If we take a local Federer family of smoothing homotopies
$a_\e$ with $da_\e = \g_\e - \g$ where $a_\e$ is supported
in an $\e$-neighborhood of $\g$ and $\lim_{\e\to 0} a_\e =0$ in the
flat topology, the first equation  above becomes
$$
\langle a_\e \rangle * \langle e \rangle \, [X]
\ \equiv\ \int_X a_\e\wedge c_1(\Omega)  + h_{\g}(e) \ \ \ \
\left(\text{mod}\  \bbz \right)   
$$
and we see that as $\e\to 0$ the smooth homomorphisms given by $a_\e$
converge to  $h_\g$. Thus the holonomy
homomorphisms lie at the boundary of the smooth ones.      \vskip .3in


\subheading{\S 16. Characters in degree 2 (Gerbes)} 
Differential characters in degree 2 are intimately connected 
to {\sl Gerbes with connection}.  These objects 
generalize the previous two cases, and have been much discussed in
the recent literature.  We shall
briefly present them and their relationship to characters. For a
rounded discussion with applications see the excellent article of
Hitchin [Hi] and the book of Brylinski [Br]. 

Just as a manifold is presented by an atlas of coordinate
charts, a gerbe  $X$ is presented by \v Cech 2-cocycle 
with values in $S^1$.  Thus the data of a Gerbe is an  open
covering $\{U_{\a}\}_{\a\in A}$ of a manifold $X$ (which for
simplicity we assume to be acyclic)  and continuous functions 
$$
g_{\a\b\g}:U_{\a}\cap U_{\b}\cap U_{\g}\ \arr\ S^1
$$
with $g_{\a\b\g}=g_{\b\a\g}^{-1}=g_{\g\b\a}^{-1}=g_{\a\g\b}^{-1}$
and satisfying the cocycle condition:
$$
g_{\b\g\d} \cdot g_{\a\g\d}^{-1}\cdot 
g_{\a\b\d} \cdot g_{\a\b\g}^{-1}\ \equiv\ 1
\qquad\text{ on }\ \ \  U_{\a}\cap U_{\b}\cap U_{\g}\cap U_{\d}.
$$
Gerbes form an abelian group under ``tensor product'', and they
can be pulled back by a continuous map $f:Y\to X$.  Gerbes arise
naturally in many contexts including Spin$^c$-geometry, the study of
Calabi-Yau manifolds, geometric quantization, and the topological
Brauer group (see [Hi], [Br]). They are most interesting when
equipped with a connection.  A {\bf connection} on the gerbe 
$g_{\a\b\g}$ defined on a manifold $X$ consists of the
following:\medskip

\hskip .5in (a)\ \ \ \  $\phi \in \ce^3(X)$\qquad 
\qquad\qquad (called the {\sl
curvature} of the gerbe),

\hskip .5in (b)\ \ \ \ $A_\a \in  \ce^2(U_\a)$\qquad\  \ \ \ \qquad for
each $\a$

\hskip .5in (c)\ \ \ \ $A_{\a\b} \in \ce^1(U_\a\cap U_\b)$ \qquad\, for
each $\a,\b$.

\medskip\noindent
with the property that:\medskip

\hskip .5in (i)\ \ \ \   $\phi \ =\ d A_\a$
\hskip 2.in  \ on $U_\a$

\hskip .5in (ii)\ \ \  $A_\b-A_\a \ =\ d A_{\a\b}$\hskip 1.6in on
$U_\a\cap U_{\b}$

\hskip .5in (iii)\ \  $A_{\a\b}+A_{\b\g}+A_{\g\a}
 \ =\ -id \log g_{\a\b\g}$\qquad\qquad\ \ \,
on $U_{\a}\cap U_{\b}\cap U_{\g}$

\medskip\noindent
{\bf Note 16.1.} This is ``one step up'' from a connection on a line
bundle. A line bundle $L$ is presented by a \v Cech 1-cocycle
$\{g_{\a\b}\}$ with values in $S^1$. A connection on $L$ consists of: 
\medskip

\hskip .7in (a$'$)\ \ \ \  $\phi \in \ce^2(X)$\qquad 
\qquad
  (called the {\sl curvature} of the line bundle),

\hskip .7in (b$'$)\ \ \ \ $A_\a \in  \ce^1(U_\a)$\qquad\  \ \ 
for each $\a$

\medskip\noindent
with the property that:\medskip

\hskip .7in (i$'$)\ \ \ \  $\phi=dA_\a$ \hskip 1.5in\ \, \  on $U_\a$,
and

\hskip .7in (ii$'$)\ \ \  $A_\b-A_\a = -id \log g_{\a\b}$
\hskip .8in  on $U_\a\cap
U_\b$.

\medskip\noindent
{\bf A slight modifiction.}  To make  stronger analogy with our de
Rham-Federer  theory, we  modify the above notion of a gerbe with
connection to one which is essentially equivalent.  We now consider
a gerbe with connection to be a triple
$$
A=(\{A_{\a}\},\{A_{\a\b}\},\{A_{\a\b\g}\})
\tag16.2$$
where $\{A_{\a}\},\{A_{\a\b}\}$ are as above and 
$A_{\a\b\g} \in \ce^0(U_\a\cap U_\b\cap U_{\g})$
is some choice
$$
A_{\a\b\g}=-i \log  g_{\a\b\g}
$$

\medskip\noindent
{\bf Gauge Equivalence.} Two gerbes with connection  
$A$ and $A'$ as in (16.2) are said to be {\bf gauge
equivalent} (written $A\sim A'$) if there exist \medskip

\hskip 1in (a)\ \ \ \ $B_\a \in  \ce^1(U_\a)$\qquad\  \ \ \ \qquad
for each $\a$

\hskip 1in (b)\ \ \ \ $B_{\a\b} \in \ce^0(U_\a\cap U_\b)$ \qquad\,
for each $\a,\b$.

\medskip\noindent
with the property that:\medskip

\hskip .2in (i)\ \ \ \   $dB_{\a} \ =\  A_\a-A_{\a}'$
\hskip 2.in \ \ \ \ \  on $U_\a$

\hskip .2in (ii)\ \ \  $B_\b-B_\a -d B_{\a\b}
\ =\ A_{\a\b}-A_{\a\b}'$\hskip 1.1in \ \ \ \ \ \  \, on
$U_\a\cap U_{\b}$

\hskip .2in (iii)\ \ 
$B_{\a\b}+B_{\b\g}+B_{\g\a}+ 2\pi n_{\a\b\g}
 \ =\ A_{\a\b\g}-A_{\a\b\g}' $  \  \qquad\ \ \ \ \  
on $U_{\a}\cap U_{\b}\cap U_{\g}$
\medskip\noindent
where $ n_{\a\b\g}$ is a locally constant integer-valued
function. This directly generalizes the notion of gauge equivalence
for line bundles with connection. 

\define\Dtot{\DD} 
\medskip\noindent
{\bf The \v Cech-deRhamForm bicomplex.} The ideas above can be nicely
packaged by considering the double complex
$
(\cc^*(\cu,\ce^*),\, \Dtot)
$
of the open cover $\cu =\{U_{\a}\}$ where $\Dtot=d+(-1)^k\d$ and $\d$
denotes the \v Cech differential. A gerbe with connection is then
simply an element 
$$
A\in \bigoplus_{j+k=2} \cc^{j}(\cu,\ce^k)
$$
which satisfies the equation
$$
\Dtot A = \phi - R
\tag16.3$$
where 
$$
\phi \in \cz^0(\cu, \cz^3)\subset\cc^0(\cu, \ce^3)  \and  
R \in \cz^3(\cu, 2\pi\bbz)\subset\cc^3(\cu, \ce^0)
$$
Thus $\phi$
is a globally defined closed 3-form and $R=\delta \{A_{\a\b\g}\}$ is
an integral \v Cech 3-cocycle.

Two such creatures $A$ and $A'$ are gauge equivalent if there exist
$$
B \in \bigoplus_{j+k=1} \cc^{j}(\cu,\ce^k)
\qquad\text{and}\qquad S\in \cc^2(\cu, 2\pi\bbz)
$$ 
such that
$$
A-A'\ =\ \Dtot B +S
\tag16.4
$$

\noindent{\bf Note 16.5.} The reader will certainly have noticed the
analogy of equation (16.3)  with the spark equation (0.1), and the
analogy of (16.4) with the equivalence relation on sparks introduced
in 1.6 and 2.4.  Indeed, as one might guess at this point, gauge
equivalence classes of gerbes with connection coincide with
differential characters of degree-2.  The key to this observation is
the following lemma whose proof is left to the reader.

\medskip\noindent
{\bf Lemma 16.6.} \ {\sl Let $A$ be a gerbe with connection
with $\Dtot A = \phi - R$ as above.

(i)\ \ \ If $R =0$, then $A\sim (F,0,0)$ where $F\in \ce^2(X)$ is a
global smooth 2-form.
 
(ii)\ \ If $\phi =0$, then $A\sim (0,0,T)$ where $T\in
\cc^2(\cu,\bbr)$ is a cochain with constant coefficients.
}

\medskip\noindent
{\bf The holonomy map.} Just as connections on line bundles have
holonomy on each oriented closed curve in $X$, gerbes with connection
have holomony on each oriented closed  surface in $X$. It is defined
as follows. Let $\Sigma\subset X$ be such a surface and $A$ a gerbe
with connection as above.  The restriction of $A$ to $\Sigma$
satisfies the hypothesis of Lemma 16.6(ii) since
$\phi\bigl|_{\Sigma}=0$. Hence $A\bigl|_{\Sigma} \sim (0,0,T)$
where $\d T\equiv 0$ ( mod  $2\pi\bbz$). Now
$T$ determines \v Cech 2-cocycle on $\Sigma$ with values in the
constant sheaf $2\pi\bbr/2\pi\bbz = S^1$, which is unique up to 
\v Cech coboundaries.
Evaluating on the fundamental cycle $[\Sigma]$ gives the {\bf
holonomy}  
$$
 h_A(\Sigma)
\in H^2(\Sigma;S^1)=S^1. 
$$
One verifies that the holonomy $h_A$ {\sl depends only on the gauge
equivalence class of } $A$.

\Prop{16.7} {\sl There is a natural equivalence of functors:}
$$
\hH^2(X)\ \cong\ \{\text{gauge equivalence classes of gerbes on } X\}
$$
In the Cheeger-Simons picture this equivalence associates
to a gerbe with connection $A$ the homomorphism from integral 2-cycles
to $S^1$ given by the holonomy $h_A$. 

\Remark{16.8} The {\sl  \v Cech-deRhamForm)}\  presentation of
differential characters sketched here for degree-2 extends easily to
all degrees.  A theory of this type was first mentioned in [FW, \S 6].
Full details will appear in a sequel [HL$_7$] to this paper.

\medskip\noindent
{\bf Obvious question.} {\sl Is there a direct relationship between
the \v Cech-deRhamForm theory and the de Rham-Federer theory of
differential characters?}\ The answer is yes -- there exists a large 
{\sl \v Cech-deRhamCurrent)}\ double complex which presents
differential characters and from which one can distill these two
theories as special cases.  Details of this will appear in [HL$_7$].      \vskip .3in


\subheading{\S 17. Some historical remarks on the complex analogue}
Jeff Cheeger realized in 1972 that characters could be represented by
forms with singularities, and explicit reference is made to this in
[CS]. However, the first appearance of the spark
equation (0.1) occured in the work of Gillet-Soul\'e on intersection
theory for arithmetic varieties [GS$_2$]. They began with
a smooth complex projective variety $X$ defined over $\bbz$, and
introduced currents $G$ of bidegree $(k,k)$, called {\sl Green
currents}, on  $X$ with the property that   
$$
 d d^C G\ =\ \phi - C
\tag17.1
$$
where $\phi$ is a smooth $2k$-form and $C$ is an algebraic cycle 
of codimension-k on $X$.  They then divided the group of  Green
currents by  the image of $\partial$ and  $\overline\partial$ and a
certain subgroup which yields linear equivalence on the cycles $C$.
The resulting quotient was called the $k$th  {\sl arithmetic
Chow group of $X$}.  These groups formed a contravariant functor
with a non-trivial ring structure which allowed Gillet and Soul\'e 
to extend previous results of Arakelov, Deligne, Beilinson and
others to a full arithmetic intersection theory.  Gillet and Soul\'e
went on to show that algebraic bundles with hermitian metrics
have  refined characteristic classes in these arithmetic Chow groups
[GS$_3$]. Eventually  much deep work of Bismut, Gillet, Soul\'e and
Lebeau  led to a  refined ``arithmetic''  version of the
Riemann-Roch-Grothendieck Theorem [BGS$_*$], [BL], [GS$_4$]  and to
arithmetic analogues of the standard conjectures [GS$_5$].

The close similarity of arithmetic Chow groups with 
differential characters and  the parallel nature of the exact
sequences in the two theories led 
Gillet-Soul\'e and B. Harris   to explore the
relationship. They realized that the natural analogue of equation
(17.1) in real geometry is the spark equation, and they took this
approach  to differential characters.  They  all realized that for 
algebraic cycles of degree 0, differential characters should be 
related to the Abel-Jacobi mapping into intermediate Jacobians.
In [GS$_1$] the relationship between the two theories is
sketched and,  among other things,  a refined version of
Riemann-Roch-Grothendieck is proved for K\"ahler fibrations. In [H] it
is shown that differential characters gave a simplified and more
conceptual approach to certain  formulas in algebraic
geometry  originally proved by Harris via Chen's iterated integrals.

The authors came independently  to the spark equation
while developing the theory of singular connections and characterstic
currents ([HL$_{1-4}$] and [HZ]) where they arise quite naturally.      \vskip .3in



 

\redefine\D{\Delta}
\redefine\G{\Gamma}
\def\wt{\widetilde}

\centerline{\bf Appendix A.   Slicing currents by sections of a
bundle}

\bigskip

Let $E\to X$ be a smooth real vector bundle of rank $n$ and let $R$ 
be $p$-dimensional current on $X$ which is either flat or
rectifiable.  Let $v_1,...,v_k \in \G(E)$ be smooth sections such
that 
$$
\text{span}\{v_1(x),...,v_k(x)\} \ = \ E_x
\qquad \text{for  all $x\in X$.}
\tag{A.1}$$
Define a map
$$
\psi :X\times \bbr^k\ \to\ E
$$
by $\psi(x,\xi) = \sum_i \xi_i v_i(x)$.  This gives a short exact
sequence of vector bundles over $X$:
$$
0\to F \to X\times \bbr^k\ @>{\psi}>>\ E @>>> 0.
$$
Consider the current $\widetilde R \equiv p^*R$ where $p : F\to X$ is
the bundle projection.  Then $\widetilde R$ is a flat or
rectifiable current (depending on $R$) of dimension $p+k-n$ in 
$X\times \bbr^k$. Applying Federer slicing theory [F, 4.3] to  the
projection  $\pi:X\times \bbr^k\to \bbr^k$ gives the following. 

\Theorem{A.2}  {\sl   The slice   
$$
R_{\x}\ = \ \langle \widetilde R, \pi, \xi\rangle
$$
exists in $\cf_{p-n}(X)$ for almost all $\xi \in \bbr^k$ and is 
rectifiable if $R$ is rectifiable.  Furthermore, for almost all 
$\xi \in \bbr^k$  we have 
$$
dR_{\x}\ = \ \langle d \widetilde R, \pi, \xi\rangle.
\tag{A.3}
$$
In particular, by (1.1) if $R$ is integrally flat so is $R_{\x}$ 
for almost all $\xi$. }

\medskip
Note that
$$
\text{supp}(R_{\x})\ \subseteq\ \left \{ (x,\x) \in 
\text{supp}(R)\times \bbr^k \,:\, \sum\x_iv_i(x) = 0\right \}.
$$
Note also that the maps $\psi$ and $t\psi$ have the same kernel for
all $t\in \bbr-\{0\}$. Thus we may reduce our parameter space from
$\bbr^k$ to $\BP_{\bbr}^{k-1}$.

\vskip .3in

\centerline{\bf Appendix B.  Intersecting currents}

\bigskip

Suppose now that $R$ and $S$ are currents of dimensions $p$ and $q$
respectively on a compact  $n$-manifold $X$.  Consider $R\times S
\subset X\times X$.  Let $\Delta\subset X\times X$ be the diagonal, and
choose an identification  of a tubular
neighborhood $\Delta_{\epsilon}$ of $\Delta$ in $X\times X$ with its
normal bundle $N \equiv N(\Delta) \cong T(X)$. Let $p:N\to \D \cong
X$ be the bundle projection. Then  $\D \subset N$ is the divisor
$\Div(v_0)$ of the tautological cross-section of $p^*N \to N$.

We now apply Appendix A to the  current
$(R\times S) \cap \D_{\epsilon} \cong (R\times S) \cap N$ and the
bundle $p^*N \to N$.  We choose sections $v_1,...,v_k \in \G (N)$
which satisfy (A.1) and proceed as above using the family
$(v_0,v_1,...,v_k)$.  We find that for almost all $\x =
(\x_1,...,\x_k) \in \bbr^k$ sufficiently small 
the slice $(R\times S)_{\x}$ of $R\times S$ by  
$$
v_{\x} = v_0 + \sum_{i=1}^k \x_i v_i
$$
exists in $\cf_{p+q-n}(N)$ and is rectifiable if $R$ and $S$ are
rectifiable. Furthermore, for almost all $\x$ we have
$$
d\{ (R\times S)_{\x}\}\ =\ \{d(R\times S)\}_{\x}
\ =\ \{(dR)\times S +(-1)^{n-p}R\times dS\}_{\x}.
\tag{B.1}$$

Now since $v_0$ is a transversal cross section of $p^*N$ we have
that for all $\|\x\|$ sufficiently small,  $\Div(v_{\x})$
is a small perturbation of the diagonal, and thus
$$
\Div(v_{\x})\ =\ \text{graph}(f_{\x})
$$
where $f_{\x} : X \to X$ is a diffeomorphism close to the identity.

\Def{B.2} We define the {\bf intersection}
$$
R\wedge (f_{\x})_* S
$$
of the currents $R$ and $(f_{\x})_* S$ in $X$ to be the $\x$-slice of
the current $R\times S$ in $X\times X$. This exists for almost all
$\x\in\bbr^k$ with $\|\x\|$ sufficiently small. Furthermore for
almost all such $\x$ we have by (B.1) that
$$
d\{R\wedge(f_{\x})_* S\}\ =\ (dR)\wedge(f_{\x})_* S
+(-1)^{n-p}R\wedge(f_{\x})_*(d S)
\tag{B.3}
$$

\vskip .3in

\centerline{\bf Appendix C. The de Rham-Federer approach to
integral cohomology}

\bigskip

On any manifold $X$ the complex of integrally flat currents
computes the integral cohomology $H^*(X;\,\bbz)$.  This follows
from the fact that the sheaves of such currents give an acyclic
resolution of the constant sheaf $\bbz$ (See [HZ$_1$]).
These comments apply also to the subcomplex of integral currents.

\Theorem{C.1 (The intersection product on de Rham-Federer
 $H^*(X;\,\bbz)$)} {\sl  The intersection of currents
introduced in B.2 gives a densely defined pairing on integral and on
integrally flat currents which descends to an associative,
graded-commutative multiplication on $H^*(X;\,\bbz)$.}

\pf  This is a straightforward consequence of 
Theorem A.2. \qed

\Note{C.2}  The intersection product on $H^*(X;\,\bbz)$ coincides
with the cup product.  Classical facts show that for certain
nice subcomplexes of the de Rham-Federer complex, where intersection
is well-defined, the two products agree on cohomology.   Hence they
must agree for general de Rham-Federer currents. However, proving the
coincidence of these products in any setting is not trivial.  The
basic reason is that the intersection of currents, when defined, is
local in nature, generalizing the wedge product on smooth currents,
whereas the cup product of cochains is not.
This important difference makes Cheeger's
definition [C] of the $*$-product on differential characters
non-trivial - it involves an infinite series like the one in \S 4. 
This difference also accounts for the difficulty in \S 4  of proving
the coincidence of his product with the one introduced here.

\Note{C.3}  Federer proves in [F, 4.4] that the complex of
integral (or integrally flat) currents with compact support
computes the {\bf homology} of $X$.  Thus Theorem A.1 also yields
an intersection ring structure on $H_*(X;\,\bbz)$.  When $X$
is compact and oriented, the implied isomorphism 
$H_{n-k}(X;\,\bbz) \cong H^k(X;\,\bbz)$ is Poincar\'e duality.

\Note{C.4} The contravariant functoriality of de Rham-Federer
cohomology is an immediate consequence of its sheaf-theoretic
definition. However, this pull-back can be defined more directly,
as we did for characters in 2.9.

\vskip .3in

\centerline{\bf Appendix D. The spark product}

\bigskip

\Theorem{D.1} {\sl For given sparks $\a\in \hH^k(X)$ and $\b\in
\hH^{\ell}(X)$ with $k+\ell \leq n$, there exist representatives $a\in
\cf^k(X)$ and $b\in\cf^{\ell}(X)$ for $\a$ and $\b$ respectively, with
$$
da = \phi - R \qquad\text{and}\qquad db = \psi - S
\tag{D.2}$$
with $\phi, \psi$ smooth and $R, S$ rectifible, such that the products
$a\wedge  b$, $a\wedge  S$, $R\wedge  b$
and $R\wedge  S$ are well-defined and  flat and $R\wedge S$ is
rectifiable. 
}

\pf
Let $a\in\a$ and $b\in\b$ be any flat representatives with exterior
derivatives given as in (D.2). 
Choose $\x$ as above so that 
$a\wedge f_{\x*}b$, $a\wedge f_{\x*}S$, $R\wedge f_{\x*}b$
and $R\wedge f_{\x*}S$ are well-defined, flat and the last is
rectifiable.  Note that 
$$
d(f_{\x*}b) = f_{\x*}\psi - f_{\x*}S.
$$ 
Now $\psi - f_{\x*}\psi = d\chi$ where $\chi$ is a smooth
$(\ell-1)$-form. Therefore, $\wt b = f_{\x*}b + \chi$ satisfies the
equation $d\wt b = \psi -f_{\x*}S$.  Then, since $S$ and 
$f_{\x*}S$ are homologous integral cycles, we have $\d\langle \wt b
\rangle =  \d\langle  b \rangle = \d\b$.  
Therefore, by the exact sequence 1.14,  $\b - \langle \wt b
\rangle = d\langle \eta \rangle$ where $\eta$  is a smooth
$d$-closed  $\ell$-form.  Hence,  
$$
\b \ =\ \langle  f_{\x*} b + \chi +\eta \rangle
$$
and the flat form $b' \equiv f_{\x*} b + \chi +\eta$ has the property
that the products $a\wedge b'$, $a\wedge d_2b'$, $d_2a\wedge b'$ and
$d_2a\wedge d_2b'$ are well defined flat currents and the last is
rectifiable. \qed

 

\centerline{\bf References}

\vskip .3in

\nobreak

\ref \key [BGS]\by J.-M. Bismut, H. Gillet and C. Soul\'e \paper 
Complex immersions and Arakelov geometry
\pages 249-331  \inbook The Grothendieck Festschrift {\bf 1},  Progress in
Mathematics  \vol 1 \publ Birkhauser \publaddr Boston \yr 1990  \endref

\ref \key [B]\by R. Bott \paper On a topological obstruction to 
integrability  
\pages 127-131 \jour Proc. of Symp. in Pure Math. \vol 16 \yr 1970 
\endref

\ref \key [BC$_1$] \by R. Bott and S.S. Chern \paper  Hermitian vector 
bundles and the equidistribution of the zeros of their holomorphic 
sections \jour Acta. Math. 
\yr 1968 \pages 71--112 \vol 114 \endref

\ref \key [BC$_2$] \bysame \paper Some formulas related to 
complex transgression \inbook Essays on Topology and Related 
Topics 
(M\'emoires d\'edi\'es \`a Georges de Rham \yr 1970 \publ 
Springer-Verlag
\publaddr New York \pages 48--57\endref

\ref \key [BGS$_1$] \by J.-M. Bismut, H. Gillet and C. Soul\'e \paper
Analytic torsion and holomorphic determinant bundles I, II, III\jour
Commun. Math. Phys. \vol 115 \yr 1988 \pages 49-78, 79-126, 301-351\endref

\ref \key [BGS$_2$] \bysame \paper
Bott-Chern currents and coomplex immersions\jour
Duke Math. J. \vol 60 \yr 1990 \pages 255-284\endref

\ref \key [BGS$_3$] \bysame \paper
Complex immersions and Arakelov geometry\inbook
Grothendieck Festschrift I\publ Birkh\"auser\publaddr Boston \yr 1990
\pages 249-331\endref

\ref\key [BL] \by J.-M. Bismut and G. Lebeau\paper Complex immersions and
Quillen metrics \jour Instituto Nazionale di Alta Mathematica,  Symposia
Mathematica \vol XI \yr 1973 \pages 441--445 \endref

\ref\key [Br]\by  J.-L. Brylinski\book Loop Spaces, Characteristic Classes
and Geometric Quantization  \publ 
 Birkhauser\publaddr Boston \yr 1993 \endref

\ref\key [C] \by J. Cheeger\paper Multiplication of 
Differential
Characters \jour Instituto Nazionale di Alta Mathematica, 
Symposia Mathematica \vol XI \yr 1973 \pages 441--445
\endref

\ref\key [CS]\by J. Cheeger and J. Simons\paper
Differential Characters and Geometric Invariants
\inbook Geometry and Topology \publ Lect. Notes in Math. no. 1167,
Springer--Verlag \publaddr New York\yr 1985 \pages 50--80\endref

\ref\key [Ch]\by S.S. Chern \paper
A simple intrinsic proof of the Gauss-Bonnet formula for closed
Riemannian manifolds\jour Ann. Math. \vol 45 \year 1944
\pages 747-752\endref

\ref\key [Cha]\by D. S. Chatterjee \paper
On the construction of abelian gerbs\jour Ph. D. Thesis, Cambridge 
\year 1998 \endref

\ref\key [deR]\by  G. de Rham \book Vari\'et\'es Diff\'erentiables, formes,
courants, formes harmoniques\publ 
 Hermann\publaddr Paris \yr 1955\endref

\ref\key [F]\by   H. Federer\book Geometric Measure 
Theory\publ 
 Springer--Verlag\publaddr New York \yr 1969\endref

\ref\key [FW] \by D. Freed and E. Witten\paper Anomalies in string
theory with D-branes \jour Asian J. Math. \vol 3 \yr 1999  \pages 819--
851\endref

\ref \key [GS$_1$] \by H. Gillet and C. Soul\'e \paper  
Arithmetic chow groups and differential characters
 \pages 30-68 \inbook Algebraic K-theory;  Connections with Geometry and
Topology   \publ Jardine and Snaith (eds.), Kluwer Academic Publishers \yr
1989 \endref

\ref \key [GS$_2$] \bysame \paper  
Arithmetic intersection theory
 \pages 94-174 \jour Publ. I.H.E.S.  \vol 72 \yr 1990 \endref

\ref \key [GS$_3$] \bysame \paper  
Characteristic classes for algebraic vector bundles with
hermitian metrics, I and II
 \pages 163-203 and 205-238 \jour Ann. of Math.  \vol 131 \yr 1990 \endref

\ref \key [GS$_4$] \bysame \paper  
An arithmetic Riemann-Roch Theorem
 \pages 473-543  \jour Inventiones Math.  \vol 110 \yr 1992\endref

\ref \key [GS$_5$] \bysame \paper 
Arithmetic analogues of the standard conjectures
 \pages 129-140 \inbook Motives (Seattle, WA, 1991), Proc. of Symp. in
Pure Math. 55, part I    \publ Amer. Math. Soc. \publaddr Providence,
RI\yr 1994 \endref

\ref \key [H] \by B. Harris \paper  
Differential characters and the Abel-Jacobi map
 \pages 69-86 \inbook Algebraic K-theory;  Connections with Geometry and
Topology   \publ Jardine and Snaith (eds.), Kluwer Academic Publishers \yr
1989 \endref

\ref \key [HL$_1$] \by F.R. Harvey and H.B. Lawson, Jr. \book A theory of
 characteristic currents associated with a singular connection
 \yr 1993 \publ Ast\'erisque, 213
 Soci\'et\'e   Math. de France\publaddr Paris\endref

\ref\key [HL$_2$] \bysame \paper Geometric residue
 theorems\jour Amer. J. Math. \vol 117 \yr 1995  \pages 829--
873\endref

\ref\key [HL$_3$] \bysame \paper 
Singularities and Chern-Weil theory, I --
The local MacPherson Formula 
\jour Asian J. Math.  \vol 4 No 1 \yr 2000 \pages  71-96
\endref

\ref\key [HL$_4$] \bysame \paper 
Singularities and Chern-Weil theory, II --
geometric atomicity
\jour Duke Math. Journal (to appear)
\endref

\ref\key [HL$_5$] \bysame \paper 
Finite volume flows and Morse theory\jour  Annals of Math. \vol 153 No 1
\year 2001 
\pages 1-25 \endref

\ref\key [HL$_6$] \bysame \paper Lefschetz-Pontrjagin duality for
differential characters \jour Anais da Academia Brasileira de Ci\^encias
\vol 73 \yr 2001 \pages 145--159  \endref

\ref\key [HL$_7$] \bysame \paper 
From sparks to grundles -- differential characters
\jour  (to appear) \pages   \endref

\ref\key [HP]  \by F.R. Harvey and J. Polking  \paper  Fundamental solutions
in complex analysis, I and II  \jour  Duke Math. J.  \vol 46\yr  1979 \pages
253--340
\endref

\ref\key [HS]  \by F.R. Harvey and S. Semmes
\paper  Zero divisors of atomic functions
\jour  Ann. of Math.\vol 135 \yr 1992 \pages 567--600
\endref

\ref\key [HZ$_1$] \by F.R. Harvey and J. Zweck\paper Stiefel--Whitney
Currents\jour J. Geometric Analysis
\vol 8 No 5 \yr 1998 \pages  805--840
\endref

\ref\key  [HZ$_2$]  \bysame \paper Divisors and Euler sparks of atomic 
sections
\jour Indiana Univ. Math. Jour. \vol 50 \yr 2001 \pages 243--298
\endref

\ref\key [Hi] \by N. Hitchin \paper Lectures on special lagrangian
submanifolds

\jour arXiv:math.DG/9907034 \vol  \yr 1999
   \pages  \endref

\ref\key [R] \by F. Ronga \paper  Le calcul de la classe de 
cohomologie duale a $\overline{\Sigma}^k$ -- Proceedings of 
Liverpool
 Singularities  Symposium, I\jour Lecture Notes in Mathematics
\vol 192\pages 313--315\yr 1971\endref

\ref \key [S] \by J. Simons \paper Characteristic forms and
transgression: characters associated to a connection \jour Stony
Brook preprint,   \yr 1974 \endref

 \enddocument